\documentclass{amsart}

\usepackage{epsfig}
\usepackage[all]{xy}

\newtheorem{thm}{Theorem}[section]
\newtheorem{lemma}[thm]{Lemma}
\newtheorem{prop}[thm]{Proposition}
\newtheorem{cor}[thm]{Corollary}
\newtheorem*{thmintro}{Theorem}
\newtheorem*{propintro}{Proposition}

\theoremstyle{definition}
\newtheorem*{definition}{Definition}
\newtheorem{ex}[thm]{Example}

\newcounter{Tlistc}
\newenvironment{Tlist}
    	{\setcounter{Tlistc}{0}
	 \begin{list}{(T\arabic{Tlistc})}
	{\usecounter{Tlistc}}}{\end{list}}
	 
\newenvironment{romanlist}
	{\begin{enumerate}
	}
	{\end{enumerate}}

\newcommand{\Zz}{\mathbb Z}
\newcommand{\Rr}{\mathbb R}
\newcommand{\Cc}{\mathbb C}
\newcommand{\Qq}{\mathbb Q}
\newcommand{\pp}{\mathcal{P}}
\newcommand{\eps}{\varepsilon}
\newcommand{\rk}{\hbox{rk}\,}
\newcommand{\sign}{\hbox{sign}\,}
\newcommand{\nuli}{\hbox{null}\,} 
 
\begin{document}

\title{Generalized Seifert surfaces and signatures of colored links}

\author{David Cimasoni}
\address{Section de Math\'ematiques, Universit\'e de Gen\`eve, 2--4 rue du Li\`evre, 1211 Gen\`eve 24, Switzerland}
\email{David.Cimasoni@math.unige.ch}
\thanks{The first author was supported by the Swiss National Science Foundation.}

\author{Vincent Florens}
\address{Universidad de Valladolid, Dto \~Algebra, Geometr\~a y Topolog\~a, Prado de la Magdalena s/n, 47011 Valladolid, Spain}
\email{vincent$\_$florens@yahoo.fr}
\thanks{The second author was supported by Marie-Curie, MCHF-2001-0615.}


\date{\today}


\begin{abstract}
In this paper, we use `generalized Seifert surfaces' to extend the Levine-Tristram signature to colored links in $S^3$. This yields an
integral valued function on the $\mu$-dimensional torus, where $\mu$ is the number of colors of the link. The case $\mu=1$
corresponds to the Levine-Tristram signature. We show that many remarkable properties of the latter invariant extend to this
$\mu$-variable generalization: it vanishes for achiral colored links, it is `piecewise continuous', and the places of the jumps
are determined by the Alexander invariants of the colored link. Using a $4$-dimensional interpretation and the Atiyah-Singer
$G$-signature theorem, we also prove that this signature is invariant by colored concordance, and that it provides a lower bound
for the `slice genus' of the colored link.
\end{abstract}

\maketitle
 
\section{Introduction}\label{section:intro}

Several notions related to knots do not extend naturally and uniquely to links. For example, the fact that two oriented links are
isotopic can be understood in different ways: one might require the isotopy to satisfy some condition, e.g. to respect an order
on the components of the links. Here is another interesting example. A knot in $S^3$ is said to be \emph{slice\/}
if it bounds a smooth disc in the $4$-ball, or equivalently, if it is the cross-section of a smooth $2$-sphere in $S^4$. This notion of
sliceness for knots admits different generalizations to links. According to Fox \cite{Fo2}, a link is \emph{slice in the ordinary
sense\/} if it is the cross-section of a single smooth $2$-sphere in $S^4$. It is \emph{slice in the strong sense\/} if each of its
components is such a cross-section for disjoint $2$-spheres in $S^4$.

One way to take into account simultaneously this variety of possible generalizations is to consider so-called \emph{colored links}.
Roughly speaking, a $\mu$-colored link is an oriented link in $S^3$ whose components are endowed with some integer in $\{1,\dots,\mu\}$
called the \emph{color\/} of the component. Two colored links are isotopic if there is an isotopy between them which respects the color
and orientation of each
component. We shall say that a $\mu$-colored link is slice if there exists $\mu$ disjoint smooth spheres $S_1,\dots,S_\mu$ in $S^4$
such that the sublink of color $i$ is a cross-section of $S_i$. Of course, a $1$-colored link is nothing but an oriented link, and
it is slice as a $1$-colored link if it is slice in the ordinary sense. At the other end of the spectrum, a $\nu$-component $\nu$-colored
link is an ordered link. It is slice as a $\nu$-colored link if it is slice in the strong sense.

Many classical invariants of oriented links, such as the Alexander polynomial and the Levine-Tristram signature, can be constructed
using Seifert surfaces. In this paper, we introduce generalized Seifert surfaces for colored links. We use them, \emph{inter alia\/}, to
extend the Levine-Tristram signature from oriented links to $\mu$-colored links. This yields a integral valued function on the
$\mu$-dimensional torus. Among other results (see the paragraph below), we show that this function vanishes almost everywhere
if the $\mu$-colored link is slice.

\medskip

Throughout the paper, all the links are assumed to be smooth and oriented.   

\subsection*{The Levine-Tristram signature.}

Let $V$ be a Seifert matrix for a link $L$ in $S^3$.
Then, $A(t)=V-tV^T$ is a presentation matrix of the Alexander $\Zz[t^{\pm 1}]$-module of $L$. In particular,
the Alexander polynomial $\Delta_L$ of $L$ is given by the determinant of $A(t)$. If $\omega\neq 1$ is a unit modulus complex
number, then $H(\omega)=(1-\overline{\omega})A(\omega)$ is a Hermitian matrix whose signature $\sigma_L(\omega)$ and nullity
$\eta_L(\omega)$ do not depend on the choice of $V$. This yields integral valued functions $\sigma_L$ and $\eta_L$ defined on
$S^1\setminus\{1\}$. In the case of $\omega=-1$, this signature was first defined by Trotter \cite{Tro} and studied by
Murasugi \cite{Mu}. The more general formulation is due to Levine \cite{Le1} and Tristram \cite{Tris}, and referred to as the
{\em Levine-Tristram signature\/}.

The functions $\sigma_L$ and $\eta_L$ are easily seen to be locally constant on the complement in
$S^1\setminus\{1\}$ of the roots of $\Delta_L$. Also, $\eta_L$ is related to the first Betti number of the finite cyclic coverings
of the exterior of $L$. 
Moreover, when restricted to roots of unity of prime order, the signature and nullity are concordance invariants.
(The case of $\omega=-1$ is due to Murasugi, and Tristram extended it to any $\omega$ of prime order.)
Finally, the so-called \emph{Murasugi-Tristram inequality\/} imposes a condition, expressed in terms of the values of $\sigma_L$
and $\eta_L$, on the Betti numbers of a smooth oriented surface in $B^4$ spanning $L$. This inequality implies in particular
that if $L$ is {\em slice in the strong sense\/} \cite{Fo2}, then $\sigma_L$ vanishes at roots of unity of prime order.

At that point, all the methods of demonstration were purely $3$-dimensional. A new light was shed on this theory in the early seventies.
Building on ideas of Rokhlin \cite{Ro}, Viro \cite{Vi} was able to interpret the Levine-Tristram signature as a $4$-dimensional object.
Indeed, he showed that for all rational values of $\omega$, $\sigma_L(\omega)$ coincides with the signature of an
intersection form related to a cyclic cover of $B^4$ branched along a Seifert surface for $L$ pushed in the interior of $B^4$.
This $4$-dimensional approach was used by Kauffman and Taylor \cite{Kauff-Tay} to obtain a short proof of the Murasugi-Tristram
inequality, in the case $\omega=-1$. They were also able to show the following inequality: if $P$ is a closed oriented smooth
surface in $S^4$ that intersects the standardly embedded $3$-sphere in $L$, then
\[
|\sigma_L(-1)|\le genus(P) + \min(0,\eta_L(-1)-\beta_0(P)+1),\eqno{(\star)}
\]
where $\beta_0(P)$ denotes the number of connected components of $P$.
In particular, if there exists such a surface of genus $0$ (that is, according to \cite{Fo2}, if $L$ is slice in the ordinary sense),
then $\sigma_L(-1)=0$. This $4$-dimensional interpretation was used with great success by several authors \cite{Gi1,G-L-M,Le2,Le3,Smo}.
See also \cite{Cass-Gord1,Flo-Gi,Gi2,Gil-Liv}.

\subsection*{Paper outline and statement of the results.}

The aim of this paper is to generalize the Levine-Tristram signature
to colored links. A {\em $\mu$-colored link\/} $L=L_1\cup\dots\cup L_\mu$ is an oriented link in $S^3$ together with  a
surjective map assigning to each component of $L$ a color in $\{1,\dots,\mu\}$. The sublink $L_i$ is constitued by
the components of $L$ with color $i$ for $i=1,\dots,\mu$.
By isotopy of colored links, we mean orientation and color-preserving isotopy. Note that a $1$-colored link
is an ordinary link, and setting $\mu=1$ in this article gives back to the known results.

\smallskip
In Section \ref{section:def}, we consider generalized Seifert surfaces called \emph{C-complexes\/}. Roughly speaking,
a C-complex for a $\mu$-colored link $L$ consists of a collection of Seifert surfaces $S_1,\dots,S_\mu$ for the sublinks
$L_1,\dots,L_\mu$ that 
intersect only along clasps. Associated to a C-complex are so-called generalized Seifert matrices. We use them to define a matrix
$A(t_1,\dots,t_\mu)$ with coefficients in $\Lambda_\mu=\Zz[t_1^{\pm 1},\dots,t_\mu^{\pm 1}]$. Of course, this matrix depends on the
choice of the C-complex for $L$. However, if $(\omega_1,\dots,\omega_\mu)$ is an element of the $\mu$-dimensional torus
$S^1\times\dots\times S^1\subset\Cc^\mu$ with $\omega_i \neq 1$, then the matrix
\[
H(\omega_1,\dots,\omega_\mu)=\prod_{i=1}^\mu (1-\overline{\omega}_i) \cdot A(\omega_1,\dots,\omega_\mu)
\]
is Hermitian and its signature and nullity are independant of the choice of the C-complex for $L$. (We use some
`generalized S-equivalence', see Lemma \ref{lemma:Sequiv} and Theorem \ref{thm:inv}.) This allows one to define 
the {\em signature\/} and the {\em nullity\/} of the $\mu$-colored link $L$ as functions
\[
\sigma_L,\eta_L\colon(S^1\setminus\{1\})^\mu\longrightarrow\Zz.
\]
Note that this was done by Cooper \cite{Coop,Coo} in the case of a $2$-colored link with $2$ components.
The restriction of these functions to the diagonal specializes to the Tristram-Levine invariants as follows
(see Proposition \ref{prop:lk}).

\begin{propintro}
Let $L=L_1 \cup \dots \cup L_\mu$ be a $\mu$-colored link, and let $L^\prime$ be the underlying link. For all $\omega$ in
$S^1\setminus\{1\}$,
\[
\sigma_{L}(\omega,\dots,\omega) = \sigma_{L^\prime}(\omega) + \sum_{i<j}lk(L_i,L_j)\quad\hbox{and}\quad
\eta_L(\omega,\dots,\omega) = \eta_{L^\prime}(\omega),
\]
where $lk$ denotes the linking number in $S^3$.
\end{propintro}
\noindent This result can be quite useful. Indeed, it is often easier to compute a multivariable signature (corresponding to a well-chosen
coloring), and to evaluate it on the diagonal, then to compute directly the Levine-Tristram invariant.

Several interesting properties of the Levine-Tristram invariants extend to our functions $\sigma_L$ and $\eta_L$.
For example, they are additive under disjoint and connected sum. Moreover, $\sigma_L$ vanishes if the colored link $L$
is isotopic to its mirror image (see Corollary \ref{cor:mirror}).

\smallskip
Section \ref{section:Alex} is devoted to the study of the natural $\Zz^\mu$-covering $\widetilde X\to X$ of the exterior of $L$ induced by the coloring of
$L$. We show that C-complexes provide a nice geometrical description of this covering, and deduce a presentation of the Alexander
$\Lambda_\mu$-module $H_1(\widetilde X)$ in terms of the generalized Seifert forms (see Theorem \ref{thm:module}). In particular,
if $\Lambda_\mu^\prime$ denotes the localization of the ring $\Lambda_\mu$ with respect to the multiplicative system generated by $t_i-1$  
for $i=1,\dots,\mu$, we establish the following result (Corollary \ref{cor:polynomial}):
\begin{thmintro}
Let $L$ be a $\mu$-colored link. Consider a C-complex $S$ for $L$ such that $S_i$ is connected for all $i$ and
$S_i\cap S_j$ is non-empty for all $i\neq j$. Then the corresponding matrix
$A(t_1,\dots,t_\mu)$ is a presentation matrix of the $\Lambda_\mu^\prime$-module
$H_1(\widetilde X)\otimes_{\Lambda_\mu}\Lambda_\mu^\prime$.
\end{thmintro}
\noindent In particular, the Alexander polynomial of $L$ is equal to the determinant of the matrix 
$A(t_1,\dots,t_\mu)$ up to multiplication by
$t_i-1$. (This latter result was obtained by the first author in \cite{Cim} using local relations on a diagram.)
This theorem implies the following characterization of the discontinuities
of $\sigma_L$ and $\eta_L$ (Theorem \ref{thm:semi}).
 
\begin{thmintro}
Let $E_r(L)\subset\Lambda_\mu$ be the $r^{th}$ Alexander ideal of $L$, and set
\[
\Sigma_r=\{\omega\in S^1 \times \dots \times S^1 \subset \Cc^\mu \;\mid\;p(\omega)=0 \hbox{ for all } p\in E_{r-1}(L)\}. 
\]
This yields a finite sequence of algebraic subvarieties of the torus $S^1 \times \dots \times S^1$. Setting
$\Sigma_r^*= \Sigma_r \cap (S^1\setminus\{1 \})^\mu$, we obtain a finite sequence
$(S^1\setminus\{1\})^\mu=\Sigma_0^*\supset\Sigma_1^*\supset\dots\supset\Sigma_{\ell-1}^*\supset\Sigma_\ell^*=\emptyset$
such that, for all $r$, $\eta_L$ is equal to $r$ on $\Sigma_r^*\setminus\Sigma_{r+1}^*$, and $\sigma_L$ is locally constant on
$\Sigma_r^*\setminus\Sigma_{r+1}^*$.
\end{thmintro} 
\noindent This `piecewise continuity' behavior was first observed by Levine \cite{Le3} for closely related invariants.
The most interesting point of our result is the relation to the Alexander invariants. Note that
even if the Alexander polynomial is zero, the signature is locally constant. 

\smallskip
In Section \ref{section:local}, we build on an idea of Conway \cite{Con} to show that the signature satisfies several `local relations.'
(We refer to Theorem \ref{thm:local} for a precise statement.) In many cases, this leads to a purely combinatorial
computation of the signature from a diagram of the corresponding colored link.

\smallskip
The following section deals with a $4$-dimensional interpretation of $\sigma_L(\omega)$ and $\eta_L(\omega)$ for all
$\omega=(\omega_1,\dots,\omega_\mu)$ with rational coordinates in $(S^1\setminus\{1\})^\mu$. We consider a union $F$ of connected
surfaces $F_1,\dots,F_\mu$ smoothly embedded in $B^4$ such that $\partial F_i\subset\partial B^4=S^3$ is the sublink $L_i$, and the
pairwise intersections of the $F_i$'s are transverse (along a finite number of points). The first homology of the exterior $W_F$ in $B^4$
of such a `spanning
surface' is free of rank $\mu$. Therefore, any rational point $\omega \in (S^1\setminus\{1\})^\mu$ determines a character of $H_1(W_F)$
of finite order. This character induces twisted homology $\Cc$-vector spaces, denoted by $H_*^\omega(W_F;\Cc)$, and a Hermitian twisted
intersection form $\varphi^\omega_F\colon H_2^\omega(W_F;\Cc)\times H_2^\omega(W_F;\Cc)\to\Cc$.
We obtain the following result (see Lemma \ref{decomp}, Lemma \ref{ram} and Proposition \ref{ccompl}).
   
\begin{propintro} 
Consider a connected C-complex $S\subset S^3$ for a $\mu$-colored link $L$. 
Let $F\subset B^4$ be the spanning surface for $L$ obtained by pushing $S$ in $B^4$.
For any rational point $\omega\in (S^1\setminus\{1\})^\mu$, $H(\omega)$ is a matrix for $\varphi^\omega_F$.
\end{propintro} 
\noindent The proof follows from an explicit geometrical description of the finite abelian coverings of $W_F$.
We also make use of the work of Sakuma \cite{Sak} for the study of cyclic quotients of these coverings.
The Atiyah-Singer $G$-signature theorem \cite{Aty-Sin} implies that the signature of $\varphi^\omega_F$ does not depend on the choice of
the spanning surface $F$ for $L$. Moreover, the nullity of $\varphi^\omega_F$ is closely related to the twisted homology
of the exterior $X$ of $L$ in $S^3$. This leads to the following result (Theorem \ref{main1}).

\begin{thmintro}
Let $L$ be a $\mu$-colored link with exterior $X$, and $F$ be a spanning surface for $L$ in $B^4$. For all rational points $\omega$ in
$(S^1\setminus\{1\})^\mu$,
\begin{eqnarray*}
\sigma_L(\omega)&=&\sign(\varphi_F^\omega),\\
\eta_L(\omega)&=&\dim H_1^\omega(X;\Cc)=\nuli(\varphi_F^\omega)+\dim H_1^\omega(W_F;\Cc)-\dim H_3^\omega(W_F;\Cc).
\end{eqnarray*}
\end{thmintro}    
\noindent As a consequence, if the colored link $L=L_1\cup\dots\cup L_\mu$ satisfies $lk(L_i,L_j)=0$ for all $i\neq j$,
then $\sigma_L(\omega)$ and $\eta_L(\omega)$ coincide with the invariants considered by the second author in \cite{Flo}.
It also relates our signature function to various invariants introduced by Gilmer
\cite{Gi1}, Smolinsky \cite{Smo} and Levine \cite{Le3}. Note that the first equality concerning the nullity is closely
related to Libgober's \cite{Li}.
 
Combining our construction with \cite{Gi1}, we obtain the following formula for the Casson-Gordon invariant of a $3$-manifold
(Theorem \ref{as}).
\begin{thmintro}
Let $M$ be the $3$-manifold obtained by surgery on a framed link $L$ with $\nu$ components and linking matrix $\Lambda$.
Let $\chi\colon H_1(M)\to\Cc^*$ be the character mapping the meridian of the $i^{th}$ component of $L$ to $\alpha^{n_i}$,
where $\alpha=e^{2i\pi/q}$ and $n_i$ is an integer coprime to $q$. Consider $L$ as a $\nu$-colored link and set
$\omega=(\alpha^{n_1},\dots,\alpha^{n_\nu})$. Then, the Casson-Gordon invariant of the pair $(M,\chi) $ is given by
\[
\sigma(M,\chi)=\Big(\sigma_L(\omega)-\sum_{i<j}\Lambda_{ij}\Big)-\sign(\Lambda)+\frac{2}{q^2}\sum_{i,j}(q-n_i)n_j\Lambda_{ij}.
\]
\end{thmintro}

The $4$-dimensional point of view developped in Section \ref{section:4-dim} makes it possible to prove several results that would
have been horrendous to check using only $3$-dimensional techniques. Let us denote by $T^\mu_\pp$ the dense subset of
$S^1\times\dots\times S^1$ given by the elements of the form $\omega=(\omega_1,\dots,\omega_\mu)$ which satisfy the following condition:
there exists a prime $p$ such that for all $i$, the order of $\omega_i$ is a power of $p$.

\begin{thmintro}
For all $\omega \in T^\mu_\pp$, $\sigma_L(\omega)$ and $\eta_L(\omega)$ are invariant by colored concordance.
\end{thmintro}
We also extend the Murasugi-Tristram inequality to the case of surfaces that intersect transversally (Theorem \ref{thm:MT}).
This can be viewed as a specialization of \cite[Theorem 4.1]{Gi1}.
Finally, we show the following generalization of the Kauffman-Taylor inequality $(\star)$.

\begin{thmintro}
Consider a colored link $L=L_1\cup\dots\cup L_\mu$. Let us assume that there exists a smooth oriented closed surface
$P=P_1\sqcup\dots\sqcup P_\mu$ in $S^4$ such that $P_i\cap S^3=L_i$ for all $i$, where $S^3$ denotes the standard embedding of the
$3$-sphere in $S^4$. Then, for all $\omega$ in $T^\mu_\pp$,
\[
|\sigma_L(\omega)| \leq genus(P) + \min(0,\eta_L(\omega)-\mu+1).
\]
\end{thmintro}
If there exists such a surface $P$ of genus zero, we say that $L$ is a {\em slice colored link\/}. As an immediate corollary, we get:
if $L$ is a slice $\mu$-colored link, then $\sigma_L(\omega)=0$ and $\eta_L(\omega)\ge\mu-1$ for all $\omega$ in $T^\mu_\pp$.
This notion of `sliceness' is in fact a natural generalization of the definitions of Fox \cite{Fo2} stated above.
Indeed, a $1$-colored link is slice if and only if it is slice in the ordinary sense. On the other hand, a $\nu$-component link is
slice as a $\nu$-colored link if and only if it is slice in the strong sense. What we get is a spectrum of sliceness notions ranging 
from the ordinary sense to the strong sense. To each coloring of a given link, there corresponds one notion of sliceness,
and one signature function. This function vanishes if the link is slice in the corresponding sense.

\smallskip
It should be pointed out that all the results of Sections \ref{section:def} to \ref{section:PC} hold for colored links
in an arbitrary $\Zz$-homology $3$-sphere. Sections \ref{section:4-dim} and \ref{section:slice} also extend to this setting,
provided the homology sphere bounds a contractible $4$-manifold.

\smallskip
Finally, let us mention that the results of the present paper have been 
successfully applied to the study of the topology of real algebraic plane curves. 
Indeed, S. Yu. Orevkov implemented an algorithm computing the 
generalized Seifert matrices (and therefore, the signature and nullity 
functions) of a colored link given as the closure of a colored braid. Using this 
computer program, the generalized Murasugi-Tristram inequality (Theorem \ref{thm:MT} below), and his method 
developed in \cite{Ore1}, he was able to complete the classification up to isotopy of 
M-curves of degree 9 with 4 nests. We refer to the upcoming paper \cite{Ore3} for details.

\subsection*{Acknowledgments.}

A part of this paper was done while the first author visited the Institut de
Recherche Math\'ematique Avanc\'ee (Strasbourg) whose hospitality he thankfully acknowledges.
He also wishes to thank Vladimir Turaev and Mathieu Baillif for valuable discussions.
The second author thanks the Universit\'e Libre de Bruxelles and the Department of Geometry and Topology of the University of Zaragoza 
for their hospitality. He thanks in particular E. B. Artal and  J. I. Cogolludo for exciting conversations.

\section{Definition and basic properties of $\sigma_L$ and $\eta_L$}\label{section:def}

The aim of this section is to define the signature and nullity of a colored
link as a natural generalization of Levine-Tristram signature of an oriented link.
 
\subsection{C-complexes}\label{sub:C-cplx}

Recall that a Seifert surface for a link in $S^3$ is a connected compact oriented surface smoothly
embedded in $S^3$ that has the link as its oriented boundary. The notion of C-complex, as introduced
in \cite{Coo} and \cite{Cim}, is a generalization of Seifert surfaces to colored links.

\begin{figure}[b]
\begin{center}
\epsfig{figure=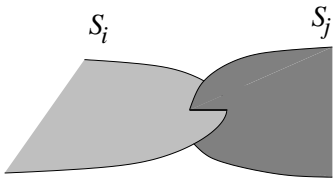,height=2.3cm}
\caption{A clasp intersection.}
\label{fig:clasp}
\end{center}
\end{figure}

\begin{definition}
A {\em C-complex\/} for a $\mu$-colored link $L=L_1\cup\dots\cup L_\mu$ is a union
$S=S_1\cup\dots\cup S_\mu$ of surfaces in $S^3$ such that:
\begin{romanlist}
\item{for all $i$, $S_i$ is a  Seifert surface for $L_i$ (possibly disconnected, but with no closed components);}
\item{for all $i\neq j$, $S_i\cap S_j$ is either empty or a union of clasps (see Figure \ref{fig:clasp});}
\item{for all $i,j,k$ pairwise distinct, $S_i\cap S_j\cap S_k$ is empty.}
\end{romanlist}
\end{definition}

The existence of a C-complex for a colored link is fairly easy, see 
\cite[Lemma 1]{Cim}. In the case $\mu=1$, a C-complex for $L$ is nothing but a (possibly disconnected)
Seifert surface for the link $L$. Let us now define the corresponding generalization of the Seifert form.
Let $N_i=S_i\times[-1,1]$ be a bicollar neighborhood of $S_i$. Given a sign $\eps_i=\pm 1$, let $S_i^{\eps_i}$ be the translated
surface $S_i\times\{\eps_i\}\subset N_i$. Also, let $X$ be the complement of an open tubular neighborhood of $L$, and let
$Y$ be the complement of $\bigcup_{i=1}^\mu int\,N_i$ in $X$. Given a sequence
$\eps=(\eps_1,\dots,\eps_\mu)$ of $\pm 1$'s, set
\[
S^\eps=\bigcup_{i=1}^\mu S_i^{\eps_i}\cap Y.
\]
Since all the intersections are clasps, there is an obvious homotopy equivalence
between $S$ and $S^\eps$ inducing an isomorphism $H_1(S)\to H_1(S^\eps)$. Let $i^\eps\colon H_1(S)\to H_1(S^3\setminus S)$ be
the composition of this isomorphism with the inclusion homomorphism $H_1(S^\eps)\to H_1(S^3\setminus S)$. Finally, let
\[
\alpha^\eps\colon H_1(S)\times H_1(S)\to\Zz
\]
be the bilinear form given by $\alpha^\eps(x,y)=lk(i^\eps(x),y)$,
where $lk$ denotes the linking number. Fix a basis of $H_1(S)$ and denote by $A_S^\eps$ (or
simply by $A^\eps$) the matrix of $\alpha^\eps$. Of course, if $\mu=1$, then $\alpha^-$ is the usual Seifert form and $A^-$ the
usual Seifert matrix. Note that for all $\eps$, $A^{-\eps}$ is equal to $(A^\eps)^T$, the transpose of the matrix $A^\eps$.  

\begin{figure}[Htb]
\begin{center}
\epsfig{figure=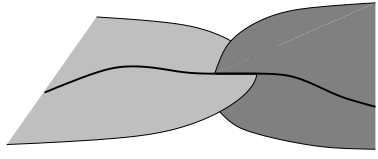,height=2cm}
\caption{A loop crossing a clasp.}
\label{fig:loop}
\end{center}
\end{figure}

For computational purposes, the following alternative definition of $i^\eps$ is more convenient.
A $1$-cycle in a C-complex is called a {\em loop\/} if it is an oriented simple closed curve which behaves as
illustrated in Figure \ref{fig:loop} whenever it crosses a clasp.
Clearly, there exists a collection of loops whose homology classes
form a basis of $H_1(S)$. Therefore, it is possible to define $i^\eps$ as follows: for any loop $x$, $i^\eps([x])$ is
the class of the $1$-cycle obtained by pushing $x$ in the $\eps_i$-normal direction off $S_i$ for $i=1,\dots,\mu$.
The fact that $x$ is a loop ensures that this can be done continuously along the clasp intersections.
We easily check that this definition of $i^\eps$ coincides with the intrinsic definition given above.

\subsection{The signature and nullity of a colored link}\label{sub:sign}

Let $L$ be a $\mu$-colored link. Consider a C-complex $S$ for $L$ and the associated Seifert matrices $A^\eps$
with respect to some fixed basis of $H_1(S)$.
Let $A(t_1,\dots,t_\mu)$ be the matrix with coefficients in $\Zz[t_1,\dots,t_\mu]$ defined by
\[
A(t_1,\dots,t_\mu)=\sum_\eps\,\eps_1\cdots\eps_\mu\,t_1^{\frac{1-\eps_1}{2}}\cdots t_\mu^{\frac{1-\eps_\mu}{2}}A^\eps,
\]
where the sum is on the $2^\mu$ possible sequences $\eps=(\eps_1,\dots,\eps_\mu)$ of $\pm 1$'s.
For $\omega=(\omega_1,\dots,\omega_\mu)$ in $T^\mu= S^1\times\dots\times S^1\subset\Cc^\mu$, set
\[
H(\omega)\;=\;\prod_{i=1}^\mu(1-\overline\omega_i)\,A(\omega_1,\dots,\omega_\mu).
\]
Using the fact that $A^{-\eps}=(A^\eps)^T$, one easily checks that $H(\omega)$ is a Hermitian matrix. Recall that the eigenvalues
of such a matrix $H$ are real. Its signature $\sign(H)$ is defined as the number of positive eigenvalues minus the number
of negative eigenvalues. The nullity $\nuli(H)$ is the number of zero eigenvalues of $H$.

\begin{definition}
Let $T^\mu_\ast$ be the open subset $(S^1\setminus\{1\})^\mu$ of the $\mu$-dimensional torus $T^\mu\subset\Cc^\mu$.
The {\em signature} and {\em nullity} of the $\mu$-colored link $L$ are the functions
\[
\sigma_L, \eta_L\colon T^\mu_\ast\longrightarrow\Zz
\]
given by $\sigma_L(\omega)=\sign(H(\omega))$ and $\eta_L(\omega)=\nuli(H(\omega))+\beta_0(S)-1$, where $\beta_0(S)$ denotes the
number of connected components of $S$.
\end{definition}

By Sylvester's theorem, $\sigma_L(\omega)$ and $\eta_L(\omega)$ do not depend on
the choice of a basis of $H_1(S)$.

\begin{thm}\label{thm:inv}
The signature $\sigma_L$ and nullity $\eta_L$ do not depend on the choice of the C-complex
for the colored link $L$. Hence, they are well-defined as isotopy invariants of the colored link $L$.
\end{thm}

This theorem relies on the following lemma (see \cite{Cim} for the proof).

\begin{figure}[Htb]
\begin{center}
\epsfig{figure=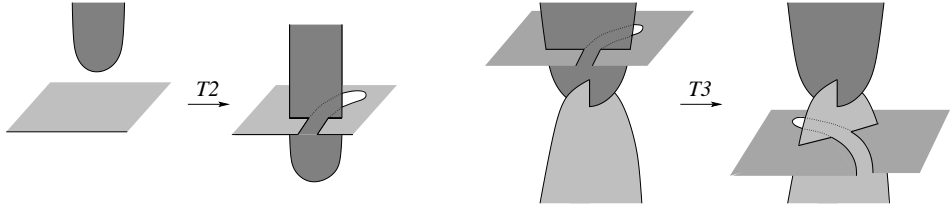,height=2.7cm}
\caption{The transformations $T2$ and $T3$ in Lemma \ref{lemma:Sequiv}.}
\label{fig:Sequiv}
\end{center}
\end{figure}

\begin{lemma}\label{lemma:Sequiv}
Let $S$ and $S^\prime$ be C-complexes for isotopic colored links.
Then, $S$ and $S'$ can be transformed into each other by a finite number of the
following operations and their inverses:
\begin{Tlist}
\setcounter{Tlistc}{-1}
\item{ambient isotopy;}
\item{handle attachment on one surface;}
\item{addition of a ribbon intersection, followed by a `push along an arc'
through this intersection (see Figure \ref{fig:Sequiv});}
\item{the transformation described in Figure \ref{fig:Sequiv}.}
\end{Tlist}
\end{lemma}

\begin{proof}[Proof of Theorem \ref{thm:inv}]
Let $H$ and $H'$ be two Hermitian matrices. We shall call $H'$ an elementary
enlargement of $H$ if
\[
H'=\begin{pmatrix}
H&\xi&0\cr \overline\xi^T&\lambda&\alpha\cr 0&\overline\alpha&0
\end{pmatrix},
\]
where $\xi$ is any complex vector, $\lambda$ any real number and
$\alpha\in\Cc^\ast$. $H$ is called an elementary reduction of $H'$.
One easily checks that the signature and nullity of a Hermitian matrix are
unchanged by elementary enlargements and reductions.
By Lemma \ref{lemma:Sequiv}, it remains to prove that the transformations $T1$
to $T3$ of a C-complex induce finite sequences
of elementary reductions and enlargements on the corresponding Hermitian
matrices (or other elementary transformations leaving the signature and nullity unchanged).

\smallskip\noindent{$(T1)$} Let $S'$ be a C-complex obtained from a C-complex
$S$ by a handle attachment on $S_k$.  If this handle connects two distinct connected components of $S$,
then $H_1(S')=H_1(S)\oplus\Zz y$, where $y$ in a $1$-cycle such that $lk(i^\eps(y),z)=0$ for all $z$ in $H_1(S')$.
Hence, the Hermitian matrices $H$ and $H'$ corresponding to $S$ and $S'$ are related by $H'=H\oplus(0)$,
so $\sign(H')=\sign(H)$ and $\nuli(H')=\nuli(H)+1$. Since $\beta_0(S')=\beta_0(S)-1$, the signature
and nullity of $L$ are unchanged.
Let us now assume that this handle attachement is performed on one connected component of $S$. In this case,
$H_1(S')=H_1(S)\oplus\Zz x\oplus\Zz y$. Moreover, the cycles $x$ and $y$ can be chosen so that the corresponding
Seifert matrices satisfy
\[
A_{S'}^\eps=\begin{pmatrix}
A_S^\eps&\ast&0\cr \ast&\ast&\pi(-\eps)\cr
0&\pi(\eps)&0
\end{pmatrix}\,,
\quad\hbox{with}\quad\pi(\eps)=\begin{cases}1&\text{if $\eps_k=+1$;}\cr 0&\text{else.}\end{cases}
\]
Such a choice of $x$ and $y$ is illustrated below.

\begin{figure}[h]
\begin{center}
\epsfig{figure=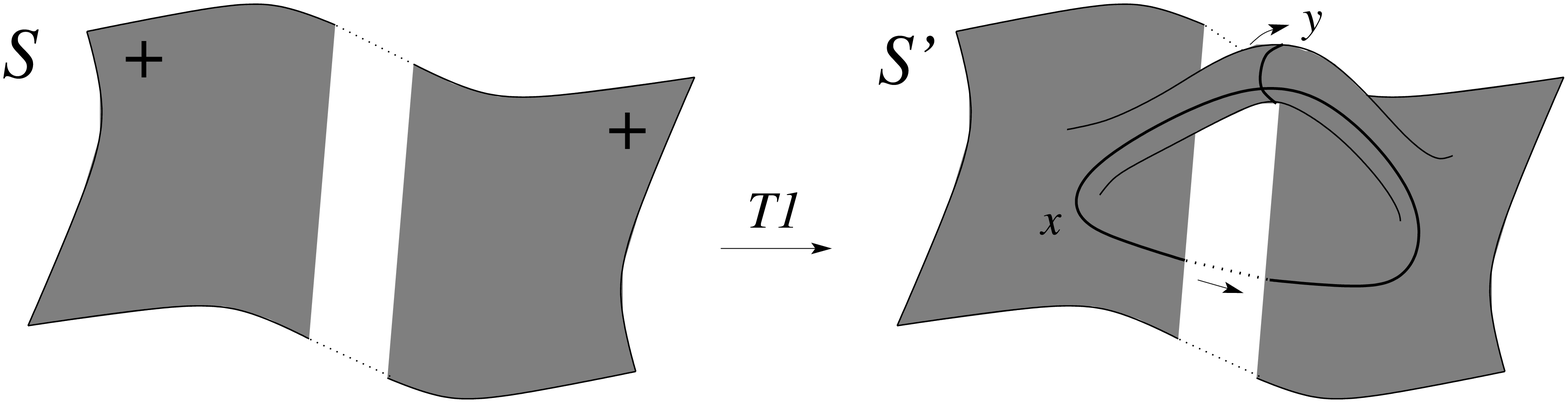,height=2.3cm}
\end{center}
\end{figure}
\noindent This time, the Hermitian matrices $H$ and $H'$ are related by
\[
H'(\omega)=\begin{pmatrix}
H(\omega)&\ast&0\cr \ast&\ast&\alpha\cr
0&\overline\alpha&0
\end{pmatrix},
\]
with $\alpha=(1-\overline\omega_k)\prod_{i\neq k}|1-\omega_i|^2$.
Since $H'(\omega)$ is Hermitian and $\omega_i\neq 1$ for all $i$,
$H'(\omega)$ is an elementary enlargement of $H(\omega)$.

\smallskip\noindent{$(T2)$} Let $S'$ be a C-complex obtained from $S$ by the
transformation $T2$. If this transformation connects two distinct connected components
of $S$, then $H_1(S')=H_1(S)\oplus\Zz z$. As above, we get $H'=H\oplus(0)$ and $\beta_0(S')=\beta_0(S)-1$,
so the signature and nullity are unchanged. On the other hand, if this transformation take place on one
connected component of $S$, then $H_1(S')=H_1(S)\oplus\Zz w\oplus\Zz z$ with $w$ and $z$ as illustrated below.

\begin{figure}[h]
\begin{center}
\epsfig{figure=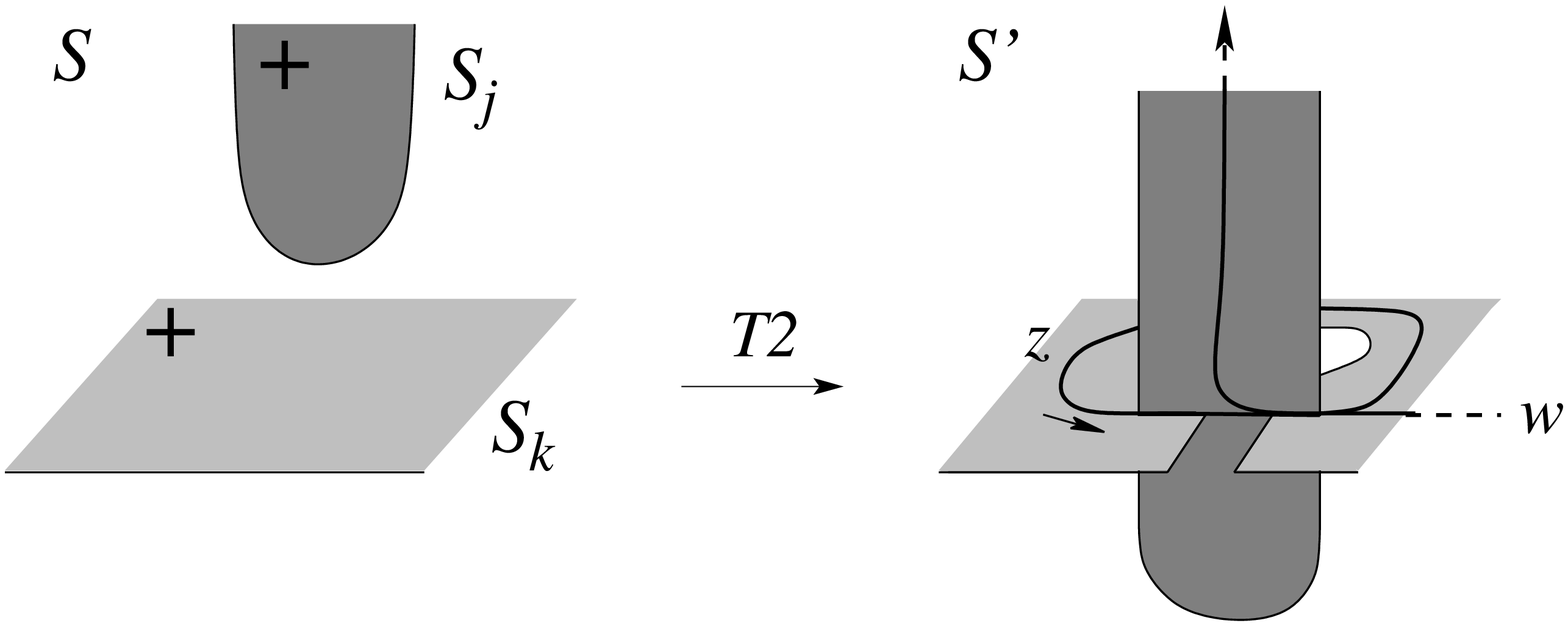,height=3cm}
\end{center}
\end{figure}
\noindent Therefore,
\[
A_{S'}^\eps=\begin{pmatrix}
A_S^\eps&\ast&0\cr \ast&\ast&\delta(-\eps)\cr
0&\delta(\eps)&0
\end{pmatrix}\,,
\quad\hbox{with}\quad\delta(\eps)=\begin{cases}1&\text{if $\eps_j=\eps_k=+1$;}\cr 0&\text{else.}\end{cases}
\]
It follows that the corresponding Hermitian matrix $H'(\omega)$ is an elementary
enlargement of $H(\omega)$ with
$\alpha=(1-\overline\omega_j)(1-\overline\omega_k)\prod_{i\neq j,k}|1-\omega_i|^2$. 

\smallskip\noindent{$(T3)$} Finally, let $S$ and $S'$ be C-complexes related by the move $T3$.
A similar computation shows that the corresponding Hermitian matrices
$H(\omega)$ and $H'(\omega)$ are both
an elementary enlargement of some Hermitian matrix. This concludes the proof.
\end{proof}

\begin{ex}
If $\mu=1$, then the colored link $L$ is just a link. Furthermore, a C-complex $S$ for $L$
is nothing but a (possibly disconnected) Seifert surface for $L$. Finally, $A^-$ is a usual Seifert matrix $A$, and $A^+=A^T$.
Hence, the corresponding Hermitian matrix is given by
\[
H(\omega)=(1-\overline\omega)(A^T-\omega A)=(1-\omega)A+(1-\overline\omega)A^T.
\]
So if $\mu=1$, the signature of $L$ is the Levine-Tristram signature of the link $L$.
\end{ex}

\begin{figure}[h]
\begin{center}
\epsfig{figure=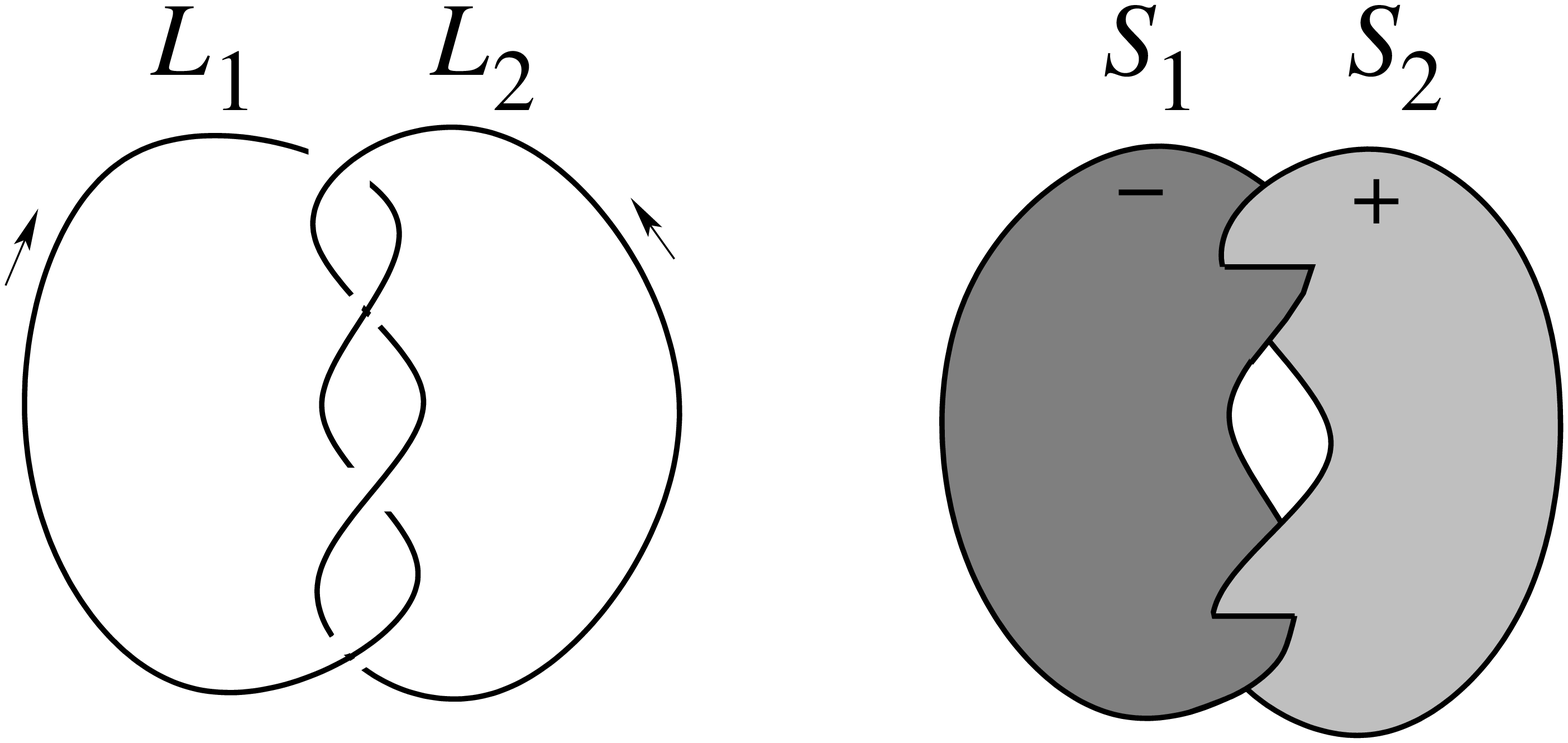,height=2.2cm}
\end{center}
\end{figure}
\begin{ex}\label{ex:lk2}
Consider the $2$-colored link $L$ illustrated above. A C-complex $S$ for $L$ is also given.
We compute $A^\eps=(-1)$ if $\eps_1=\eps_2$, and $A^\eps=(0)$ else.
Hence,
\[
H(\omega_1,\omega_2)=(1-\overline\omega_1)(1-\overline\omega_2)(-1-\omega_1\omega_2)=
-2\Re((1-\omega_1)(1-\omega_2)).
\]
\noindent\parbox{0.75\textwidth}{\noindent So $\sigma_L(\omega_1,\omega_2)$ is given by the sign of
$-\Re((1-\omega_1)(1-\omega_2))$. Furthermore,
$\eta_L(\omega_1,\omega_2)=1$ if $\omega_1\omega_2=-1$, and $\eta_L(\omega_1,\omega_2)=0$ else.
Let us draw the domain $T_\ast^2$ as a square. The value of the function
$\sigma_L$ can be represented as illustrated opposite. Note that $\sigma_L$ and $\eta_L$ are constant on the connected components
of the complement of the zeroes of $\Delta_L(t_1,t_2)=t_1t_2+1$, the Alexander polynomial of $L$.
We shall explain this fact in Section \ref{section:PC}.}
\hfill\parbox{0.22\textwidth}{\epsfig{figure=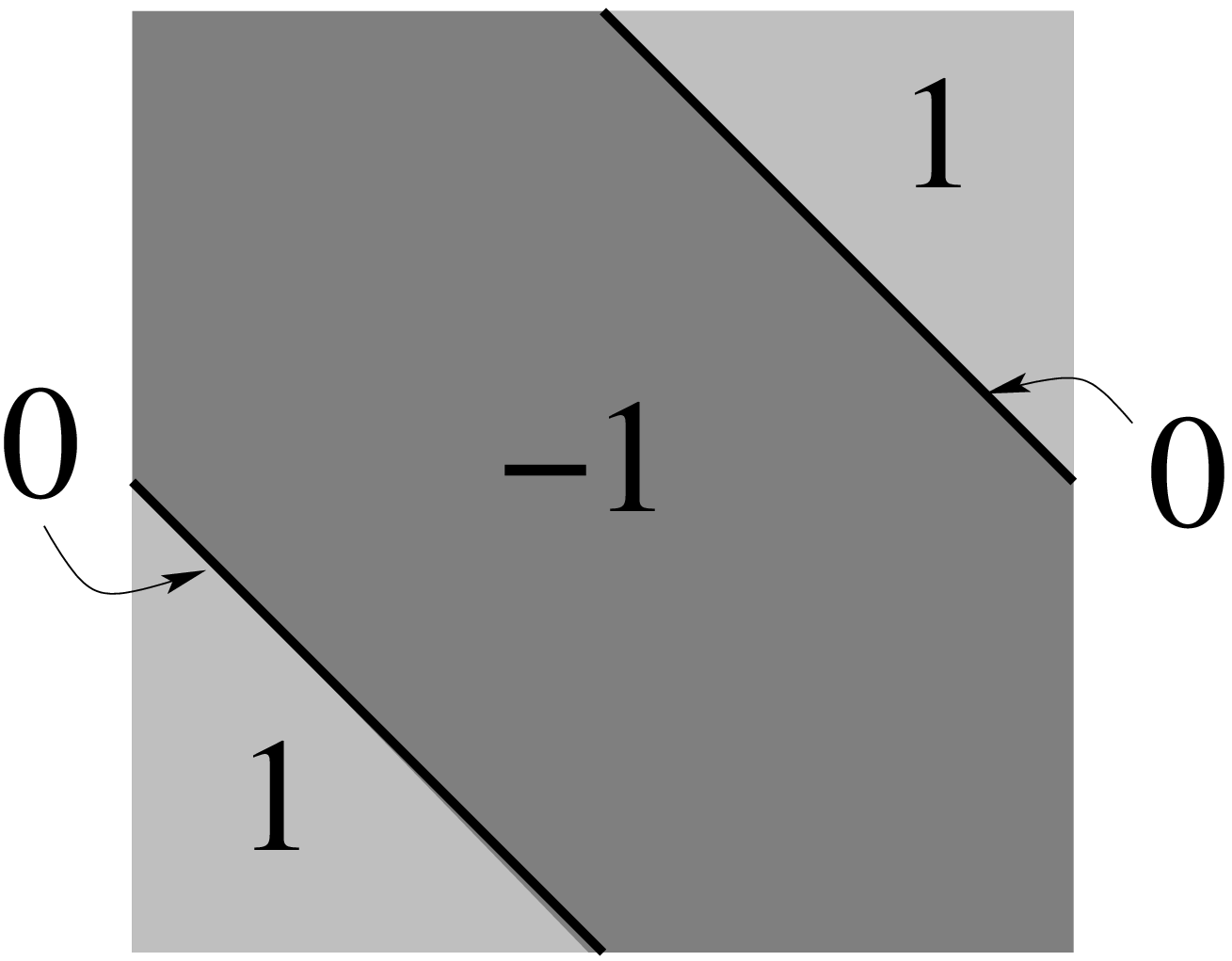,height=2cm}} 
\end{ex}

\begin{prop}\label{prop:lk}
Let $L=L_1\cup\dots\cup L_{\mu+1}$ be a $(\mu+1)$-colored link. Consider the $\mu$-colored link $L'=L'_1\cup\dots\cup L'_\mu$
given by $L'_i=L_i$ for $i<\mu$ and $L'_\mu=L_\mu\cup L_{\mu+1}$. Then, for all $(\omega_1,\dots,\omega_\mu)$ in $T^\mu_*$,
\begin{eqnarray*}
\sigma_{L'}(\omega_1,\dots,\omega_\mu)&=&\sigma_L(\omega_1,\dots,\omega_\mu,\omega_\mu)-lk(L_\mu,L_{\mu+1}),\\
\eta_{L'}(\omega_1,\dots,\omega_\mu)&=&\eta_L(\omega_1,\dots,\omega_\mu,\omega_\mu).
\end{eqnarray*}
\end{prop}

Before giving the proof of this proposition, let us point out an interesting consequence:
it is possible to compute the signature and nullity of a $\mu$-colored link by considering any
finer coloring of the same underlying link. In particular, all the signatures (corresponding to all the possible colorings)
can be computed from the signature corresponding to a coloring with the maximal number of colors. This simplifies greatly the
computations in many cases, as illustrated by the following (didactic) example.

\begin{ex}
Let us try to compute the Levine-Tristram signature of the link $L$ illustrated below.
One possibility is to choose a Seifert surface for $L$ and to compute the corresponding Seifert matrix.

\noindent\parbox{0.65\textwidth}{\noindent On the other hand,
consider a $3$-colored link $L'$ obtained by coloring the components of $L$ with three different colors.
There is an obvious contractible
C-complex for $L'$, so $\sigma_{L'}$ is identically zero. By Proposition
\ref{prop:lk}, the Levine-Tristram signature of $L$ is given  by $\sigma_L(\omega)=\sigma_{L'}(\omega,\omega,\omega)-2=-2$.}
\hfill\parbox{0.3\textwidth}{\epsfig{figure=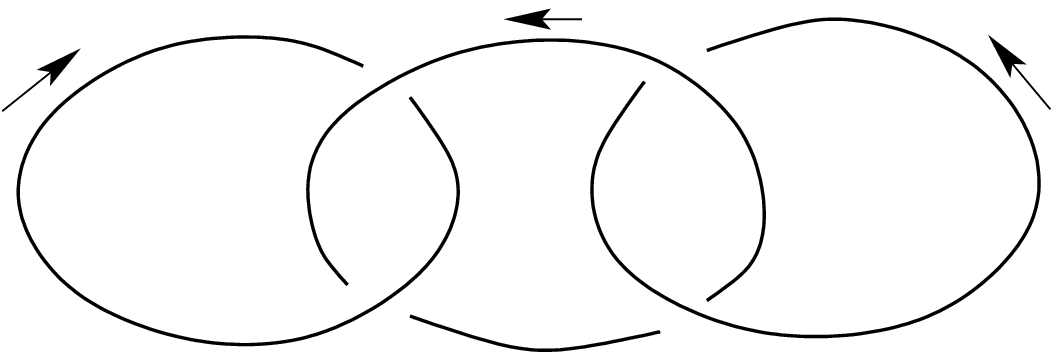,height=1.2cm}}
\end{ex}

\begin{proof}[Proof of Proposition \ref{prop:lk}]
First, note that it is sufficient to prove this statement when $L_{\mu+1}$ is a knot. Consider a C-complex
$S=S_1\cup\dots\cup S_{\mu+1}$ for $L$, and let $\ell$ be the number of clasps in $S_\mu\cap S_{\mu+1}$. 
A certain number of these clasps (say, $\ell_+$) induce a positive linking number between $L_\mu$ and $L_{\mu+1}$, while the $\ell_-$
remaining ones induce a negative linking number. By definition, $\ell=\ell_++\ell_-$ and $lk(L_\mu,L_{\mu+1})=\ell_+-\ell_-$.
Using transformation $T2$ of Figure \ref{fig:Sequiv}, it may be assumed that $\ell_+$ and $\ell_-$ are positive. Via handle
attachment, it may also be assumed that $S_\mu$ is connected. Finally, one easily checks that the C-complex $S$ can be chosen so
that the knot $L_{\mu+1}$ crosses the $\ell_-$ negative clasps first, and then the $\ell_+$ positive ones. This situation is illustrated
below.
\begin{figure}[h]
\begin{center}
\epsfig{figure=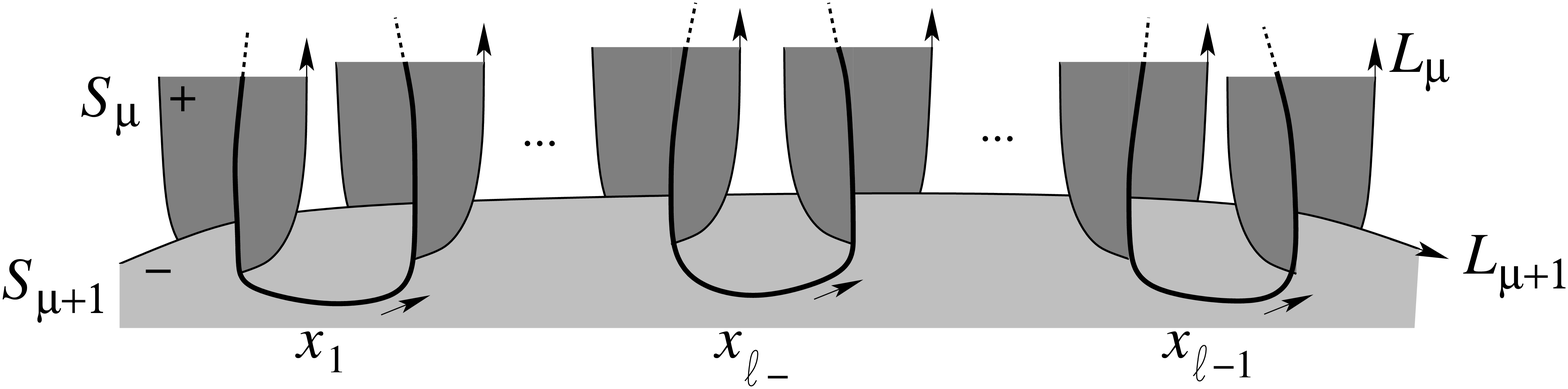,height=2.5cm}
\end{center}
\end{figure}

Let $\widetilde S$ be the C-complex obtained from $S$ by removing these $\ell$ clasps. Since $S_\mu$ and $S_{\mu+1}$ are connected,
$H_1(S)=\bigoplus_{i=1}^{\ell-1}\Zz x_i\oplus H_1(\widetilde S)$, with $x_i$ the $1$-cycles in $S_\mu\cup S_{\mu+1}$ depicted above.
Fix a basis $\mathcal B$ of $H_1(\widetilde S)$. For any sequence $\eps=(\eps_1,\dots,\eps_{\mu+1})$ of $\pm 1$'s,
the Seifert matrix $A_S^\eps$ corresponding to the basis $\left<x_1,\dots,x_{\ell-1}\right>\cup {\mathcal B}$ of $H_1(S)$ can be written
\[
A_S^\eps=\begin{pmatrix}
D^{\eps_\mu\eps_{\mu+1}}&C^{\eps_\mu}\cr(C^{-\eps_\mu})^T&A_{\widetilde S}^\eps
\end{pmatrix}. 
\] 
Fix an element $\omega'=(\omega_1,\dots,\omega_\mu)$ of $T^\mu_*$. Since $D^{\eps_\mu\eps_{{\mu+1}}}$ only depends on
$\eps_\mu$ and $\eps_{{\mu+1}}$, the upper left block of the corresponding Hermitian matrix $H_S(\omega',\omega_\mu)$ is given by
\[
\prod_{i=1}^\mu|1-\omega_i|^2(\lambda D+\overline\lambda D^T+N),
\]
where $\lambda=1-\omega_\mu$, $D=D^{--}$ and
\[
N=(D^{-+}-D)+(D^{-+}-D)^T=
\begin{pmatrix}
-2&1\cr \phantom{-}1&\ddots&\ddots\cr &\ddots&-2&\!\!\phantom{-}1\cr
&&\phantom{-}1&\!\!\phantom{-}0&-1\cr &&&\!\!-1&\phantom{-}2&\ddots\cr &&&&\ddots&\ddots&-1\cr &&&&&-1&\phantom{-}2
\end{pmatrix}.
\eqno{(\star)}
\]
Similarly, the upper right block of $H_S(\omega',\omega_\mu)$ is equal to
$\prod_{i=1}^\mu|1-\omega_i|^2(\lambda C^-+\overline\lambda C^+)$. Finally, observe that each coefficient of $A_{\widetilde S}^\eps$
is either independant of $\eps_\mu$, or independant of $\eps_{\mu+1}$. Therefore,
$H_{\widetilde S}(\omega',\omega_\mu)=|\lambda|^2H_{\widetilde S'}(\omega')$, where
$\widetilde S'=S_1\cup\dots\cup S_{\mu-1}\cup\widetilde S_\mu'$ and
$\widetilde S_\mu'=\widetilde S_\mu\sqcup\widetilde S_{\mu+1}$. To sum up, we have the equality
\[
H_S(\omega',\omega_\mu)=\prod_{i=1}^\mu|1-\omega_i|^2
\begin{pmatrix}
\lambda D+\overline\lambda D^T +N&\lambda C^-+\overline\lambda C^+\cr
(\overline\lambda C^-+\lambda C^+)^T&\prod_{i=1}^{\mu-1}|1-\omega_i|^{-2}H_{\widetilde S'}(\omega')
\end{pmatrix}.\eqno{(\star\star)}
\]
Let us now turn to the $\mu$-colored link $L'$. The C-complex $S$ can be transformed into a C-complex
$S'=S_1'\cup\dots\cup S'_\mu$ for $L'$ as follows: set $S_i'=S_i$ for $i<\mu$, and let $S'_\mu$ be the surface obtained from
$S_\mu\cup S_{\mu+1}$ by `smoothing' the $\ell$ clasps as illustrated below.

\begin{figure}[h]
\begin{center}
\epsfig{figure=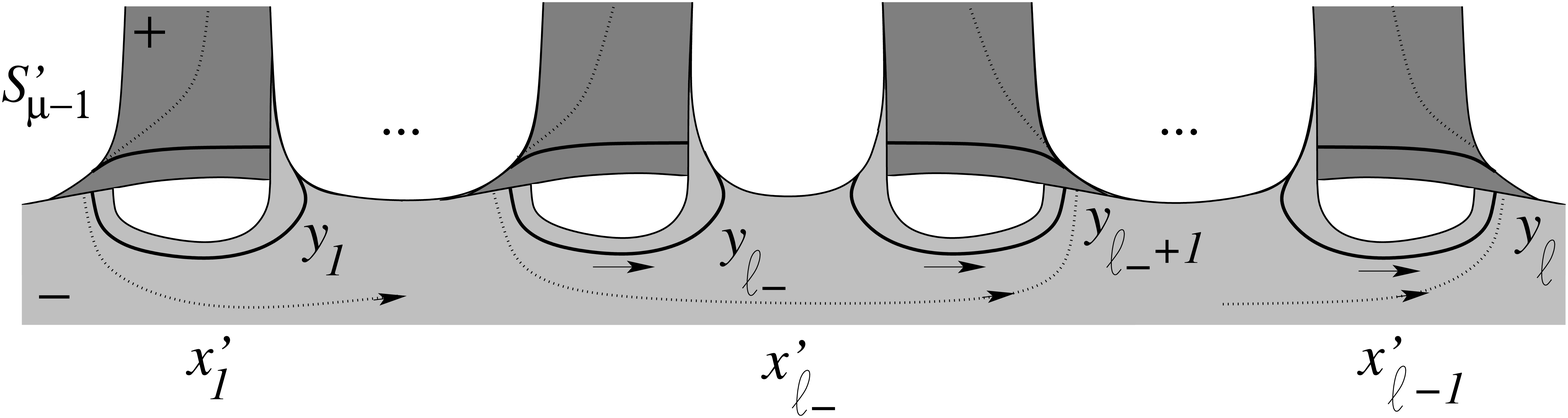,height=3cm}
\end{center}
\end{figure}
Basically, each clasp is replaced by two half-twisted bands joining $S_\mu$ and $S_{\mu+1}$.
This time, $H_1(S')=\bigoplus_{i=1}^\ell\Zz y_i\oplus\bigoplus_{i=1}^{\ell-1}\Zz x'_i\oplus H_1(\widetilde S')$ with the $1$-cycles
$y_i$ and $x_i'$ depicted above. Using the same method and notations as above, one can compute the Hermitian matrix
$H_{S'}(\omega')$ corresponding to the basis $\left<y_1,\dots,y_\ell,x_1'\dots,x_{\ell-1}'\right>\cup{\mathcal B}$ of $H_1(S')$.
It is given by
\[
\prod_{i=1}^\mu|1-\omega_i|^2
\begin{pmatrix}
|\lambda|^2(I_{\ell_-}\oplus -I_{\ell_+})&\overline\lambda M&0\cr
	\lambda M^T&\lambda D+\overline\lambda D^T&\lambda C^-+\overline\lambda C^+\cr
	0&(\overline\lambda C^-+\lambda C^+)^T&\prod_{i=1}^{\mu-1}|1-\omega_i|^{-2}H_{\widetilde S'}(\omega')
\end{pmatrix},
\]
where $I_k$ denotes the $(k\times k)$-identity matrix and
\[
M=\begin{pmatrix}
\phantom{-}1\cr -1&\ddots\cr &\ddots&\!\!\phantom{-}1\cr &&\!\!-1
\end{pmatrix}.
\]
Using the equalities $(\star)$ and $(\star\star)$ together with the fact that $\lambda\neq 0$, it is easy to check that the matrix
$H_{S'}(\omega')$ is conjugate to the Hermitian matrix $I_{\ell_-}\oplus -I_{\ell_+}\oplus H_S(\omega',\omega_\mu)$. Hence,
\begin{eqnarray*}
\sigma_{L'}(\omega')&=&\sign(H_{S'}(\omega'))=\sign(H_S(\omega',\omega_\mu))+\ell_--\ell_+\cr
	&=&\sigma_L(\omega',\omega_\mu)-lk(L_\mu,L_{\mu+1}),
\end{eqnarray*}
and $\nuli H_{S'}(\omega')=\nuli H_S(\omega',\omega_\mu)$. Since $\beta_0(S)=\beta_0(S')$, the proposition is proved.
\end{proof}

Note that we shall give an alternative proof of this result in subsection \ref{sub:proof4dim}.

\subsection{Basic properties of $\sigma_L$ and $\eta_L$}\label{sub:basic}

We conclude this section with an enumeration of several properties of the signature and nullity.
Deeper properties shall be presented in the following sections.
The first proposition follows easily from the fact that $H(\omega)$ is Hermitian.

\begin{prop}
For any $(\omega_1,\dots,\omega_\mu)\in T_\ast^\mu$,
\begin{eqnarray*}
\sigma_L(\omega^{-1}_1,\dots,\omega^{-1}_\mu)&=&\sigma_L(\omega_1,\dots,\omega_\mu),\\
\eta_L(\omega^{-1}_1,\dots,\omega^{-1}_\mu)&=&\eta_L(\omega_1,\dots,\omega_\mu).\qed
\end{eqnarray*}
\end{prop}

\begin{prop}
Let $L$ be a colored link, and let $L'$ be the colored link obtained from $L$ by reversing the orientation
of every component of the sublink $L_1$ of $L$. Then, for all $\omega=(\omega_1,\dots,\omega_\mu)$ in $T^\mu_\ast$,
\begin{eqnarray*}
\sigma_{L'}(\omega_1,\dots,\omega_\mu)&=&\sigma_L(\omega^{-1}_1,\omega_2,\dots,\omega_\mu),\\
\eta_{L'}(\omega_1,\dots,\omega_\mu)&=&\eta_L(\omega_1^{-1},\omega_2,\dots,\omega_\mu).
\end{eqnarray*}
\end{prop}
\begin{proof}
If $S=S_1\cup S_2\cup\dots\cup S_\mu$ is a C-complex for $L$, then $S'=(-S_1)\cup S_2\cup\dots\cup S_\mu$
is a C-complex for $L'$. Since $A_{S'}^{\eps'}=A_S^\eps$ with $\eps_1'=-\eps_1$ and $\eps'_i=\eps_i$
for $i>1$, the corresponding Hermitian matrices $H'$ and $H$ are related by
$H'(\omega_1,\omega_2,\dots,\omega_\mu)=H(\omega^{-1}_1,\omega_2,\dots,\omega_\mu)$.
\end{proof}

\begin{cor}
If $-L$ denotes the colored link $L$ with the opposite orientation, then $\sigma_{-L}=\sigma_L$ and $\eta_{-L}=\eta_L$.\qed
\end{cor}

So the signature cannot distinguish between a colored link and its inverse. On the other hand, it is a useful invariant
for telling apart a colored link and its mirror image.

\begin{prop}
If $\overline L$ denotes the mirror image of the colored link $L$, then $\sigma_{\overline L}=-\sigma_L$ and $\eta_{\overline L}=\eta_L$.
\end{prop}
\begin{proof}
If $S$ is a C-complex for $L$, then the mirror image $S'$ of $S$ is a C-complex for $\overline L$. Therefore,
$A_{S'}^\eps=-A_S^\eps$ and the Hermitian matrices $H$ and $H'$ satisfy $H'(\omega)=-H(\omega)$.
\end{proof}

\begin{cor}\label{cor:mirror}
If a colored link $L$ is isotopic to its mirror image, then $\sigma_L$ is identically zero.\qed
\end{cor}

Finally, the signature and nullity behave well under connected and disjoint sum.

\begin{prop}
Let $L'=L_1\cup\dots\cup L_{\nu-1}\cup L'_\nu$ and $L''=L''_\nu\cup L_{\nu+1}\cup\dots\cup L_\mu$ be two
colored links. Consider a colored link $L=L_1\cup\dots\cup L_\mu$, where $L_\nu$ is a
connected sum of $L'_\nu$ and $L''_\nu$ along any of their components.
Then, for all $\omega_1,\dots,\omega_\mu\in S^1\setminus\{1\}$,
\begin{eqnarray*}
\sigma_L(\omega_1,\dots,\omega_\mu)&=&
	\sigma_{L'}(\omega_1,\dots,\omega_\nu)+\sigma_{L''}(\omega_\nu,\dots,\omega_\mu),\\
\eta_L(\omega_1,\dots,\omega_\mu)&=&
	\eta_{L'}(\omega_1,\dots,\omega_\nu)+\eta_{L''}(\omega_\nu,\dots,\omega_\mu).
\end{eqnarray*}
\end{prop}
\begin{proof}
Given $S'$ a C-complex for $L'$ and $S''$ a C-complex for $L''$, a C-complex $S$ for $L$ is given by
the band sum of $S'$ and $S''$ along the corresponding components of $S'_\nu$ and $S''_\nu$. Since $S'$ and $S''$
have no closed components, this band can be chosen such that only its ends meet $S'$ and $S''$. Clearly,
$A^\eps_S=A_{S'}^{\eps'}\oplus A_{S''}^{\eps''}$ with $\eps'=(\eps_1,\dots,\eps_\nu)$ and
$\eps''=(\eps_\nu,\dots,\eps_\mu)$. The corresponding Hermitian matrices $H$, $H'$ and $H''$ satisfy
\[
H(\omega_1,\dots,\omega_\mu)=\prod_{i>\nu}|1-\omega_i|^2 H'(\omega_1,\dots,\omega_\nu)\oplus
\prod_{i<\nu}|1-\omega_i|^2 H''(\omega_\nu,\dots,\omega_\mu).
\]
Since $\omega_i\neq 1$ for all $i$ and $\beta_0(S)=\beta_0(S')+\beta_0(S'')-1$, this implies the proposition.
\end{proof}

\begin{prop}
Let $L'$ and $L''$ be colored links with disjoint sets of colors. Consider the colored link $L$
given by the disjoint sum of $L'$ and $L''$. Then,
\begin{eqnarray*}
\sigma_L(\omega',\omega'')&=&\sigma_{L'}(\omega')+\sigma_{L''}(\omega''),\\
\eta_L(\omega',\omega'')&=&\eta_{L'}(\omega')+\eta_{L''}(\omega'')+1.
\end{eqnarray*}
\end{prop}
\begin{proof}
Let $S'$ and $S''$ be C-complexes for $L'$ and $L''$. A C-complex $S$ for $L$ is given by the disjoint union of $S'$ and $S''$.
Clearly, $A_S^{(\eps',\eps'')}=A_{S'}^{\eps'}\oplus A_{S''}^{\eps''}$, so
\[
H(\omega',\omega'')=\prod_i|1-\omega''_i|^2 H'(\omega')\oplus\prod_j|1-\omega'_j|^2 H''(\omega'').
\]
Since $\beta_0(S)=\beta_0(S')+\beta_0(S'')$, the proposition is proved.
\end{proof}

\section{The Alexander module of a colored link}\label{section:Alex}

Associated to a $\mu$-colored link $L$ with exterior $X$ is a natural $\Zz^\mu$-covering $\widetilde X\to X$. The aim of
this section is to show how the space $\widetilde X$ can be constructed from a C-complex for $L$. This leads to a presentation of
the Alexander module of $L$, that is, the $\Zz[t_1^{\pm 1},\dots,t_\mu^{\pm 1}]$-module $H_1(\widetilde{X})$. This generalizes
a celebrated theorem of Seifert, which corresponds to the case $\mu=1$.

We shall use these results in the next section to derive relations between the signature, the nullity, and the Alexander invariants of
colored links.

\subsection{Basics}\label{sub:basics}

Let $L=L_1\cup\dots\cup L_\mu$ be a colored link, and let $X$ denote its exterior. The epimorphism $\pi_1(X)\to\Zz^\mu$
given by $\gamma\mapsto(lk(\gamma,L_1),\dots,lk(\gamma,L_\mu))$ induces a regular $\Zz^\mu$-covering $\widetilde{X}\to X$.
The homology of $\widetilde{X}$ is a natural module over the ring $\Lambda_\mu=\Zz[t_1^{\pm 1},\dots,t_\mu^{\pm 1}]$, where
$t_i$ denotes the covering transformation corresponding to an oriented meridian of $L_i$. The $\Lambda_\mu$-module
$H_1(\widetilde{X})$ is called the {\em Alexander module\/} of the colored link $L$. Of course, if $\mu=1$, then
$H_1(\widetilde{X})$ is nothing but the usual Alexander module of the link $L$.

To extract handy information from such a cumbersome invariant as a module over $\Lambda_\mu$, the standard trick is to consider
its elementary ideals. Although these objects are widely used, their definitions vary according to the authors.
Therefore, we shall now clarify the meaning of these concepts in the present work.
Let $R$ be a noetherian factorial domain, that is, an integral domain in which every ideal is finitely generated, and every non-zero
non-invertible element has a unique factorization.
Let $P$ be an $m\times n$ matrix with coefficients in $R$. Let us denote by $E_r(P)$ the ideal of $R$ generated by
all the $(m-r)\times(m-r)$ minors of $P$. By convention, $E_r(P)=(0)$ if $r<0$ and $E_r(P)=R$ if $r\ge m$. Let $\Delta_r(P)$ denote the
greatest common divisor of the elements of $E_r(P)$. (Recall that the greatest common divisor of $a_1,\dots,a_n$ in $R$ is an element
$d$ of $R$ which divides $a_i$ for all $i$, and such that if $c\in R$ divides $a_i$ for all $i$, then $c$ divides $d$.)
Since $R$ is a noetherian factorial domain, $\Delta_r(P)$ exists, and is well-defined up to multiplication
by a unit of $R$. Given $\Delta$ and $\Delta'$ in $R$, let us note $\Delta\,\dot{=}\,\Delta'\,$ if
$\,\Delta=u\Delta'$ for some unit $u$ of $R$.
Now, let $M$ be a module over a factorial ring $R$. A {\em finite presentation\/} of $M$ is an exact sequence
$F\stackrel{\varphi}{\to}E\to M\to 0$, where $E$ and $F$ are free $R$-modules with finite basis. A matrix of $\varphi$
is a {\em presentation matrix\/} of $M$. The $r^{th}$ {\em elementary ideal\/} of $M$ is the ideal of $R$
given by $E_r(M):=E_r(P)$, where $P$ is any presentation matrix of $M$. It is easy to check that these ideals do not depend on
the presentation of $M$. In particular, the element $\Delta_r(M):=\Delta_r(P)$ of $R$ is well defined up to multiplication by a unit
of $R$.

Let us now turn back to colored links. Given a $\mu$-colored link $L$, we just defined the Alexander module of $L$ as some module
$H_1(\widetilde{X})$ over the ring $\Lambda_\mu$. The $r^\mathrm{th}$ elementary ideal $E_r(H_1(\widetilde{X}))$ is the
$r^{th}$ {\em Alexander ideal\/} of $L$, and is denoted by $E_r(L)$. The polynomial $\Delta_r(L):=\Delta_r(H_1(\widetilde{X}))$ is
called the $r^{th}$ {\em Alexander polynomial\/} of $L$; $\Delta_0(L)$ is called \underline{the} {\em Alexander polynomial\/}
of $L$, and is denoted by $\Delta_L$. Again, note that $\Delta_r(L)$ is only defined up to multiplication by a unit
of $\Lambda_\mu$, that is, by $\pm t_1^{m_1}\cdots t_\mu^{m_\mu}$ with integers $m_i$.

\subsection{A presentation of the Alexander module using a C-complex}\label{sub:pres}

Fix a colored link $L=L_1\cup\dots\cup L_\mu$. Consider a C-complex $S=S_1\cup\dots\cup S_\mu$ for $L$
such that each $S_i$ is connected and $S_i\cap S_j\neq\emptyset$ for all $i\neq j$.
(Such a C-complex exists by transformations $T1$ and $T2$ of Lemma \ref{lemma:Sequiv}.) For $i=1,\dots,\mu$, choose
some interior point $v_i$ of $S_i\setminus\bigcup_{j\neq i}S_i\cap S_j$. Given a clasp in $S_i\cap S_j$ with $i<j$, consider
an oriented edge in $S_i\cup S_j$ joining $v_i$ and $v_j$ and passing through this single clasp as described in Figure \ref{fig:loop}.
This leads to a collection of oriented edges $\{e_{ij}^1,\dots,e_{ij}^{c(i,j)}\}$, where $c(i,j)$ denotes the number of clasps
in $S_i\cap S_j$ (that is: the number
of connected components of $S_i\cap S_j$). Let $K_{ij}\subset S_i\cup S_j$ denote the graph given by the union of these edges.
Finally let $K_\mu$ be the complete graph with vertices $\{v_i\}_{1\le i\le\mu}$ and edges $\{e_{ij}^1\}_{1\le i<j\le\mu}$.

\begin{lemma}\label{lemma:dec}
The homology of $S=S_1\cup\dots\cup S_\mu$ is equal to
\[
H_1(S)=\bigoplus_{1\le i\le\mu}H_1(S_i)\;\oplus\bigoplus_{1\le i<j\le\mu}H_1(K_{ij})\;\oplus\;H_1(K_\mu).
\]
Furthermore, a basis of $H_1(K_{ij})$ is given by $<\beta_{ij}^\ell>_{1\le\ell\le c(i,j)-1}$, where
$\beta_{ij}^\ell=e_{ij}^\ell-e_{ij}^{\ell+1}$. Finally, a basis of $H_1(K_\mu)$ is given by
$<\gamma_{1ij}>_{2\le i<j\le\mu}$, where $\gamma_{ijk}=e^1_{ij}-e^1_{ki}+e^1_{jk}$.
\end{lemma}
\begin{proof}
The C-complex $S$ can be constructed as follows. Consider the complete graph $K_\mu$. Add the graphs $K_{ij}$ one by one, for
$1\le i<j\le\mu$. Finally, paste $S_1,S_2,\dots,S_\mu$. Note that at each step, the pasting is done along a contractible space.
A recursive use of the Mayer-Vietoris exact sequence therefore leads to the first statement of the lemma. The fact that 
$H_1(K_{ij})=\bigoplus_{\ell=1}^{c(i,j)-1}\Zz\beta_{ij}^\ell$ is clear. Finally, we have the relation
$\gamma_{ijk}=\gamma_{1jk}-\gamma_{1ik}+\gamma_{1ij}$. Hence the family $<\gamma_{1ij}>_{2\le i<j\le\mu}$ generates
$H_1(K_\mu)$. Furthermore, $1-\rk H_1(K_\mu)=\chi(K_\mu)=\mu-\genfrac(){0cm}{1}\mu 2$, so $\rk H_1(K_\mu)=\genfrac(){0cm}{1}{\mu-1}2$
and this family is a basis of $H_1(K_\mu)$.
\end{proof}

Recall the homomorphism $i^\eps\colon H_1(S)\to H_1(S^3\setminus S)$ of subsection \ref{sub:C-cplx}.
\begin{thm}\label {thm:module}
Let $L=L_1\cup\dots\cup L_\mu$ be a colored link, and consider a C-complex $S=S_1\cup\dots\cup S_\mu$ for $L$ such that
each $S_i$ is connected and $S_i\cap S_j\neq\emptyset$ for all $i\neq j$.
Let $\alpha\colon H_1(S)\otimes_\Zz\Lambda_\mu\to H_1(S^3\setminus S)\otimes_\Zz\Lambda_\mu$
be the homomorphism of $\Lambda_\mu$-modules given by
\[
\alpha=\sum_\eps\eps_1\cdots\eps_\mu\,t_1^{\frac{\eps_1+1}{2}}\cdots t_\mu^{\frac{\eps_\mu+1}{2}}i^\eps,
\]
where the sum is on all sequences $\eps=(\eps_1,\dots,\eps_\mu)$ of $\pm 1$'s. Then, the
Alexander module $H_1(\widetilde{X})$ of $L$ admits the finite presentation
\[
\widehat{H}\otimes_\Zz\Lambda_\mu\stackrel{\widehat\alpha}{\longrightarrow}H_1(S^3\setminus S)\otimes_\Zz\Lambda_\mu\longrightarrow
H_1(\widetilde{X})\longrightarrow 0,
\]
where $\widehat{H}=\bigoplus_{1\le i\le\mu}H_1(S_i)\oplus\bigoplus_{1\le i<j\le\mu}H_1(K_{ij})\oplus
\bigoplus_{1\le i<j<k\le\mu}\Zz\gamma_{ijk}$ and $\widehat\alpha$ is given by
\begin{itemize}
\item{$\widehat\alpha=\prod_{n\neq i}(t_n-1)^{-1}\alpha$ on $H_1(S_i)$ for $1\le i\le\mu$;}
\item{$\widehat\alpha=\prod_{n\neq i,j}(t_n-1)^{-1}\alpha$ on $H_1(K_{ij})$ for $1\le i<j\le\mu$;}
\item{$\widehat\alpha(\gamma_{ijk})=\prod_{n\neq i,j,k}(t_n-1)^{-1}\alpha(\gamma_{ijk})$ for $1\le i<j<k\le\mu$.}
\end{itemize}
\end{thm}

We postpone the proof of this theorem to the end of the section.

For colored links with $1$, $2$ or $3$ colors, this result provides a square presentation matrix of the Alexander module expressed
in terms of the Seifert matrices $A^\eps$. More precisely, we have the following corollaries.

\begin{cor}[Seifert \cite{Sei}]
Let $A$ be a Seifert matrix for the link $L$. Then, $tA-A^T$ is a presentation matrix of the Alexander module of $L$.
\end{cor}
\begin{proof}
Theorem \ref{thm:module} gives the finite presentation
\[
H_1(S)\otimes_\Zz\Lambda\stackrel{\widehat\alpha}{\to}H_1(S^3\setminus S)\otimes_\Zz\Lambda\to H_1(\widetilde{X})\to 0,
\]
where $\widehat\alpha=\alpha=ti^+-i^-$. By Alexander duality, a matrix of $i^+$ (resp. $i^-$) is given by the transpose of $A^+$
(resp. $A^-$). Therefore, $tA-A^T$ is a presentation matrix.
\end{proof}

Similarly, we get the following result.

\begin{cor}[Cooper \cite{Coo}]\label{cor:mod2}
Let $S=S_1\cup S_2$ be a C-complex for a colored link $L=L_1\cup L_2$. Let $A$ (resp. $B$) be a matrix of the form $\alpha^{--}$
(resp. $\alpha^{-+}$) with respect to a basis of $H_1(S)$ adapted to the decomposition
$H_1(S)=H_1(S_1)\oplus H_1(S_2)\oplus H_1(K_{12})$. Then, a presentation matrix of the Alexander module of $L$ is given by
\[
(t_1t_2A-t_1B-t_2B^T+A^T)\cdot D,
\]
where $D=(D_{ij})$ is the diagonal matrix given by
\[
D_{ii}=\begin{cases}(t_2-1)^{-1}&\text{if $1\le i\le\beta_1(S_1)$;}\cr
	(t_1-1)^{-1}&\text{if $\beta_1(S_1)<i\le\beta_1(S_1)+\beta_1(S_2)$;}\cr
	1&\text{if $\beta_1(S_1)+\beta_1(S_2)<i\le\beta_1(S)$,}
	\end{cases} 
\]
and $\beta_1(\cdot)$ denotes the first Betti number.\qed
\end{cor}

The case $\mu=3$ is similar, but the complete statement is a little cumbersome. Instead, let us give an example of such a computation.

\begin{ex}\label{ex:3-color}
Consider the $3$-colored link $L$ given in Figure \ref{fig:3-color}. A C-complex $S$ for $L$ is also drawn. Clearly, $H_1(S)=\Zz\gamma$.
Furthermore, $lk(\gamma^\eps,\gamma)=1$ if $\eps_1=\eps_2=\eps_3$, and all the other linking numbers are zero.
Hence, a matrix of $\widehat\alpha=\alpha$ is given by $(t_1t_2t_3-1)$, and $H_1(\widetilde{X})=\Lambda_3/(t_1t_2t_3-1)$.
\begin{figure}[Htb]
\begin{center}
\epsfig{figure=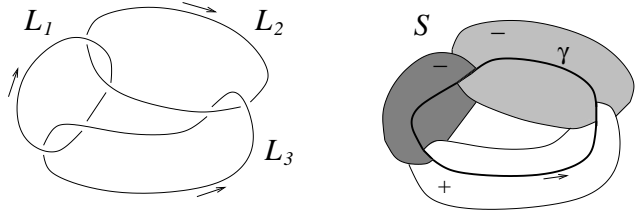,height=2.8cm}
\end{center}
\label{fig:3-color}
\caption{The $3$-colored link of Example \ref{ex:3-color}.}
\end{figure}
\end{ex}

For colored links with $\mu\ge 4$ colors, the presentation of Theorem \ref{thm:module} has deficiency $\genfrac(){0cm}{1}{\mu-1}3$, so
the corresponding presentation matrix is not square. This is not a surprise. Indeed, Crowell and Strauss \cite{C-S}
proved that if an ordered link has $\mu\ge 4$ components and if $\Delta_L\neq 0$, then its Alexander
module does not admit any square presentation matrix. Their proof applies to colored links as well.
Therefore, it is not possible to get a presentation matrix of $H_1(\widetilde{X})$ using the Seifert matrices $A^\eps$ if
$\mu>3$. Nevertheless, it is possible to compute the Alexander invariants up to some indeterminacy.
More precisely, let $\Lambda_\mu'$ denote the localization of the ring $\Lambda_\mu$ with respect to the multiplicative system generated by
$\{t_i-1\}_{1\le i\le\mu}$.

\begin{cor}\label{cor:polynomial}
Let $L$ be a $\mu$-colored link. Consider a C-complex $S$ for $L$ such that $S_i$ is connected for all $i$ and
$S_i\cap S_j\neq\emptyset$ for all $i\neq j$. Then the corresponding matrix
\[
A(t_1,\dots,t_\mu)=\sum_\eps\eps_1\cdots\eps_\mu\,t_1^{\frac{1-\eps_1}{2}}\cdots t_\mu^{\frac{1-\eps_\mu}{2}}A^\eps
\]
is a presentation matrix of the $\Lambda'_\mu$-module $H_1(\widetilde X)\otimes_{\Lambda_\mu}\Lambda_\mu'$.
In particular, for all $r$, there are non-negative integers $m_i$ such that the following equality holds in $\Lambda_\mu$:
\[
\prod_{i=1}^\mu(1-t_i)^{m_i}\,\Delta_r(L)\;\dot{=}\;\Delta_r(A(t_1,\dots,t_\mu)).
\] 
\end{cor}
\begin{proof}
Let $\mathcal B$ be a basis of $H_1(S)$ adapted to the decomposition given in Lemma \ref{lemma:dec}. Consider the dual basis ${\mathcal B}^\ast$
of $H_1(S^3\setminus S)$ obtained by Alexander duality. Clearly, the matrix of $i^\eps\colon H_1(S)\to H_1(S^3\setminus S)$ with
respect to the bases ${\mathcal B},{\mathcal B}^\ast$ is the transposed matrix of $\alpha^\eps\colon H_1(S)\times H_1(S)\to\Zz$ with respect to
$\mathcal B$, that is, a Seifert matrix $A^{-\eps}$. Consider the basis $\widehat{\mathcal B}={\mathcal B}\cup<\gamma_{ijk}>_{2\le i<j<k\le\mu}$ of
$\widehat{H}$, and the extension $\widehat{i}^\eps\colon\widehat{H}\to H_1(S^3\setminus S)$ of $i^\eps$. The matrix
$\widehat A^{-\eps}$ of $\widehat{i}^\eps$ with respect to the bases $\widehat{\mathcal B},{\mathcal B}^\ast$ is a matrix with
$\beta_1(S)$ rows and $\beta_1(S)+\genfrac(){0cm}{1}{\mu-1}3$ columns. $A^{-\eps}$ is made of the first $\beta_1(S)$ columns of
$\widehat A^{-\eps}$. Set
\[
B=\sum_\eps\eps_1\cdots\eps_\mu\,t_1^{\frac{\eps_1+1}{2}}\cdots t_\mu^{\frac{\eps_\mu+1}{2}}A^{-\eps}\quad\hbox{and}\quad
\widehat B=\sum_\eps\eps_1\cdots\eps_\mu\,t_1^{\frac{\eps_1+1}{2}}\cdots t_\mu^{\frac{\eps_\mu+1}{2}}\widehat A^{-\eps}.
\]
By Theorem \ref{thm:module}, a presentation matrix of the $\Lambda_\mu$-module $H_1(\widetilde X)$ is given by $P=\widehat B D^{-1}$,
where $D$ is a diagonal matrix whose diagonal entries are products of factors of the form $t_i-1$. Therefore, $\widehat B$ is a
presentation matrix of the $\Lambda'_\mu$-module $H_1(\widetilde X)\otimes_{\Lambda_\mu}\Lambda_\mu'$.
Now, let $\ell_{ijk}$ denote the column in $\widehat B$ corresponding to the element $\gamma_{ijk}$ of $\widehat{\mathcal B}$,
for $1\le i<j<k\le\mu$. If $i>1$, then $\gamma_{ijk}=\gamma_{1ij}-\gamma_{1ik}+\gamma_{1jk}$ in $H_1(S)$. Therefore,
$\ell_{ijk}=\ell_{1ij}-\ell_{1ik}+\ell_{1jk}$ for all $2\le i<j<k\le\mu$. So the $\genfrac(){0cm}{1}{\mu-1}3$ last columns of $\widehat B$ are
linear combinations of the other columns. Hence, $B$ is also a presentation matrix of $H_1(\widetilde X)\otimes_{\Lambda_\mu}\Lambda_\mu'$.
Since $B=(-1)^\mu A(t_1,\dots,t_\mu)$, the corollary is proved.
\end{proof}

Using this result, it is possible to give an intrinsic interpretation of the nullity $\eta_L(\omega)$ for $\omega$ with rational
coordinates. Indeed, consider the character $\chi_\omega\colon H_1(X)\to\Cc^*$ which maps a meridian of $L_i$ to $\omega_i$. If $\omega$
has rational coordinates, then the image of $\chi_\omega$ is a finite cyclic group $C_q$. The homology of the corresponding cyclic
covering $X^q\to X$ is a module over $\Zz[C_q]$, and so is $\Cc$ via the inclusion $C_q\subset\Cc^*$. Then,
\[
\eta_L(\omega)=\dim_\Cc H_1(X^q)\otimes_{\Zz[C_q]}\Cc.
\]
We shall give a $4$-dimensional proof of this result in Section \ref{section:4-dim}.

\begin{figure}[Htb]
\begin{center}
\epsfig{figure=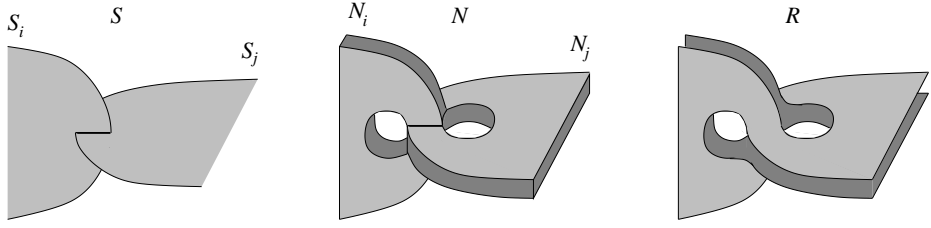,height=2.9cm}
\caption{The C-complex $S$ and the spaces $N$ and $R$ near a clasp.}
\label{fig:FNR}
\end{center}
\end{figure}

\subsection{Proof of Theorem \ref{thm:module}}\label{sub:proof}

As in Section \ref{section:def}, let $X$ denote the exterior of $L$, and let $N_i=S_i\times [-1,1]$ be a closed bicollar neighbourhood of $S_i$. Set
$N=X\cap\bigcup_{i=1}^\mu N_i$, $Y=X\setminus\bigcup_{i=1}^\mu int\,N_i$ and $R=N\cap Y$. (See Figure \ref{fig:FNR} for an illustration
of these spaces near a clasp.) Recall that the regular $\Zz^\mu$-covering $\widetilde X\stackrel{p}{\to}X$ is induced by the epimorphism
$\pi_1(X)\to\Zz^\mu$ given by $\gamma\mapsto(lk(\gamma,L_1),\dots,lk(\gamma,L_\mu))$.
The decomposition $X=N\cup Y$ leads to a Mayer-Vietoris exact sequence of $\Lambda_\mu$-modules
\[
H_1(\widetilde{R})\stackrel{(\varphi,-\psi)}{\longrightarrow}H_1(\widetilde{N})\oplus H_1(\widetilde{Y})\longrightarrow H_1(\widetilde X)
\longrightarrow H_0(\widetilde{R})\stackrel{(\varphi_0,-\psi_0)}{\longrightarrow}H_0(\widetilde{N})\oplus H_0(\widetilde{Y}),
\]
where $\widetilde{R}$ (resp. $\widetilde{N}$, $\widetilde{Y}$) stands for $p^{-1}(R)$ (resp. $p^{-1}(N)$, $p^{-1}(Y)$), and
$\varphi$, $\psi$, $\varphi_0$, $\psi_0$  are inclusion homomorphisms. Note that $R\cap S=\emptyset$. Therefore, given any loop
$\gamma$ in $R$, $lk(\gamma,L_i)=\gamma\cdot S_i=0$ for $1\le i\le\mu$. Hence, $\widetilde{R}\stackrel{p}{\to}R$ is the trivial
$\Zz^\mu$-covering and $H_\ast(\widetilde{R})$ is isomorphic to $H_\ast(R)\otimes\Lambda$, where $\otimes=\otimes_\Zz$ and
$\Lambda=\Lambda_\mu$. For the same reason, we can identify $H_\ast(\widetilde{Y})$ with $H_\ast(Y)\otimes\Lambda$.

We now claim that the homomorphism $(\varphi_0,-\psi_0)$ is injective. Indeed, $S$ is connected, so $R$ is connected as well provided
$\mu>1$. In this case, $\psi_0=id_{\Lambda}$ and $(\varphi_0,-\psi_0)$ is injective. If $\mu=1$, then
$H_0(R)\otimes\Lambda=\Lambda\oplus\Lambda$ and $H_0(\widetilde{N})=\Lambda=H_0(Y)\otimes\Lambda$. The matrix of
$(\varphi_0,-\psi_0)$ with respect to some well-chosen bases is equal to
$\begin{pmatrix}\phantom{-}1&\phantom{-}t\cr -1&-1\end{pmatrix}$,
so the homomorphism is injective. By this claim, we get the exact sequence
\[
H_1(R)\otimes\Lambda\stackrel{(\varphi,-\psi)}{\longrightarrow}
H_1(\widetilde{N})\oplus (H_1(Y)\otimes\Lambda)\longrightarrow H_1(\widetilde X)\to 0. 
\]
Let us assume momentarily that the inclusion homomorphism $\varphi$ is onto. In this case, the following sequence is exact
\[
\ker\varphi\stackrel{\psi}{\longrightarrow}H_1(Y)\otimes\Lambda\longrightarrow H_1(\widetilde X)\longrightarrow 0.
\]
Since $H_1(Y)=H_1(S^3\setminus S)$, we are left with the computation of the homomorphism
$\varphi\colon H_1(R)\otimes\Lambda\to H_1(\widetilde{N})$. We shall divide this tedious computation into several steps.
\medskip

\noindent\parbox{0.6\textwidth}{\noindent{\em Computation of $H_1(R)$.} Let $v_i$ be an interior point of
$(S_i\cap X)\setminus\bigcap_{j\neq i}S_i\cap S_j$. For
$\eps_i=\pm 1$, set $S_i^{\eps_i}=S_i\times\{\eps_i\}\subset S_i\times[-1,1]=N_i$, and let $v_i^{\eps_i}\in R$ denote the point
$v_i\times\{\eps_i\}\in N_i$. Fix a clasp in $S_i\cap S_j$ with $i<j$ and two signs $\eps_i$, $\eps_j$. Consider an oriented edge in $R$
joining $v_i^{\eps_i}$ to $v_j^{\eps_j}$, and passing near this single clasp (see the illustration opposite).}
\hfill\parbox{0.35\textwidth}{\epsfig{figure=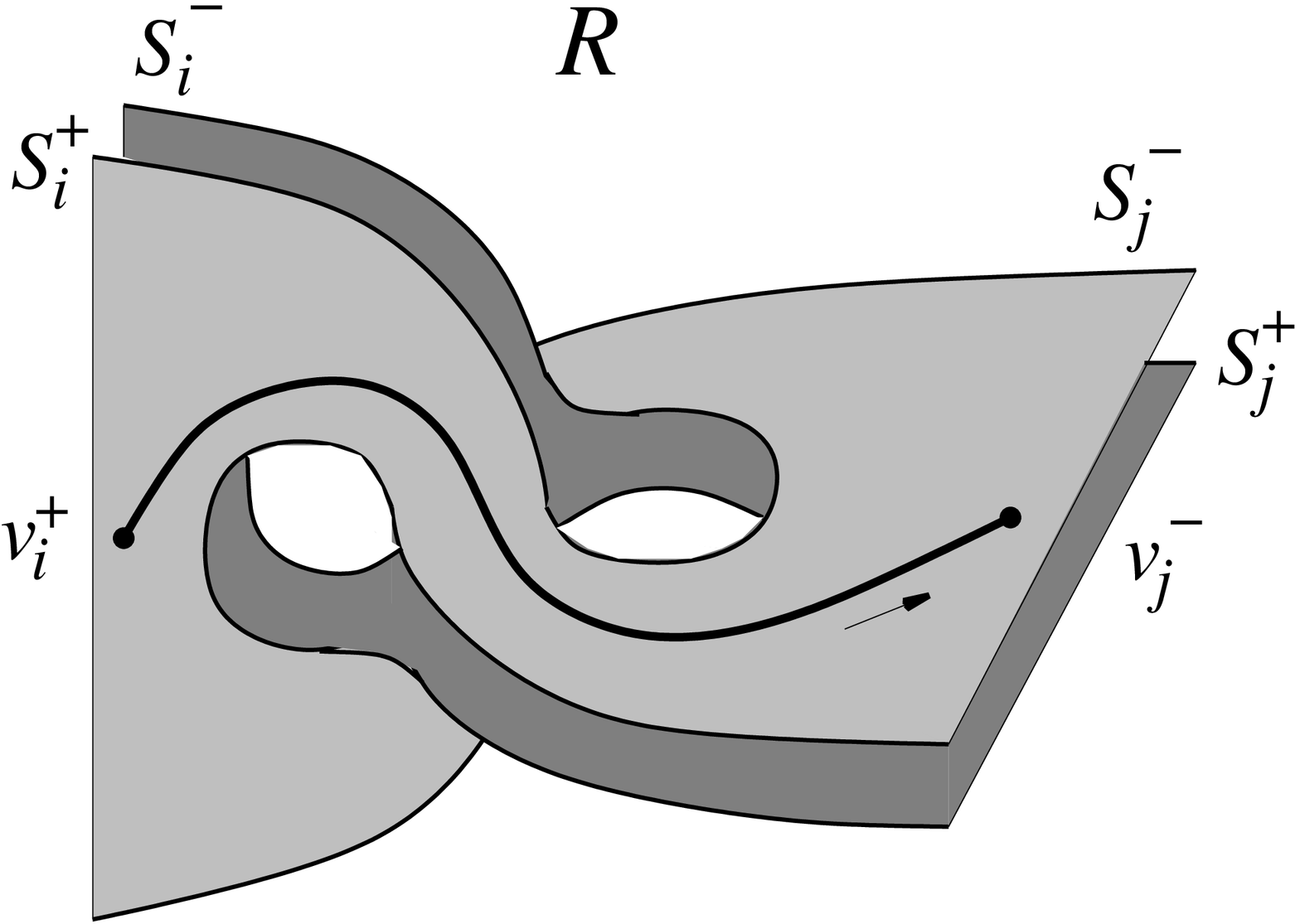,height=3cm}}
This leads to a collection
$\mathcal{E}^{\eps_i\eps_j}_{ij}$ of $c(i,j)$ oriented edges, where $c(i,j)$ is the number of clasps in $S_i\cap S_j$.
Let $K_{ij}^{\eps_i\eps_j}\subset R$ denote the graph given by the union of these edges. Finally, let $K_\mu^{\pm}$ denote the graph
with vertices $\{v_i^{\eps_i}\,;\,1\le i\le\mu,\;\eps_i=\pm 1\}$ and edges
$\{e_{ij}^{\eps_i\eps_i}\,;\,1\le i<j\le\mu,\;\eps_i,\eps_j=\pm 1\}$, where $e_{ij}^{\eps_i\eps_i}$ is one distinguished element
of $\mathcal{E}^{\eps_i\eps_j}_{ij}$. As in the proof of Lemma \ref{lemma:dec}, observe that $R$ can be constructed as follows.
Consider the graph $K_\mu^{\pm}$. Add the graphs $K_{ij}^{--}$ one by one, for $1\le i<j\le\mu$, then the graphs $K_{ij}^{-+}$,
$K_{ij}^{+-}$ and $K_{ij}^{++}$. Finally, paste $S_i^-$ and $S_i^+$ for $1\le i\le\mu$. At each step, the pasting is done along
a contractible space. Therefore,
\begin{eqnarray*}
H_1(R)&=&\bigoplus_i\left(H_1(S_i^-)\oplus H_1(S_i^+)\right)\\
&\oplus&\bigoplus_{i<j}\left(H_1(K_{ij}^{--})\oplus H_1(K_{ij}^{-+})\oplus H_1(K_{ij}^{+-})\oplus H_1(K_{ij}^{++})\right)\\
&\oplus&H_1(K_\mu^{\pm}).
\end{eqnarray*}
\medskip

\noindent\parbox{0.61\textwidth}{\noindent{\em Computation of $H_1(\widetilde{N})$.}
Consider the complete oriented graph $K_\mu$ defined above Lemma \ref{lemma:dec}. Add a vertex $v_{ij}$ in the interior of the
edge joining $v_i$ with $v_j$, and paste two oriented loops $x_i$ and $x_j$ based at $v_{ij}$.
The resulting graph, denoted by $\Gamma_\mu$, is naturally embedded in $S\cap X\subset N$ as illustrated opposite.}
\hfill\parbox{0.35\textwidth}{\epsfig{figure=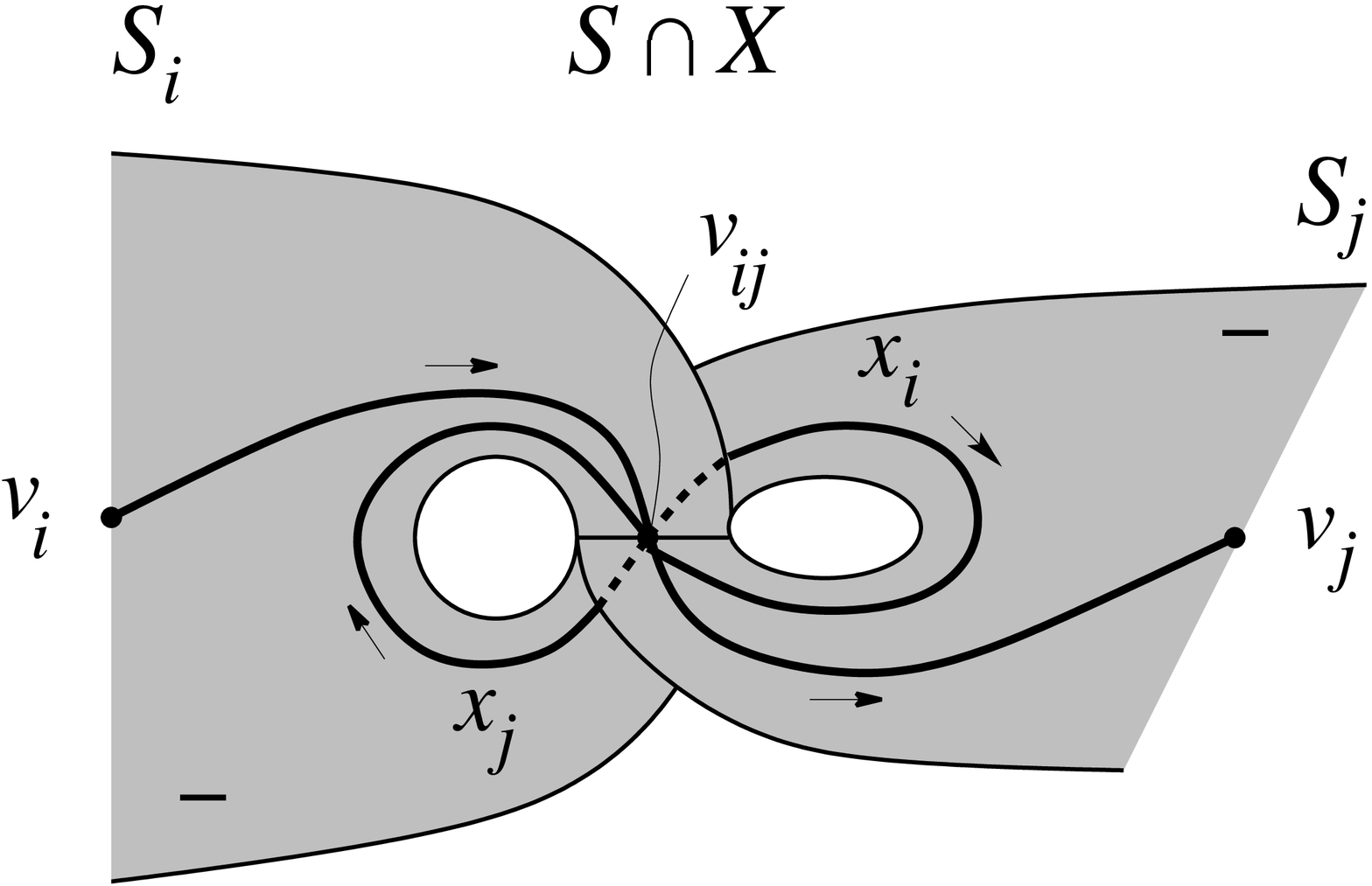,height=2.9cm}}
Similarly, let $\Gamma_{ij}\subset S_i\cap S_j\cap X $ be the
oriented graph obtained from $K_{ij}\subset S_i\cap S_j$ by adding a vertex in the interior of each edge of $K_{ij}$, and two
oriented loops at this vertex. For $i<j$, consider the loop at $v_{ij}$ given by $x_ix_jx_i^{-1}x_j^{-1}$. This loop lifts to a loop
$\gamma_{ij}$ in $\widetilde\Gamma_{ij}:=p^{-1}(\Gamma_{ij})$ and to a loop $\gamma'_{ij}$ in $\widetilde\Gamma_\mu:=p^{-1}(\Gamma_\mu)$.
We claim that $H_1(\widetilde{N})$ splits as follows
\[
H_1(\widetilde{N})=\bigoplus_i\left(H_1(S_i)\otimes\Lambda\right)\oplus
\Big(\bigoplus_{i<j} H_1(\widetilde\Gamma_{ij})\oplus H_1(\widetilde\Gamma_\mu)\Big)\Big/(\gamma_{ij}=\gamma'_{ij})_{i<j}.
\]
\noindent\parbox{0.6\textwidth}{Indeed, $N$ has the homotopy type of $S\cap X$, which can be contructed as
follows. Consider the graph $\Gamma_\mu$. Add the graphs $\Gamma_{ij}$ for $1\le i<j\le\mu$, and then the surfaces
$S'_i$ for $1\le i\le\mu$, where $S'_i$ is obtained from $S_i\cap X$ by removing an open neighbourhood of the
clasps as illustrated opposite.}
\hfill\parbox{0.35\textwidth}{\epsfig{figure=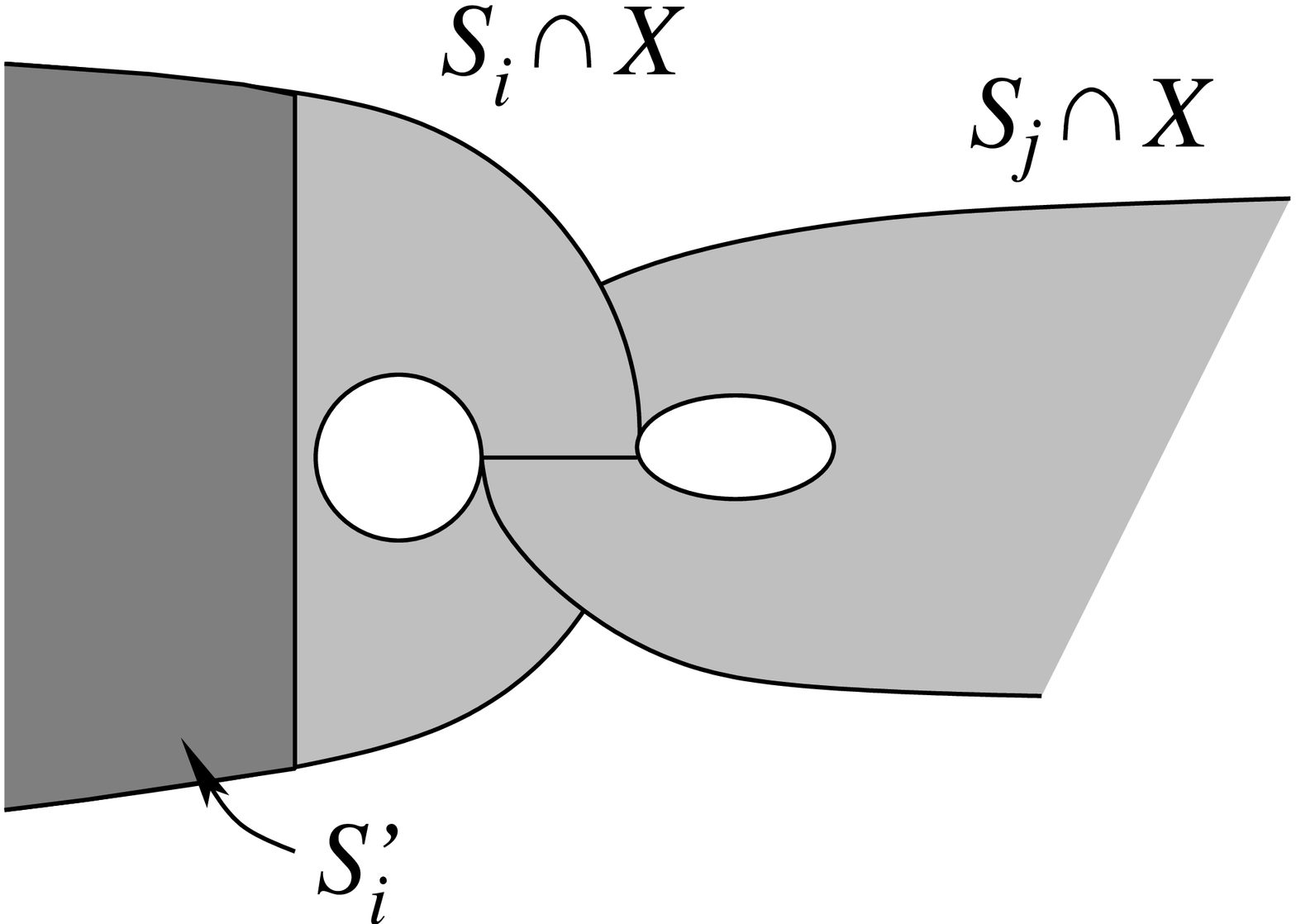,height=3cm}}
Hence, $\widetilde{N}$ has the homotopy type of
$\widetilde\Gamma_\mu\cup\bigsqcup_{i<j}\widetilde\Gamma_{ij}\cup\bigsqcup_{i}p^{-1}(S'_i)$.
The first union $\Gamma_\mu\cup\Gamma_{12}$ is made along $\Gamma_\mu\cap\Gamma_{12}$ which has the homotopy type of the wedge of the
circles $x_1$ and $x_2$. This yields a Mayer-Vietoris exact sequence of
$\Lambda$-modules
\[
0\to\Lambda\gamma_{12}\to H_1(\widetilde\Gamma_\mu)\oplus H_1(\widetilde\Gamma_{12})\to
H_1(\widetilde\Gamma_\mu\cup\widetilde\Gamma_{12})\to
\Lambda/\mathcal{I}_{12}\to\Lambda/\mathcal{I}\oplus\Lambda/\mathcal{I}_{12},
\]
where $\mathcal I$ is the augmentation ideal $(t_1-1,\dots,t_\mu-1)$ and $\mathcal{I}_{12}$ the ideal $(t_1-1,t_2-1)$ of $\Lambda$.
Hence $H_1(\widetilde\Gamma_\mu\cup\widetilde\Gamma_{12})
=\left(H_1(\widetilde\Gamma_\mu)\oplus H_1(\widetilde\Gamma_{12})\right)/(\gamma_{12}=\gamma'_{12})$. Using this argument inductively,
we get
\[
H_1(\widetilde\Gamma_\mu\cup\bigsqcup_{i<j}\widetilde\Gamma_{ij})=\Big(\bigoplus_{i<j} H_1(\widetilde\Gamma_{ij})\oplus
H_1(\widetilde\Gamma_\mu)\Big)\Big/(\gamma_{ij}=\gamma'_{ij})_{i<j}.
\]
The pasting of $S'_i$ is performed along a contractible
space. Since $p^{-1}(S'_i)=S'_i\times\Zz^\mu$ and $S'_i$ has the homotopy type of $S_i$, the
claim follows from yet another application of the Mayer-Vietoris sequence.
\medskip

Using the computation of $H_1(R)$ and $H_1(\widetilde{N})$ above, we see that the inclusion homomorphism
$\varphi\colon H_1(R)\otimes\Lambda\to H_1(\widetilde{N})$ splits into
a direct sum $\varphi=\bigoplus_i\varphi_i\oplus\bigoplus_{i<j}\varphi_{ij}\oplus\varphi_K$. We shall now check that each homomorphism
in this sum is onto, and then compute its kernel. Each kernel will translate into a family of
relations via the inclusion homomorphism $\psi\colon H_1(R)\otimes\Lambda\to H_1(S^3\setminus S)\otimes\Lambda$. Note that
up to multiplication by a unit of $\Lambda$, the latter homomorphism is given by $\psi=\psi'\otimes id_\Lambda$ where
$\psi'\colon H_1(R)\to H_1(S^3\setminus S)$ is the inclusion homomorphism.
\medskip

\noindent{\em Computation of $\varphi_i$.} The homomorphism
\[
\varphi_i\colon\left(H_1(S_i^-)\oplus H_1(S_i^+)\right)\otimes\Lambda\to H_1(S_i)\otimes\Lambda
\]
is given by $\varphi_i(x^-,y^+)=t_ix+y$ for $x,y\in H_1(S_i)$. Therefore, the homomorphism $\varphi_i$ is onto,
and its kernel is generated
by $\{(-x^-,t_ix^+)\,;\,x\in H_1(S_i)\}$ as a $\Lambda$-submodule of $\left(H_1(S_i^-)\oplus H_1(S_i^+)\right)\otimes\Lambda$.
Given $x\in H_1(S_i)\subset H_1(S)$, we have $\psi(x^+)=i^+(x)$ and $\psi(x^-)=i^-(x)$, where $i^+$ (resp.
$i^-$) is the homomorphism $i^\eps$ with $\eps$ any sequence of $\pm 1$'s such that $\eps_i=+1$ (resp. $\eps_i=-1$).
Therefore, $\psi(\ker(\varphi_i))=(t_ii^+-i^-)(H_1(S_i))$. The restriction to $H_1(S_i)$ of the homomorphism $i^\eps$ only depends on
$\eps_i$. Hence, the restriction to $H_1(S_i)$ of the homomorphism
$\alpha=\sum_\eps\eps_1\cdots\eps_\mu\,t_1^{\frac{\eps_1+1}{2}}\cdots t_\mu^{\frac{\eps_\mu+1}{2}}i^\eps$ is equal to
$(t_ii^+-i^-)\prod_{n\neq i}(t_n-1)$. This gives the first set of relations stated in the theorem.
\medskip

\noindent{\em Computation of $\varphi_{ij}$.} Let us now consider the inclusion homomorphism
\[
\varphi_{ij}\colon\left(H_1(K_{ij}^{--})\oplus H_1(K_{ij}^{-+})\oplus H_1(K_{ij}^{+-})\oplus H_1(K_{ij}^{++})\right)\otimes\Lambda\to
H_1(\widetilde\Gamma_{ij})/(\gamma_{ij}).
\]
We need to compute the homology of $\widetilde\Gamma_{ij}$. Number the clasps of $S_i\cap S_j$ from $1$ to $c(i,j)$. For
$1\le\ell<c(i,j)$, consider the loop in $\Gamma_{ij}$ starting at $v_i$, going to the vertex corresponding to the clasp number $\ell$,
then around $x_i$, back to $v_i$, to the vertex corresponding to the clasp number $\ell+1$, around $x_i^{-1}$ and back to $v_i$. This
loop lifts to a loop $\delta_i^\ell$ in $\widetilde\Gamma_{ij}$. One checks by induction on $c(i,j)\ge 1$ that
\[
H_1(\widetilde\Gamma_{ij})=
(H_1(K_{ij})\otimes\Lambda)\oplus\bigoplus_{\ell=1}^{c(i,j)-1}(\Lambda\delta_i^\ell\oplus\Lambda\delta_j^\ell)\oplus\Lambda\gamma_{ij}.
\]
Recall the basis $\{\beta_{ij}^\ell\}_\ell$ of $H_1(K_{ij})$ given in Lemma \ref{lemma:dec}. This yields a basis
$\{(\beta_{ij}^\ell)^{\eps_i\eps_j}\}_\ell$ of $H_1(K_{ij}^{\eps_i\eps_j})$ for $\eps_i,\eps_j=\pm 1$. The homomorphism $\varphi_{ij}$ is
given by
\[
\begin{array}{ll}
(\beta_{ij}^\ell)^{--}\mapsto t_it_j\beta_{ij}^\ell,&
(\beta_{ij}^\ell)^{++}\,\mapsto\,(t_i+t_j-1)\beta_{ij}^\ell+\delta_i^\ell-\delta_j^\ell,\\
(\beta_{ij}^\ell)^{-+}\mapsto t_i(t_j\beta_{ij}^\ell-\delta_j^\ell),&
(\beta_{ij}^\ell)^{+-}\,\mapsto\,t_j(t_i\beta_{ij}^\ell+\delta_i^\ell).
\end{array}
\]
Therefore, $\varphi_{ij}$ is onto and its kernel is generated by
\[
\{(t_it_j\beta^{++},-t_i\beta^{+-},-t_j\beta^{-+},\beta^{--})\,;\,\beta\in H_1(K_{ij})\}.
\]
For $\beta\in H_1(K_{ij})\subset H_1(S)$, $\psi(\beta^{++})=i^{++}(\beta)$ where $i^{++}=i^\eps$ with $\eps$ any sequence of $\pm 1$'s
such that $\eps_i=\eps_j=+1$. Similarly, $\psi(\beta^{+-})=i^{+-}(\beta)$, $\psi(\beta^{-+})=i^{-+}(\beta)$ and
$\psi(\beta^{--})=i^{--}(\beta)$. Hence, the image of $\ker(\varphi_{ij})$ under $\psi$ is the image of $H_1(K_{ij})$ under the
homomorphism $t_it_ji^{++}-t_ii^{+-}-t_ji^{-+}+i^{--}$. Since the restriction to $H_1(K_{ij})$ of the homomorphism $\alpha$ is equal to
$(t_it_ji^{++}-t_ii^{+-}-t_ji^{-+}+i^{--})\prod_{n\neq i,j}(t_n-1)$, we get the second set of relations stated in the theorem. 
\medskip

\noindent{\em Computation of $\varphi_K$.} Finally, let us deal (a little faster) with the inclusion homomorphism
\[
\varphi_K\colon H_1(K_\mu^\pm)\otimes\Lambda\to H_1(\widetilde\Gamma_\mu).
\] 
By an Euler characteristic argument, the rank of $H_1(K_\mu^\pm)$ is equal to $4{\genfrac(){0cm}{1}\mu 2}-2\mu+1$. Moreover, a recursive use of
the Mayer-Vietoris exact sequence leads to the following fact: it is possible to present the $\Lambda$-module $H_1(\widetilde\Gamma_\mu)$ 
with $4{\genfrac(){0cm}{1}\mu 2}-2\mu+1$ generators and $\genfrac(){0cm}{1}\mu 3$ relations. (The $\genfrac(){0cm}{1}\mu 3$ relations correspond to all the
possible $3$-dimensional `cubes' in $\Zz^\mu$.) One then checks that $\varphi_K$ is onto. Therefore, each relation in
$H_1(\widetilde\Gamma_\mu)$ yields an element in $\ker(\varphi_K)$. One shows that the relation coming from the cube with coordinates
$1\le i<j<k\le\mu$ corresponds to the element
\[
t_it_jt_k\gamma_{ijk}^{+++}-t_it_j\gamma_{ijk}^{++-}-t_it_k\gamma_{ijk}^{+-+}-t_jt_k\gamma_{ijk}^{-++}+t_i\gamma_{ijk}^{+--}
+t_j\gamma_{ijk}^{-+-}+t_k\gamma_{ijk}^{--+}-\gamma_{ijk}^{---}
\]
in $\ker(\varphi_K)$, where $\gamma_{ijk}^{\eps_i\eps_j\eps_k}$ is the $1$-cycle in $H_1(K_\mu^\pm)$ joining $v_i^{\eps_i}$,
$v_j^{\eps_j}$ and $v_k^{\eps_k}$ in this order. The image under $\psi$ of this element is equal to
$\prod_{n\neq i,j,k}(t_n-1)^{-1}\alpha(\gamma_{ijk})$, completing the proof of the theorem.
\qed

\section{Piecewise continuity of the signature and nullity}\label{section:PC}

Let $L$ be a $\mu$-colored link. The signature and nullity of $L$ can be understood as functions
\[
\sigma_L,\eta_L\colon T^\mu_\ast\colon\longrightarrow\Zz,
\]
where $T^\mu_\ast=(S^1\setminus\{1\})^\mu\subset\Cc^\mu$. In this subsection, we prove the following `piecewise continuity' result.

\begin{thm}\label{thm:semi}
Let $E_r(L)\subset\Lambda_\mu$ be the $r^{th}$ Alexander ideal of $L$, and set
\[
\Sigma_r=\{\omega\in T^\mu_\ast\;\mid\;p(\omega)=0 \hbox{ for all } p\in E_{r-1}(L)\}. 
\]
This yields a finite sequence $T^\mu_\ast=\Sigma_0\supset\Sigma_1\supset\dots\supset\Sigma_{\ell-1}\supset\Sigma_\ell=\emptyset$
such that, for all $r$, $\eta_L$ is equal to $r$ on $\Sigma_r\setminus\Sigma_{r+1}$, and $\sigma_L$ is constant on the
connected components of $\Sigma_r\setminus\Sigma_{r+1}$.
\end{thm}

\begin{proof}
Let $S=S_1\cup\dots\cup S_\mu$ be a C-complex for $L$. Since $\sigma_L$ and $\eta_L$ remain unchanged by transformations $T1$ and $T2$
(recall Lemma \ref{lemma:Sequiv} and the proof of Theorem \ref{thm:inv}), it may be assumed that $S_i$ is connected for all $i$, and
$S_i\cap S_j\neq\emptyset$ for all $i\neq j$. In this case, Corollary \ref{cor:polynomial} implies that
\[
E_r(L)\otimes_{\Lambda_\mu}\Lambda_\mu'=E_r(A(t))\otimes_{\Lambda_\mu}\Lambda_\mu',
\]
where $A(t)=A(t_1,\dots,t_\mu)=\sum_\eps\eps_1\cdots\eps_\mu\,t_1^{\frac{1-\eps_1}{2}}\cdots t_\mu^{\frac{1-\eps_\mu}{2}}A^\eps$.
Hence, a given $\omega\in T^\mu_\ast$ belongs to $\Sigma_r$ if and only if $p(\omega)=0$ for all $p\in E_{r-1}(A(t))$.
Now, recall that $\sigma_L(\omega)$ and $\eta_L(\omega)$ are defined as the signature and nullity of the Hermitian matrix
$H(\omega)=\prod_{i=1}^\mu(1-\overline{\omega}_i)A(\omega)$. By this equality, $p(\omega)=0$ for all $p\in E_{r-1}(A(t))$
if and only if all the $(n-r+1)\times(n-r+1)$ minors of $H(\omega)$ vanish, where $n$ is the dimension of $H(\omega)$.
This occurs if and only if the
nullity of $H(\omega)$ is greater or equal to $r$. Therefore, $\omega\in T^\mu_\ast$ belongs to $\Sigma_r$ if and only if
$\eta_L(\omega)\ge r$. This gives the first part of the theorem. The second
part is a consequence of the first point and of the following claim.\\
{\em Claim:\/} Let $x\mapsto H(x)$ be a continuous path in the space of $n$-dimensional Hermitian matrices. If the nullity of $H(x)$
is constant along the path, then the signature of $H(x)$ is constant as well.\\
Indeed, consider the continuous path $p\colon[0,1]\to\Cc[\lambda]$, where $p(x)$ is the characteristic polynomial of $H(x)$.
Since $H(x)$ is Hermitian, the roots of $p(x)$ are real. Furthermore, they depend continuously on the coefficients of $p(x)$.
In other words, the eigenvalues of $H(x)$ are real continuous functions of $x$.
The nullity of $H(x)$ counts the number of these eigenvalues which are zero. Hence, if this
number is constant, the sign of the eigenvalues cannot change. Therefore, the signature of $H(x)$ is constant.
\end{proof}

As a corollary, we obtain the following result that extends a well known property of the Levine-Tristram signature.
Note that in the case $lk(L_i,L_j)=0$, it is due to the second author (see Section \ref{section:4-dim} and
\cite[theorems 3.7 and 3.8]{Flo}).
\begin{cor}\label{cor:semiC}
Let $L$ be a $\mu$-colored link, and let $C$ denote the complement in $T^\mu_\ast$ of the zeroes of its Alexander
polynomial. Then, $\sigma_L$ is constant on the connected components of $C$ and $\eta_L$ 
vanishes on $C$.
\end{cor}
\begin{proof}
By Theorem \ref{thm:semi}, we just need to check that $C=T^\mu_\ast\setminus\Sigma_1$, that is: if $\omega\in T^\mu_\ast$, then
$p(\omega)=0$ for all $p\in E_0(L)$ if and only if $\Delta_L(\omega)=0$.
Corollary \ref{cor:polynomial} implies that $E_0(L)\otimes_{\Lambda_\mu}\Lambda_\mu'=E_0(A(t))\otimes_{\Lambda_\mu}\Lambda_\mu'$.
Since $A(t)$ is a square matrix, $E_0(A(t))$ is the principal ideal generated by
$\det(A)\,\dot{=}\,\prod_i(1-t_i)^{m_i}\cdot\Delta_L$. This proves the corollary.
\end{proof}

Before concluding this section, let us look back at the $2$-colored link $L$ given in Example \ref{ex:lk2}.
By Corollary \ref{cor:mod2}, a presentation matrix of its Alexander
module is given by $(t_1t_2+1)$. Therefore, $E_0(L)=(t_1t_2+1)$ and $E_r(L)=\Lambda_2$ for $r\ge 1$, leading to
\[
\Sigma_r=\begin{cases}T^2_\ast&\text{if $r\le 0$;}\cr
	\{(\omega_1,\omega_2)\in T^2_\ast\,\mid\,\omega_1\omega_2+1=0\}&\text{if $r=1$;}\cr
	\emptyset&\text{if $r\ge 2$.}
	\end{cases}
\]
By Theorem \ref{thm:semi}, $\sigma_L$ is constant on the connected components of $T^2_\ast\setminus\Sigma_1$ and of $\Sigma_1$.
Furthermore, $\eta_L$ is equal to $0$ on $T^2_\ast\setminus\Sigma_1$ and equal to $1$ on $\Sigma_1$. This coincides with the
computations made in Example \ref{ex:lk2}.

\section{On the computation by local moves}\label{section:local}

We now present and generalize an idea of J. Conway (see \cite[\S 7.10]{Coop} and \cite[Lemma 3.1]{Ore2}) leading to a purely
combinatorial computation of $\sigma_L$ for many colored links.
Recall that for any $\mu$-colored link $L$, there exists a well-defined invariant
$\nabla_L(t_1,\dots,t_\mu)\in\Zz(t_1,\dots,t_\mu)$ called the {\em Conway potential function of $L$\/},
which satisfies:
\[
\nabla_L(t_1,\dots,t_\mu)\,\;\dot{=}\,\begin{cases}\frac{1}{t_1-t_1^{-1}}\,\Delta_L(t_1^2)&\text{if $\mu=1$;}\cr  
            	\Delta_L(t_1^2,\dots,t_\mu^2)&\text{if $\mu>1$.}\end{cases}
\]
This normalization of the Alexander polynomial was first introduced by J. Conway \cite{Con}, and formally defined
by R. Hartley \cite{Har}. Note that it is possible to compute this invariant from a colored link diagram using only
combinatorial methods (see \cite{Mur}).

Given a complex number $z=e^{i\theta}$ with $0<\theta<2\pi$, we shall denote by $z^{1/2}$ the
complex number $e^{i\theta/2}$. If $\omega=(\omega_1,\dots,\omega_\mu)\in T^\mu_\ast$, then
$\omega^{1/2}=(\omega^{1/2}_1,\dots,\omega^{1/2}_\mu)$. Also, we shall denote by $sgn(\lambda)$ the sign of the
real number $\lambda$.

\begin{figure}[Hb]
\begin{center}
\epsfig{figure=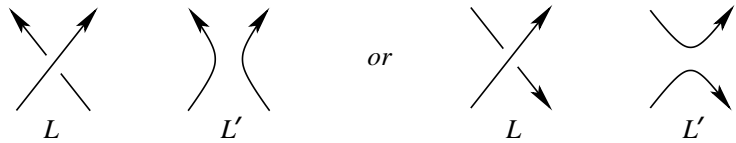,height=1.8cm}
\caption{The local move in part a) of Theorem \ref{thm:local}. The two strands are assumed to have the same color.}
\label{fig:M1}
\end{center}
\end{figure}

\begin{figure}[Ht]
\begin{center}
\epsfig{figure=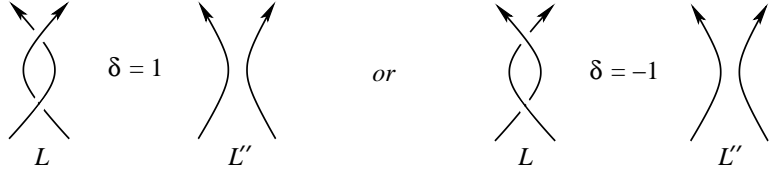,height=2.2cm}
\caption{The local moves in part b) of Theorem \ref{thm:local}. The colors of the two strands are assumed to be distinct.}
\label{fig:M2}
\end{center}
\end{figure}

\begin{thm}\label{thm:local}
$a)$ Let $L$ and $L'$ be two colored links given by diagrams related by a single change as illustrated in Figure \ref{fig:M1}.
For any $\omega\in T^\mu_\ast$ such that $\nabla_{L'}(\omega^{1/2})\neq 0$,
$i\;\frac{\nabla_L(\omega^{1/2})}{\nabla_{L'}(\omega^{1/2})}$ is a real number and
\[
\sigma_L(\omega)=\sigma_{L'}(\omega)+ sgn\left(i\;\frac{\nabla_L(\omega^{1/2})}{\nabla_{L'}(\omega^{1/2})}\right).
\]
Conversely, for any $\omega\in T^\mu_\ast$ such that $\nabla_L(\omega^{1/2})\neq 0$,
\[
\sigma_{L'}(\omega)=\sigma_L(\omega)+sgn\left(i\;\frac{\nabla_{L'}(\omega^{1/2})}{\nabla_L(\omega^{1/2})}\right).
\]
$b)$ Consider colored links $L$ and $L''$ which differ by one of the local moves described in Figure \ref{fig:M2}.
If $\omega\in T^\mu_\ast$ satisfies $\nabla_{L''}(\omega^{1/2})\neq 0$, then
$\frac{\nabla_{L}(\omega^{1/2})}{\nabla_{L''}(\omega^{1/2})}$ is a real number and
\[
\sigma_L(\omega)=\sigma_{L''}(\omega)+\delta\cdot sgn\left(\frac{\nabla_L(\omega^{1/2})}{\nabla_{L''}(\omega^{1/2})}\right),
\]
with $\delta=\pm 1$ as in Figure \ref{fig:M2}. Conversely, if $\omega\in T^\mu_\ast$ satisfies $\nabla_L(\omega^{1/2})\neq 0$, then
\[
\sigma_{L''}(\omega)=\sigma_{L}(\omega)-\delta\cdot sgn\left(\frac{\nabla_{L''}(\omega^{1/2})}{\nabla_{L}(\omega^{1/2})}\right).
\]
\end{thm}
 
Given $L$ a link, let $\Delta_L(t)$ denote its Alexander-Conway polynomial
$\Delta_L(t)=(t^{1/2}-t^{-1/2})\nabla_L(t^{1/2})$. The theorem above implies the following result for the
Levine-Tristram signature.

\begin{cor}\label{cor:1col}
Let $L$ and $L'$ be two links related by a single change as illustrated in Figure \ref{fig:M1}.
For any $\omega\in S^1\setminus\{1\}$ such that $\Delta_{L'}(\omega)\neq 0$,
$i\;\frac{\Delta_{L'}(\omega)}{\Delta_L(\omega)}\in\Rr$ and
\[
\sigma_L(\omega)=\sigma_{L'}(\omega)+sgn\left(i\;\frac{\Delta_L(\omega)}{\Delta_{L'}(\omega)}\right).
\]
Furthermore, for any $\omega\in T^\mu_\ast$ such that $\Delta_L(\omega)\neq 0$,
\[
\sigma_{L'}(\omega)=\sigma_L(\omega)+sgn\left(i\;\frac{\Delta_{L'}(\omega)}{\Delta_L(\omega)}\right).\qed
\]
\end{cor}

Let us postpone the proof of Theorem \ref{thm:local} to the end of the section. We shall first discuss to what extend
this result leads to an algorithm for the computation of the signature of a colored link. Let us start with the case of a knot $K$.
Corollary \ref{cor:1col} provides the following combinatorial algorithm for the computation of $\sigma_K(\omega)$ for all
but a finite number of $\omega$ in $S^1$. Consider a diagram for $K$ with $n$ double points and no `nugatory crossings'
(that is, a diagram that remains connected after the move of Figure \ref{fig:M1}).
Apply the local move of Figure \ref{fig:M1} to any crossing. We get a
2-component link $L'$ given by a connected diagram. Apply the same local move to any crossing between the two components of $L'$.
(Such a crossing exists since the diagram is connected.) This yields a knot $K''$ given by a diagram with $n-2$ crossings. Corollary
\ref{cor:1col} gives the equality
$$
\sigma_K(\omega)=\sigma_{K''}(\omega)-sgn\left(i\;\frac{\Delta_{L'}(\omega)}{\Delta_K(\omega)}\right)
+sgn\left(i\;\frac{\Delta_{L'}(\omega)}{\Delta_{K''}(\omega)}\right)
$$
for all $\omega$ such that $\Delta_{K}(\omega)\Delta_{K''}(\omega)\neq 0$. Since both $K$ and $K''$ are knots, $\Delta_K\Delta_{K''}$
is non-zero. Hence, the relation above holds for all but a finite number of $\omega$ in $S^1$. We are done by induction
on $n$.
\begin{ex}
Let $K$ be the right-hand trefoil knot. Consider the following transformations of $K$.

\begin{figure}[h]
\begin{center}
\epsfig{figure=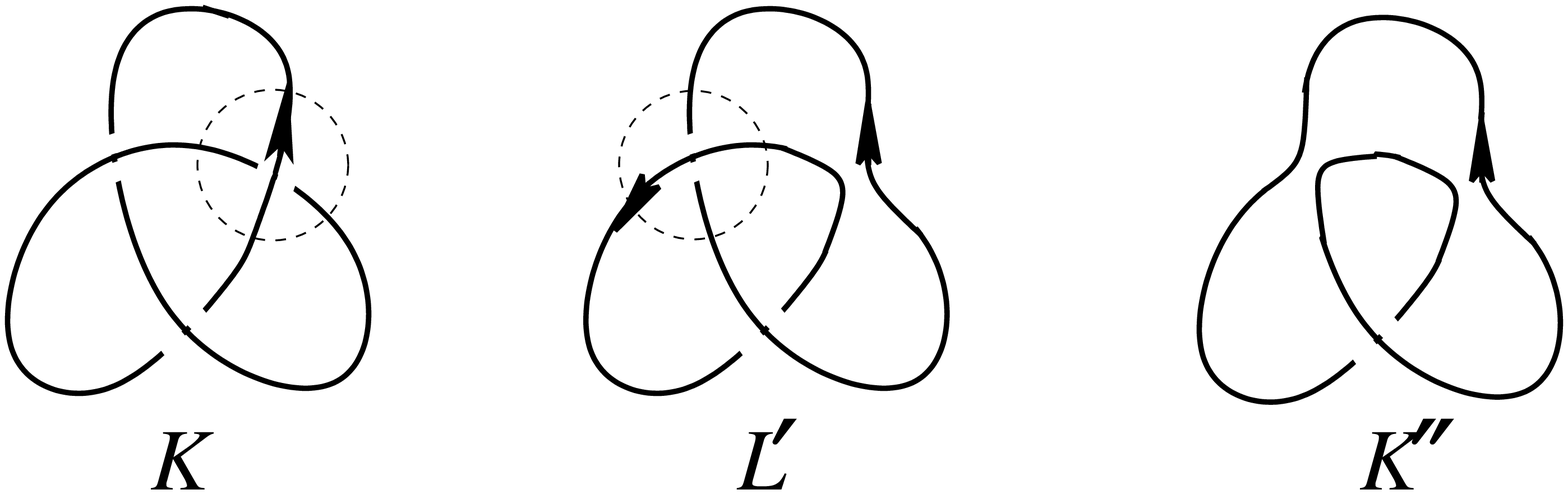,height=2.2cm}
\end{center}
\end{figure}
\noindent One easily computes $\Delta_K(t)=t-1+t^{-1}$. Moreover, $L'$ is the positive Hopf link and $K''$ the trivial knot, so
$\Delta_{L'}(t)=t^{1/2}-t^{-1/2}$ and $\Delta_{K''}(t)=1$. This gives
\begin{eqnarray*}
\sigma_K(\omega)&=&\sigma_{K''}(\omega)-sgn\left(i\,\frac{\omega^{1/2}-\overline{\omega^{1/2}}}{\omega-1+\overline\omega}\right)
+sgn\left(i\left(\omega^{1/2}-\overline{\omega^{1/2}}\right)\right)\\
&=& sgn\left(\frac{1}{\omega-1+\overline\omega}\right)-1.
\end{eqnarray*}
On this particularly simple example, $\Delta_{L'}$ never vanishes on $S^1\setminus\{1\}$. Hence, Corollary \ref{cor:1col}
in fact leads to the formula
\[
\sigma_K(\omega)=sgn(\omega-1+\overline\omega)-1
\]
which holds for any value of $\omega$.
\end{ex}

Now, let $L=K_1\cup\dots\cup K_\mu$ be a $\mu$-component $\mu$-colored link. Theorem \ref{thm:local} b) leads to the following
result: there exists a $\mu$-variable polynomial $\Delta$ such that, for all $\omega=(\omega_1,\dots,\omega_\mu)$ in $T^\mu$ with
$\Delta(\omega)\neq 0$, the computation of $\sigma_L(\omega)$ boils down to the computation of
$\sigma_{K_i}(\omega_i)$ for $i=1,\dots,\mu$. And we just saw that for all $i$, $\sigma_{K_i}(\omega_i)$ can be computed
for all but a finite number of $\omega_i$ in $S^1$. However, for some links $L$, the polynomial $\Delta$ is always zero and the
result above is useless. For example, consider an irreducible boundary link $L$ with $\mu\ge 3$ components. Let $L'$ be any
$\mu$-colored link obtained from $L$ by one of the local moves described in figures \ref{fig:M1} and \ref{fig:M2}. Then, it is easy to
check that both $\nabla_L$ and $\nabla_{L'}$ are identically zero. Therefore, Theorem \ref{thm:local} is of no use for the computation
of $\sigma_L$ in this particular case.

However, the method described above does lead to an algorithm for links that are not `algebraically split'.
Recall that a $\mu$-component link $L$ is said to be {\em algebraically split\/} if there is an ordering of its components
$K_1,\dots,K_\mu$ and an integer $1\le k<\mu$ such that $lk(K_i,K_j)=0$ for all $1\le i\le k<j\le\mu$. Using the Torres formula,
one checks that if $L$ is not algebraically split, then $\nabla_L$ is non-zero. Using this result, one easily sees that if $L$
is not algebraically split, then the polynomial $\Delta$ described above can be chosen to be non-zero.
Therefore, $\sigma_L(\omega)$ can be computed for all $\omega$ in the complement in $T^\mu$ of the zeroes
of some non-zero polynomial. Let us illustrate this with an example.

\begin{ex}
Consider the colored link $L$ given in Example \ref{ex:lk2}. Also, let $L''$ be the $2$-colored positive Hopf link.
Since $\nabla_L(t_1,t_2)=t_1t_2+t^{-1}_1t^{-1}_2$ and $\nabla_{L''}(t_1,t_2)=1$, Theorem \ref{thm:local} b) gives
\[
\sigma_L(\omega_1,\omega_2)=
\sigma_{L''}(\omega_1,\omega_2)+sgn\left(\omega^{1/2}_1\omega^{1/2}_2+\overline{\omega^{1/2}_1\omega^{1/2}_2}\right).
\]
By Theorem \ref{thm:local} b) again (or since there is a contractible C-complex for $L''$), the signature of $L''$ is zero.
Setting $\omega_1=e^{i\theta_1}$ and $\omega_2=e^{i\theta_2}$, we get
\[
\sigma_L(\theta_1,\theta_2)=sgn\left(\cos\left((\theta_1+\theta_2)/2\right)\right).
\]
This coincides with the result given in Example \ref{ex:lk2}.
\end{ex}

The proof of Theorem \ref{thm:local} rests upon several lemmas.

\begin{lemma}\label{lemma:herm}
Let $H$ be a $d$-dimensional Hermitian matrix. Then its signature $\sigma$ and nullity $\eta$ satisfy
$\sigma+\eta\equiv d\pmod{2}$. Furthermore, if $\eta=0$, then
$\sigma\equiv d\pmod{4}$ if and only if $\det(H)>0$. 
\end{lemma}
\begin{proof}
Let $m$ be the number of negative eigenvalues of $H$. We have
\[
\sigma+\eta = (d-m-\eta)-m+\eta = d-2m \equiv d\pmod{2}.
\]
If $\eta=0$, then $\det(H)>0$ if and only if $m$ is even, that is, if and only if $\sigma=d-2m\equiv d\pmod{4}$.
\end{proof}

\begin{lemma}\label{lemma:mod2}
Let $L$ be a colored link with $\nu$ connected components. Then, for all $\omega\in T^\mu_\ast$,
\[
\sigma_L(\omega)+\eta_L(\omega)\equiv \nu + \sum_{i<j}lk(L_i,L_j)+1\pmod{2}.
\]
\end{lemma}
\begin{proof}
Let $S$ be a connected C-complex for $L$. Applying Lemma \ref{lemma:herm} to the matrix $H(\omega)$, we have
$\sigma_L(\omega)+\eta_L(\omega)\equiv\rk H_1(S)\pmod{2}$. An elementary Euler characteristic argument shows that
$\rk H_1(S)\equiv \nu+\#\{\hbox{clasps of $S$}\}+1\pmod{2}$.
Finally, the number of clasps of $S$ clearly has the same parity as $\sum_{i<j}lk(L_i,L_j)$.
\end{proof}

\begin{lemma}\label{lemma:mod4}
Let $L$ be a colored link with $\nu$ connected components. Then, for any $\omega\in T^\mu_\ast$ such
that $\eta_L(\omega)=0$, we have
\[
\sigma_L(\omega)\equiv \nu+\sum_{i<j}lk(L_i,L_j)-sgn(i^{\nu}\nabla_L(\omega^{1/2}))\pmod{4}.
\]
\end{lemma}
\begin{proof}
Let $S$ be a connected C-complex for $L$.
Fix $\omega\in T^\mu_\ast$ such that $\eta_L(\omega)=0$.
If $z_j$ denotes the complex number $\omega_j^{1/2}$, we have
\begin{eqnarray*}
\overline{H(\omega)}&=&\prod_{j=1}^\mu(1-\omega_j)\,
\sum_\eps\,\eps_1\cdots\eps_\mu\,\omega_1^{\frac{\eps_1-1}{2}}\cdots\omega_\mu^{\frac{\eps_\mu-1}{2}}A^\eps\\
&=&\prod_{j=1}^\mu(\overline z_j^2-1)\,
\sum_\eps\,\eps_1\cdots\eps_\mu\,z_1^{\eps_1+1}\cdots z_\mu^{\eps_\mu+1}A^\eps\\
&=&\prod_{j=1}^\mu(\overline z_j-z_j)\,B(z_1,\dots,z_\mu),
\end{eqnarray*}
where $B(z_1,\dots,z_\mu)$ denotes the matrix $\sum_\eps\,\eps_1\cdots\eps_\mu\,z_1^{\eps_1}\cdots z_\mu^{\eps_\mu}A^\eps$. On
the other hand, we know by \cite{Cim} that the Conway potential function of $L$ is given by
\[
\nabla_L(t_1,\dots,t_\mu)=(-1)^{\frac{c-\ell}{2}}\prod_{j=1}^\mu(t_j-t_j^{-1})^{\chi(S\setminus S_j)-1}\det(-B(t_1,\dots,t_\mu)),
\]
where $c$ is the number of clasps of $S$ and $\ell=\sum_{i<j}lk(L_i,L_j)$.
Let $\beta_1$ denote the first Betti number of $S$.
By Lemma \ref{lemma:herm}, $\sigma_L(\omega)\equiv\beta_1\pmod{4}$ if and only if
\begin{eqnarray*}
0&<&\det(H(\omega))=\det(\prod_{j=1}^\mu(\overline z_j-z_j)\,B(z_1,\dots,z_\mu))\\
&=&(-1)^{\frac{c-\ell}{2}+\beta_1+\mu\beta_1}
\prod_{j=1}^\mu(z_j-\overline z_j)^{1-\chi(S\setminus S_j)+\beta_1}\nabla_L(z_1,\dots,z_\mu).
\end{eqnarray*}
Since $z_j=e^{i\theta}$ with $0<\theta<\pi$, $z_j-\overline{z}_j=\lambda i$ with $\lambda\in\Rr_+^\ast$. Therefore,
$\sigma_L(\omega)\equiv\beta_1\pmod{4}$ if and only if $0<i^\alpha\,\nabla_L(z_1,\dots,z_\mu)$, with
\[
\alpha=c-\ell+2\beta_1-2\mu\beta_1+\sum_{j=1}^\mu(1-\chi(S\setminus S_j)+\beta_1).
\]
An Euler characteristic argument shows that $\alpha=1-\ell+\beta_1$. We know from Lemma \ref{lemma:mod2}
that $\beta_1\equiv\nu+\ell+1\pmod{2}$. Via the transformation $T2$ (recall Figure \ref{fig:Sequiv}),
it may be assumed that $\beta_1\equiv\nu+\ell-1\pmod{4}$. Therefore, $\sigma_L(\omega)\equiv\nu+\ell-1\pmod{4}$
if and only if $i^{\nu}\nabla_L(z_1,\dots,z_\mu)>0$. This fact, together with Lemma \ref{lemma:mod2}, gives the result.
\end{proof}

The following lemma is a direct consequence of Theorem \ref{thm:semi}. Nevertheless, we now give an alternative proof.
\begin{lemma}\label{lemma:null}
Given any colored link $L$ and any $\omega\in T^\mu_\ast$, $\eta_L(\omega)=0$ if and only if $\nabla_L(\omega^{1/2})\neq 0$.
\end{lemma}
\begin{proof}
We saw in the proof of Lemma \ref{lemma:mod4} that
\[
\det(H(\omega))=\pm\prod_j(\omega^{1/2}_j-\overline{\omega^{1/2}_j})^{m_j}\nabla_L(\omega^{1/2})
\]
for some integers $m_j$.
Note that $\eta_L(\omega)\neq 0$ if and only if $\det(H(\omega))=0$. Since $\omega^{1/2}_j\neq\overline{\omega^{1/2}_j}$, the
lemma is checked. 
\end{proof}

\begin{proof}[Proof of Theorem \ref{thm:local}]
$a)$ Let $S'$ be a connected C-complex for $L'$. A C-complex $S$ for $L$ is obtained from $S'$ by attaching a band with a half-twist.
Fix $\omega\in T^\mu_\ast$. If $H'(\omega)$ and $H(\omega)$ denote the associated Hermitian matrices, then clearly
\[
H(\omega)=\begin{pmatrix}H'(\omega)&v\cr\overline v^T&\lambda\end{pmatrix}
\]
for some complex vector $v$ and real number $\lambda$. Therefore,
\[
|\eta_{L'}(\omega)-\eta_L(\omega)|+|\sigma_{L'}(\omega)-\sigma_L(\omega)|=1.\eqno{(\star)}
\]
Let us now assume that $\nabla_{L'}(\omega^{1/2})\neq 0$. By Lemma \ref{lemma:null}, $\eta_{L'}(\omega)=0$. If
$\nabla_L(\omega^{1/2})=0$, then $\eta_L(\omega)>0$. By equation $(\star)$, $\sigma_{L'}(\omega)=\sigma_L(\omega)$, so
the theorem holds in this case. Let us now assume that $\nabla_L(\omega^{1/2})\neq 0$. By Lemma \ref{lemma:null} and
equation $(\star)$, we have $\sigma_L(\omega)=\sigma_{L'}(\omega)+\eps$ for some $\eps=\pm 1$. Reducing this equation modulo $4$,
Lemma \ref{lemma:mod4} implies
\[
\nu+\ell-sgn(i^{\nu}\nabla_L(\omega^{1/2}))\equiv\nu'+\ell'-sgn(i^{\nu'}\nabla_{L'}(\omega^{1/2}))+\eps\pmod{4},
\]
where $\ell=\sum_{i<j}lk(L_i,L_j)$ and $\ell'=\sum_{i<j}lk(L'_i,L'_j)$. Hence,
\[
\eps\equiv(\nu-\nu')+(\ell-\ell')+sgn(i^{\nu'}\nabla_{L'}(\omega^{1/2}))-sgn(i^{\nu}\nabla_L(\omega^{1/2}))\pmod{4}.
\]
Clearly, $\ell=\ell'$ and $\nu=\nu'+\tau$ with $\tau=\pm 1$. Therefore,
\[
\eps\equiv\tau+sgn(i^{\nu'}\nabla_{L'}(\omega^{1/2}))-sgn(i^{\nu}\nabla_L(\omega^{1/2}))\pmod{4}.
\]
Since $\eps=\pm 1$ and $\tau=\pm 1$, this implies
\[
\eps=sgn\left(i^{1-\tau}\,\frac{i^{\nu}\,\nabla_L(\omega^{1/2})}{i^{\nu'}\,\nabla_{L'}(\omega^{1/2})}\right)
=sgn\left(i\,\frac{\nabla_L(\omega^{1/2})}{\nabla_{L'}(\omega^{1/2})}\right),
\]
so the first equality is checked. The second equality can be derived from this one using the following observation.

\begin{figure}[h]
\begin{center}
\epsfig{figure=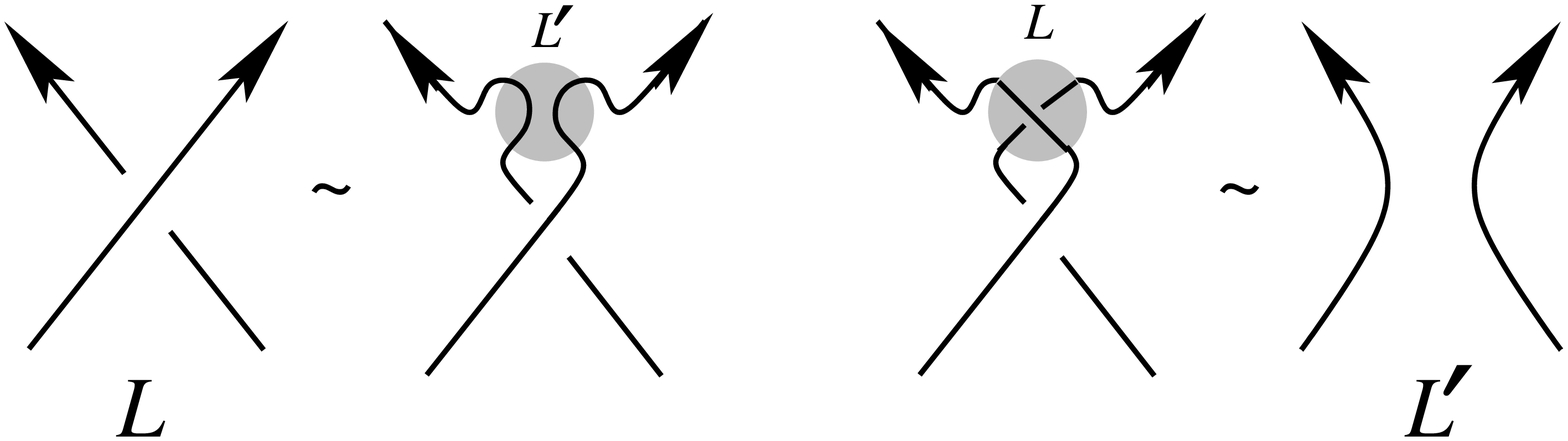,height=2cm}
\end{center}
\end{figure}
\noindent
$b)$ Given $S''$ a connected C-complex for $L''$, a C-complex $S$ for $L$ is obtained from $S''$ by attaching a clasp. As in the proof of
$a)$, this leads to the equation $\sigma_L(\omega)=\sigma_{L''}(\omega)+\eps$ for some $\eps=\pm 1$ satisfying
\[
\eps\equiv(\nu-\nu'')+(\ell-\ell'')+sgn(i^{\nu''}\nabla_{L''}(\omega^{1/2}))-sgn(i^{\nu}\nabla_L(\omega^{1/2}))\pmod{4}.
\]
This time, $\nu=\nu''$ and $\ell=\ell''+\delta$, with $\delta=\pm 1$ as described in the statement of the theorem. Hence,
\[
\eps\equiv\delta+sgn(i^{\nu}\nabla_{L''}(\omega^{1/2}))-sgn(i^{\nu}\nabla_L(\omega^{1/2}))\pmod{4}.
\]
Since $\eps=\pm 1$ and $\delta=\pm 1$, this implies
\[
\eps=\delta\cdot sgn\left(\frac{i^{\nu}\,\nabla_L(\omega^{1/2})}{i^{\nu}\,\nabla_{L''}(\omega^{1/2})}\right)
=\delta\cdot sgn\left(\frac{\nabla_L(\omega^{1/2})}{\nabla_{L''}(\omega^{1/2})}\right).
\]
The last equality is a consequence of this equation.
\end{proof}

\section{The $4$-dimensional viewpoint}\label{section:4-dim}

In this section, we present the signatures $\sigma_L(\omega)$ from the point of view of coverings and
intersection forms. Rougthly speaking, we show that they can be constructed as the Atiyah-Singer invariant of a
finite abelian covering of the link exterior in $S^3$. In doing this, we also relate our signatures to invariants introduced by
Gilmer \cite{Gi1}, Smolinski \cite{Smo}, Levine \cite{Le3} and the second author \cite{Flo}.
 
\begin{definition}
Let $F$ be a union of compact connected oriented smooth surfaces $F_1,\dots,F_\mu$ in $B^4$.
We shall say that  $F$ has \emph{boundary} $L$ if the following holds:
\begin{romanlist}
\item{For all $i$, $F_i$ is smoothly embedded in $B^4$ and $\partial F_i=L_i$.}
\item{For all $i\neq j$, $F_i$ and $F_j$ intersect transversally in a finite number of points, possibly empty.}
\item{For all $i,j,k$ pairwise distinct, $F_i\cap F_j\cap F_k$ is empty.}
\end{romanlist}
\end{definition}

The existence of a surface $F$ with boundary $L$ is obvious. For such an $F$, denote by $W_F$ (or simply by $W$) the 
complement of an open tubular neighborhood of $F$ in $B^4$. By the exact sequence of the pair $(B^4,B^4\setminus F)$ and duality, one 
shows that $H_1(W_F)$ is the free abelian group generated by the meridians $m_1,\dots,m_\mu$ of the components $F_1,\dots,F_\mu$ of $F$.
Therefore, the group of characters Hom$(H_1(W_F),S^1)$ can be identified with the $\mu$ dimensional torus $T^\mu$. To the element
$\omega=(\omega_1,\dots,\omega_\mu)$ of $T^\mu$ corresponds the character $\chi_\omega$ given by $\chi_\omega(m_i)=\omega_i$.
 
Let $T^\mu_\Qq$ be the subset of $T^\mu_*$ constitued by the points with rational coordinates, i.e.
\[
T^\mu_\Qq=\left\{(\omega_1,\dots,\omega_\mu)\in T^\mu\;;\;\omega_j=e^{2i\pi\theta_j}\hbox{ with }\theta_j\in\Qq\,\cap\,]0,1[\,\right\}.
\]
Consider an element $\omega$ of $T^\mu_\Qq$. The image of the corresponding character $\chi_\omega$ is the subgroup of $\Cc^*$
generated by $\alpha=e^{2i\pi/q}$, where $q$ is the least common multiple of the orders of the $\omega_i$'s. Hence, $\chi_\omega$ 
induces a $q$-fold cyclic covering $W^q\to W$ with a canonical deck transformation $\tau$ generating $C_q$, the group of the covering.
The cellular chain complex $C_*(W^q)$ is a $\Zz[C_q]$-module, and so is $\Cc$ via the homomorphism $j\colon C_q\to\Cc^*$ which sends
$\tau$ to $\alpha$. The \emph{twisted homology} of $(W,\omega)$, denoted by $H_*^\omega(W;\Cc)$, is the 
homology of the chain complex $C_*(W^q)\otimes_{\Zz[C_q ]}\Cc$. Note that $\Cc$ is flat over $\Zz[C_q ]$, so
\[
H_*^\omega(W;\Cc) = H_*(W^q)\otimes_{\Zz[C_q ]}\Cc.
\]
The twisted homology $H_*^\omega(X;\Cc)$ of the exterior of $L$ in $S^3$ is defined similarly,
using the character of $H_1(X)$ induced by $\chi_\omega$.

\begin{definition}
For $\omega$ in $T^\mu_\Qq$, let $\varphi_F^\omega$ (or simply $\varphi^\omega$) denote the twisted intersection form defined by 
\[
\begin{array}{ccl}
\varphi_F^\omega \;\colon\;H_2^\omega(W_F;\Cc)\times H_2^\omega(W_F;\Cc)&\longrightarrow& \Cc\\
\phantom{\varphi_\omega(F)\;\colon\;}(x\otimes u, y\otimes v)&\longmapsto&u\overline{v}\sum_{i=1}^{q}<x,\tau^i y>\alpha^i
\end{array}
\]
for $x,y\in H_2(W^q)$ and $u,v\in\Cc$. Here, $<\;,\;>$ is the intersection form induced by the orientation of $W$ lifted to $W^q$,
and $\lambda\mapsto\overline{\lambda}$ is the complex conjugation.
\end{definition}

Note that $\varphi^\omega$ is a well-defined Hermitian form. In fact, it is conjugate to the ordinary intersection form on
$H_2(W^q;\Cc)$ restricted to the eigenspace of $\tau$ with eigenvalue $\alpha$ (see Lemma \ref{eigen} below).

This section is mainly devoted to the proof of the following theorem.

\begin{thm}\label{main1}
Let $L$ be a $\mu$-colored link. For any surface $F$ in $B^4$ with boundary $L$, and for all $\omega\in T^\mu_\Qq$, 
\begin{eqnarray*}
\sigma_L(\omega)&=&\sign(\varphi_F^\omega),\\
\eta_L(\omega)&=&\dim H_1^\omega(X;\Cc)\;=\;\nuli(\varphi_F^\omega)+\dim H_1^\omega(W_F;\Cc)-\dim H_3^\omega(W_F;\Cc).
\end{eqnarray*}
\end{thm}

\subsection{Finite abelian coverings}\label{sub:fac}
   
In this subsection, we relate the twisted intersection forms of finite abelian coverings and of their cyclic quotients.   

As above, let $W$ be the exterior of a surface $F$ with boundary $L$ (or more generally, any compact connected oriented
$4$-dimensional manifold with $H_1(W)=\bigoplus_{i=1}^\mu\Zz$, where the generator of the $i^\mathrm{th}$ summand is denoted $m_i$).   
Consider integers $q_1,\dots,q_\mu$ with $q_i>1$ for all $i$. The natural projection
$\gamma\colon H_1(W)\to G=C_{q_1}\times\dots\times C_{q_\mu}$ mapping $m_i$ to a preferred generator $\tau_i$ of
$C_{q_i}$ induces a finite abelian covering $W^\gamma\to W$. As usual,
$H_*(W^\gamma)$ is a module over the group ring $\Zz[G]$. Let $S\colon H_2(W^\gamma)\times H_2(W^\gamma)\to\Zz[G]$
be the pairing given by
\[
S(x,y)=\sum_{g\in G}<x,gy>g,
\]
where $<\;,\;>$ is the intersection form induced by the orientation of $W$ lifted to $W^\gamma$. Note that $S$ is sesquilinear
with respect to the involution of $\Zz[G]$ induced by $g\mapsto g^{-1}$.

Now, consider an element $\omega$ of $T^\mu_\Qq$ with $\omega_i$ of order $q_i$ for all $i$. As above, let $q$ denote the least common
multiple of the $q_i$'s. The corresponding character $\chi_\omega\colon H_1(W)\to\Cc^*$ can be written
$\chi_\omega=j\circ r_\omega\circ\gamma$, where $j\colon C_q\to\Cc^*$ satisfies
$j(\tau)=\alpha$ and $r_\omega\colon G\to C_q$ is such that $j\circ r_\omega$ maps $\tau_i=\gamma(m_i)$ to $\omega_i$.
Let $s_\omega\colon G\to\Cc^*$ denote the character given by $j\circ r_\omega$. It induces a homomorphism
of rings $\Zz[G]\to\Cc$ compatible with the involutions of these rings. In particular, it endows $\Cc$ with a structure of
$\Zz[G]$-module. Let $\psi^\omega$ (or $\psi_F^\omega$) denote the Hermitian form on
$H_2(W^\gamma)\otimes_{\Zz[G]}\Cc$ induced by the $\Zz[G]$-sesquilinear form $S$ via the ring homomorphism $\Zz[G]\to\Cc$.
In other words, $\psi^\omega$ is given by
\[
\begin{array}{ccl}
\psi^\omega \;\colon\;H_2(W^\gamma)\otimes_{\Zz[G]}\Cc\times H_2(W^\gamma)\otimes_{\Zz[G]}\Cc&\longrightarrow& \Cc\\
\phantom{\psi^\omega \;\colon\;}(x\otimes u, y\otimes v)&\longmapsto&u\overline{v}\sum_{g\in G}<x,gy>s_\omega(g).
\end{array}
\]

\begin{lemma}\label{eigen}
Given a character $s\colon G\to\Cc^*$, set
\[
E_s=\{x\in H_2(W^\gamma;\Cc)\,;\,\hbox{$gx=s(g)x$ for all $g\in G$}\}.
\]
Then, $\psi^\omega$ is conjugate to the restriction to $E_{s_\omega}$ of the Hermitian intersection form on $H_2(W^\gamma;\Cc)$.
\end{lemma}
\begin{proof} 
First, observe that $H_2(W^\gamma)\otimes_{\Zz[G]}\Cc$ is equal to $H_2(W^\gamma;\Cc)\otimes_{\Cc[G]}\Cc$. For all
$s\in\mathcal{R}_G=\hbox{Hom}(G,\Cc^*)$, let $\vartheta_{s}$ be the element of $\Cc[G]$ given by
\[
\vartheta_s=\frac{1}{|G|}\sum_{g\in G}\overline{s(g)}g.
\]
One easily checks that $g\vartheta_s=s(g)\vartheta_s$ for all $g\in G$ and $s\in\mathcal R_G$. This gives the inclusion
$\vartheta_s H_2(W^\gamma;\Cc)\subset E_s$. On the other hand, $\vartheta_s\vartheta_s=\vartheta_s$
and $\vartheta_s\vartheta_{s'}=0$ if $s\neq s'$. Hence, if $x$ belongs to $E_s$, then $\vartheta_{s'}x=0$
for all $s\neq s'$. Since $\sum_{s\in\mathcal{R}_G}\vartheta_s=1$ (see \cite[Lemma 2.1]{Sak}), we have $x=\vartheta_sx$ and the
inclusion $E_s\subset\vartheta_s H_2(W^\gamma;\Cc)$ is proved. In other words, there is a canonical  
isomorphism of $\Cc[G]$-modules
\[
H_2(W^\gamma;\Cc)\simeq\bigoplus_{s\in\mathcal R_G}\vartheta_s H_2(W^\gamma;\Cc)=\bigoplus_{s\in\mathcal R_G}E_s.
\]
Since $G$ acts by isometries, this decomposition is orthogonal with respect to the intersection form $<\;,\;>$ on $H_2(W^\gamma;\Cc)$.
Now, observe that $E_s\otimes_{\Cc[G]}\Cc=0$ for all $s\neq s_\omega$, and $E_{s_\omega}\otimes_{\Cc[G]}\Cc=E_{s_\omega}$.
Therefore, the multiplication
$H_2(W^\gamma;\Cc)\to\vartheta_\omega H_2(W^\gamma;\Cc)$ by $\vartheta_\omega:=\vartheta_{s_\omega}$ induces an isomorphism
\[
H_2(W^\gamma;\Cc)\otimes_{\Cc[G]}\Cc\simeq E_{s_\omega}.
\]
Since $\psi^\omega(x\otimes u,y\otimes v)=|G|u\overline{v}<x,\vartheta_\omega y>=
|G|u\overline{v}<\vartheta_\omega x,\vartheta_\omega y>=|G|<\vartheta_\omega ux,\vartheta_\omega vy>$ for all
$x,y\in H_2(W^\gamma;\Cc)$ and $u,v\in\Cc$, the lemma is proved.
\end{proof}

If the $q_i$'s are pairwise coprime, then $r_\omega\colon G\to C_q$ is an isomorphism, so $W^\gamma=W^q$. In this case,
$H_2(W^\gamma)\otimes_{\Zz[G]}\Cc=H_2^\omega(W;\Cc)$ and the forms $\psi^\omega$ and $\varphi^\omega$ are equal. In general, we have
the following result, very much in the spirit of \cite{Sak}. 

\begin{lemma}\label{decomp}
Let $\omega$ be an element of $T^\mu_\Qq$ and consider $\gamma$ and $s_\omega$ such that $\chi_\omega=s_\omega\circ\gamma$.
Then, the Hermitian forms $\varphi^\omega$ and $\psi^\omega$ are conjugate.
\end{lemma}
\begin{proof}
By the proof of Lemma \ref{eigen}, the multiplication by $\vartheta_\omega=\vartheta_{s_\omega}$ induces an isomorphism
\[
H_2(W^\gamma;\Cc)\otimes_{\Cc[G]}\Cc\to\vartheta_\omega H_2(W^\gamma;\Cc)\otimes_{\Cc[G]}\Cc,
\]
so that $\psi^\omega$ is conjugate to  $|G|<\;,\;>$ restricted to $\vartheta_\omega H_2(W^\gamma;\Cc)\otimes_{\Cc[G]}\Cc$
(which is nothing but $E_{s_\omega}$). For the same reason, the multiplication by
$r_\omega(\vartheta_\omega)=\frac{1}{q}\sum_{i=1}^q\overline{\alpha}^i\tau^i$ induces an isomorphism
\[
H_2(W^q;\Cc)\otimes_{\Cc[C_q]}\Cc\to r_\omega(\vartheta_\omega)H_2(W^q;\Cc)\otimes_{\Cc[C_q]}\Cc,
\]
so that $\varphi^\omega$ is conjugate to  $q<\;,\;>$ restricted to $r_\omega(\vartheta_\omega)H_2(W^q;\Cc)\otimes_{\Cc[C_q]}\Cc$.

Set $B=\ker(s_\omega)$ and $\mathcal R_G(B)=\{s\in\mathcal R_G | B\subset \ker(s)\}$. By \cite[p. 205]{Sak}, we have
natural isomorphisms of $\Cc[G]$-modules
\[
H_2(W^q;\Cc)\simeq(\sum_{b\in B}b)H_2(W^\gamma;\Cc)\simeq|B|\bigoplus_{s\in\mathcal R_G(B)}\vartheta_s H_2(W^\gamma;\Cc),
\]
where $H_2(W^q;\Cc)$ is a $\Cc[G]$-module via the homomorphism $\Cc[G]\to\Cc[C_q]$ induced by $r_\omega\colon G\to C_q$.
Since $s_\omega$ belongs to $\mathcal R_G(B)$, this gives an isomorphism
\[
r_\omega(\vartheta_\omega)H_2(W^q;\Cc)\simeq|B|\vartheta_\omega H_2(W^\gamma;\Cc).
\]
Tensoring this isomorphism by $\Cc$, we see that the intersection form on the space
$r_\omega(\vartheta_\omega)H_2(W^q;\Cc)\otimes_{\Cc[C_q]}\Cc$ is conjugate to the intersection form on the space\break
$\vartheta_\omega H_2(W^\gamma;\Cc)\otimes_{\Cc[C_q]}\Cc$ multiplied by $|B|^2$.
Summing up everything, $\varphi^\omega$ is conjugate to $\psi^\omega$ multiplied by $\frac{q|B|^2}{|G|}=|B|$. 
\end{proof}

\subsection{Intersection forms and Seifert forms}\label{sub:int}

For technical reasons, we shall work with branched coverings instead of regular coverings. As above, let $F\subset B^4$ be a surface with
boundary $L$, and let $W_F$ denote its exterior. Given $\omega\in T^\mu_\Qq$, let $q_i$ be the order of $\omega_i$.
The regular covering $W_F^\gamma\to W_F$ induced by the projection $\gamma\colon H_1(W_F)\to G=C_{q_1}\times\dots\times C_{q_\mu}$
can be extended to a covering $p\colon\overline{W}_F^\gamma\to B^4$ branched along $F$, such that $p|\colon p^{-1}(F_i)\to F_i$
is a $G/C_{q_i}$-covering branched along $F_i\cap F_j$ with branch index $q_j$. The construction of $\overline{W}_F^\gamma$ is quite
clear except near the intersection points of two components of $F$, where it can be described as follows. An intersection point of
$F_i$ and $F_j$ has a neighborhood $N$ in $B^4$ such that the triple $(N,F_i,F_j)$ is diffeomorphic to
$(B^2\times B^2,B^2\times\{0\},\{0\}\times B^2)$. Then, $p^{-1}(N)$ is a disjoint union of balls, and the covering projection
restricted to each of these balls is given by
$B^2\times B^2\to B^2\times B^2\simeq N$, $(z_1,z_2)\mapsto (e^{2i\pi/q_i}z_1,e^{2i\pi/q_j}z_2)$.
Note that the boundary of $\overline{W}_F^\gamma$ is the covering of $S^3$ branched along $L$ induced by the restriction of $\gamma$ to
the homology of the exterior of $L$.

Clearly, we can define a twisted intersection form $\overline{\psi}_F^\omega$ on
$H_2(\overline{W}_F^\gamma)\otimes_{\Zz[G]}\Cc$ as explained in the previous subsection. Furthermore, we have the following result.

\begin{lemma}\label{ram}
The forms  $\psi^\omega_F$ and $\overline{\psi}_F^\omega$ are conjugate.
\end{lemma}
\begin{proof}
By Lemma \ref{eigen}, we just need to check that the inclusion $W_F^\gamma\subset \overline W_F^\gamma$ induces an isomorphism of
$\Cc$-vector spaces $\vartheta_\omega H_2(W_F^\gamma;\Cc)\simeq\vartheta_\omega H_2(\overline W_F^\gamma;\Cc)$. This follows from
a standard Mayer-Vietoris argument (see e.g. \cite[Lemma 5]{Smo} and \cite[Lemma 6.3]{Coop}).
\end{proof}
  
This subsection is devoted to the proof of the following proposition, which generalizes \cite{Vi} and \cite[Proposition 6.1]{Coop}.
Note that this result implies that our signatures coincide with invariants introduced by Smolinski \cite{Smo}.

\begin{prop} \label{ccompl}
Consider a connected C-complex $S\subset S^3$ for a $\mu$-colored link $L$. For 
$\omega\in T_\Qq^\mu$, let $H(\omega)$ be the corresponding Hermitian matrix (recall Section \ref{section:def}). If 
$F$ is the surface with boundary $L$ obtained by pushing $S$ in $B^4$, then $H(\omega)$ is a matrix for 
$\overline\psi^\omega_F$.
\end{prop}
\begin{proof}
As the surfaces $S_i$ are pushed into $B^4$, they trace out $3$-manifolds $M_i$ homeomorphic to $S_i\times[0,1]$ with $L_i\times[0,1]$
collapsed to a single copy of $L_i$. It may be assumed that the $M_i$'s intersect transversally, so each clasp in $S_i\cap S_j$ gives
a $2$-disc in $M_i\cap M_j$. Now, split $B^4$ along $M=M_1\cup\dots\cup M_\mu$. The boundary of the resulting manifold $B$ contains
two copies $M^+$ and $M^-$ of $M$ which intersect along $F$. Let $\{gB\}_{g\in G}$ be disjoint copies of $B$.
Consider the $4$-dimensional $G$-manifold $E$ obtained by pasting these copies along $gM_i^\pm\subset gB$ as follows:
\[
E=\bigsqcup_{g\in G}gB\Big/ (\tau_i gM_i^+=gM_i^-)_{1\le i\le\mu,\,g\in G}.
\]
(Recall that $\tau_i$ is a preferred generator of $C_{q_i}$.) Clearly, the projection $E\to E/G$ is nothing but the branched
covering $\overline W_F^\gamma\to B^4$. Furthermore, there is a deformation retract from $B$ to the cone $C(M^+\cup M^-)$ over
$M^+\cup M^-$, which itself retracts by deformation on the cone $CF$ over $F$. Therefore, $\overline W_F^\gamma=E$ is homotopy
equivalent to the $G$-space
\[
Y=\bigsqcup_{g\in G}gCF\Big/(\tau_i gF_i=gF_i)_{1\le i\le\mu,\,g\in G}.
\]
We shall now use this description of $\overline W_F^\gamma$ to compute the $\Zz[G]$-module $H_2(\overline W_F^\gamma)=H_2(Y)$.
Consider a basis of $H_1(F)$ whose elements are given by $1$-cycles $\{e_\alpha\}_{\alpha\in A}$. For $\alpha\in A$,
let $I(\alpha)$ denote the set of indices $i$ such that $e_\alpha$ meets $F_i$. Given such an $e_\alpha$ in $F$,
let $Ce_\alpha\subset CF$ denote the cone over $e_\alpha$, and $\Sigma e_\alpha$ the `suspension' defined by
\[
\Sigma e_\alpha=\prod_{i\in I(\alpha)}(1-\tau_i)Ce_\alpha\subset Y.
\] 
By construction of $Y$, $\Sigma e_\alpha$ is a $2$-cycle in $Y$. Furthermore, a massive use of the Mayer-Vietoris exact sequence leads
to the following fact: the $\Zz[G]$-module $H_2(Y)$ is given by
\[
H_2(Y)=\bigoplus_{\alpha\in A}\left(\Zz[G]\Big/(1+\tau_i+\dots+\tau_i^{q_i-1})_{i\in I(\alpha)}\right)\Sigma e_\alpha.
\]
Using the fact that $\omega_i\neq 1$ for $i\in I(\alpha)$, we have
\begin{eqnarray*}
H_2(Y)\otimes_{\Zz[G]}\Cc&=&
\bigoplus_{\alpha\in A}\left(\left(\Zz[G]\Big/(1+\tau_i+\dots+\tau_i^{q_i-1})_{i\in I(\alpha)}\right)
\otimes_{\Zz[G]}\Cc\right)\Sigma e_\alpha\\
&=&\bigoplus_{\alpha\in A}\Cc\Sigma e_\alpha.
\end{eqnarray*}
Furthermore, since $\omega_i\neq 1$ for $i\notin I(\alpha)$, $x_\alpha=\prod_{i=1}^\mu(1-\tau_i)Ce_\alpha$ is a non-zero multiple
of $\Sigma e_\alpha$ in $H_2(Y)\otimes_{\Zz[G]}\Cc$. Therefore, a basis of the latter space is given by $\{x_\alpha\}_{\alpha\in A}$.

We are left with the computation of the form $\overline\psi^\omega_F$ on the elements of this basis.
To do this, we shall deform $x_\alpha$ into another $2$-cycle $\widetilde x_\alpha$ as follows. First, note that
\[
x_\alpha=\sum_\eps \eps_1\cdots\eps_\mu\tau_1^{\frac{1-\eps_1}{2}}\cdots\tau_\mu^{\frac{1-\eps_\mu}{2}}Ce_\alpha,
\]
where the sum is on all sequences $\eps=(\eps_1,\dots,\eps_\mu)$ of $\pm 1$'s. Set
\[
\widetilde x_\alpha=
\sum_\eps \eps_1\cdots\eps_\mu\tau_1^{\frac{1-\eps_1}{2}}\cdots\tau_\mu^{\frac{1-\eps_\mu}{2}}C(i^\eps(e'_\alpha)),
\]
where $e_\alpha'$ is the $1$-cycle in $S$ corresponding to $e_\alpha$ in $F$, and $i^\eps$ is the map defined in Section \ref{section:def}. Clearly,
$\widetilde x_\alpha$ is a $2$-cycle and $<\widetilde x_\alpha,gx_\beta>=<x_\alpha,gx_\beta>$ for all $g\in G$ and $\beta\in A$.
Furthermore, the intersection number between $C(i^\eps(e'_\alpha))$ and $gCe_\beta$ is given by
\[
C(i^\eps(e'_\alpha))\cdot (gCe_\beta)=\begin{cases}lk(i^\eps(e'_\alpha),e'_\beta)&\text{if $g=1$;}\cr 0&\text{else.}\end{cases}
\]
Since $G$ acts by isometries, this is sufficient in order to determine
$\overline\psi^\omega_F(x_\alpha,x_\beta)=\sum_{g\in G}<\widetilde x_\alpha,gx_\beta>s_\omega(g)$.
The result of this tedious computation is
\[
\overline\psi^\omega_F(x_\alpha,x_\beta)=
\prod_{i=1}^\mu(1-\overline{\omega}_i)
\sum_\eps\eps_1\cdots\eps_\mu\,\omega_1^{\frac{1-\eps_1}{2}}\cdots\omega_\mu^{\frac{1-\eps_\mu}{2}}\,lk(i^\eps(e'_\alpha),e'_\beta),
\]
concluding the proof.
\end{proof}

\subsection{Proof of Theorem \ref{main1}}\label{sub:proof4dim}

We need one more result.

\begin{prop}
For all $\omega$ in $T^\mu_\Qq$, $\sign(\varphi_F^\omega)$ does not depend on the choice of the surface $F$ in $B^4$ with boundary $L$.
\label{indep}
\end{prop}
\begin{proof}
By Lemmas \ref{eigen}, \ref{decomp} and \ref{ram}, $\varphi_F^\omega$ is conjugate to the restriction to
$E_{s_\omega}$ of the intersection form on $H_2(\overline W_F^\gamma;\Cc)$. Moreover, the signature of the latter form
is easily seen to be given by a linear combination of the $g$-signatures of $\overline W_F^\gamma$, denoted by
$\sigma(g,\overline W_F^\gamma)$ (see \cite{Gor2}).
So, we are left with the proof that $\sigma(g,\overline W_F^\gamma)$ does not depend on the choice of $F$. Consider two surfaces $F$
and $F^\prime$ with boundary $L$, and glue together $\overline{W}_F^\gamma$ and $-\overline{W}_{F^\prime}^\gamma$ along their common
boundary. By Novikov additivity, the resulting manifold $Z$ satisfies
$\sigma(g,Z)=\sigma(g,\overline{W}_F^\gamma)-\sigma(g,\overline{W}_{F^\prime}^\gamma)$. Using the Atiyah-Singer $G$-signature theorem
\cite{Aty-Sin}, one easily checks that $\sigma(g,Z)$ is zero for all $g\in G$ (see \cite[p. 216]{Smo}). This completes the proof.
\end{proof}

\begin{proof}[Proof of Theorem \ref{main1}]
The equality $\sigma_L(\omega)=\sign(\varphi_F^\omega)$ is a direct consequence of Lemma \ref{decomp}, Lemma \ref{ram},
Proposition \ref{ccompl} and Proposition \ref{indep}. On the other hand, consider the exact sequence of the pair $(W_F,\partial W_F)$
with (twisted) complex coefficients. Clearly, a matrix for $\varphi^\omega_F$ is (the transpose of) a matrix for
$H_2^\omega(W_F;\Cc)\to H_2^\omega(W_F,\partial W_F;\Cc)$. This leads to
\[
\nuli(\varphi^\omega_F)=\dim H_1^\omega(\partial W_F;\Cc)-\dim H_1^\omega(W_F;\Cc)+\dim H_1^\omega(W_F,\partial W_F;\Cc).
\]
By duality, $\dim H_1^\omega(W_F,\partial W_F;\Cc)=\dim H_3^\omega(W_F;\Cc)$. Finally, a standard Mayer-Vietoris argument shows that
the inclusion $X\subset\partial W_F$ induces an isomorphism  $H_1^\omega(X;\Cc)=H_1^\omega(\partial W_F;\Cc)$. This gives
\[
\dim H_1^\omega(X;\Cc)=\nuli(\varphi^\omega_F)+\dim H_1^\omega(W_F;\Cc)-\dim H_3^\omega(W_F;\Cc).
\]
Note that this equation holds for any surface $F$ in $B^4$ with boundary $L$.

We now turn to the proof of the equality $\eta_L(\omega)=\dim H_1^\omega(X;\Cc)$, assuming the notation of the previous subsections.
Let $F$ be a connected C-complex for $L$ pushed in $B^4$. By Lemma \ref{decomp}, Lemma \ref{ram} and Proposition \ref{ccompl},
$\eta_L(\omega)=\nuli(\varphi^\omega_F)$, so we are left with the proof that $H_1^\omega(W_F;\Cc)=H_3^\omega(W_F;\Cc)=0$ for such an $F$.
By the proof of Lemma \ref{decomp}, $H_*^\omega(W_F;\Cc)=H_*(W_F^\gamma)\otimes_{\Zz[G]}\Cc$. Furthermore, a Mayer-Vietoris argument
shows that $H_*(W^\gamma_F)=H_*(\overline W^\gamma_F)$. Finally, recall that $\overline W^\gamma_F$ has the homotopy type of a $G$-space
obtained by gluing cones over $F$ along their boundary (see the proof of Proposition \ref{ccompl}). Since $F$ is connected and has no
closed component, $H_1(\overline W^\gamma_F)=H_3(\overline W^\gamma_F)=0$. This concludes the proof.
\end{proof}

\smallskip
This $4$-dimensional interpretation of $\sigma_L(\omega)$ can be used to give an alternative proof (for $\omega\in T^\mu_\Qq$) of
several results obtained in Section \ref{section:def}. For example, Propositions \ref{ccompl} and \ref{indep} imply the independance of
$\sigma_L(\omega)$ on the C-complex (Theorem \ref{thm:inv}).

As promised in Section \ref{section:def}, we shall now give another proof of Proposition \ref{prop:lk}. Consider
$\omega=(\omega_1,\dots,\omega_\mu,\omega_\mu)\in T^{\mu+1}_\Qq$, and let $F=F_1\cup\dots\cup F_\mu\cup F_{\mu+1}$ be a surface in $B^4$
with boundary $L$. In a closed neighborhood of an intersection point of $F_\mu$ and $F_{\mu+1}$, $F$ is given by 
two transversal disks with boundary a $2$-colored Hopf link. Let $F'=F'_1\cup\dots\cup F'_\mu$ be the surface obtained by
replacing these disks by a cylinder, whose boundary is a $1$-colored Hopf link. Clearly, $F'$ has boundary the $\mu$-colored link $L'$.
Note that the signature of a $2$-colored Hopf link is zero, while the signature of the positive (resp. negative) $1$-colored Hopf link
is $-1$ (resp. $+1$). By the additivity of the signature, these local transformations decrease $\sigma_L(\omega)$ by
the algebraic intersection number $F_\mu\cdot F_{\mu+1}=lk(L_\mu,L_{\mu+1})$.

\subsection{On the Casson-Gordon invariants of 3-manifolds}\label{sub:CG}

Let $M$ be an oriented closed $3$-manifold, and let $\chi\colon H_1(M)\to\Cc^*$ be a character of finite order. The following
reformulation of the Atiyah-Singer invariant \cite{Aty-Sin} of $(M,\chi)$ is due to Casson and Gordon \cite{Cass-Gord1,Cass-Gord2}.
Since $\chi$ is of finite order, its image is the cyclic subgroup of $\Cc^*$ generated by $\alpha= e^{2i\pi/q}$, for some $q$.
It induces a $q$-fold cyclic covering $M^q\to M$ with a canonical deck transformation $\tau$ generating the group $C_q$ of the covering.
Since the bordism group $\Omega_3(BC_q)$ is equal to $C_q$, there is a positive integer $n$ such that $n$ disjoint copies of
$M$ bound a compact oriented $4$-manifold $W$ over $BC_q$. Let $W^q\to W$ be the induced $q$-fold covering. The deck transformation
$\tau$ extends to a deck transformation of $W^q$, also denoted by $\tau$. As above, let $H^\chi_*(W;\Cc)$ denote the homology of the
chain complex $C_*(W^q)\otimes_{\Zz[C_q]}\Cc$, where the structure of $\Zz[C_q]$-module on $\Cc$ is given by the map $C_q\to\Cc^*$
which sends $\tau$ to $\alpha=e^{2i\pi/q}$. Finally, let $\varphi^\chi$ be the twisted intersection form on $H^\chi_2(W;\Cc)$.
\begin{definition} 
The \emph{Casson-Gordon invariant\/} of $(M,\chi)$ is $\sigma(M,\chi)=\frac{1}{n}(\sign(\varphi^\chi)-\sign(W))$. The related
nullity is defined by $\eta(M,\chi)=\dim H_1^\chi(M;\Cc)$.
\end{definition}
The fact that $\sigma(M,\chi)$ depends only on the pair $(M,\chi)$ is a consequence of the Atiyah-Singer $G$-signature theorem,
and Novikov additivity. Note that one may also consider branched coverings in order to define or compute this invariant (see e.g.
\cite[Proposition 3.5]{Gi1}).

\smallskip
The aim of this subsection is to relate the Casson-Gordon invariant of a manifold obtained by surgery on a framed link $L$
to the signature $\sigma_L$ of this link. We have the following result.

\begin{thm}\label{as}
Let $M$ be the $3$-manifold obtained by surgery on a framed link $L$ with $\nu$ components and linking matrix $\Lambda$.
Let $\chi\colon H_1(M)\to\Cc^*$ be the character mapping the meridian of the $i^{th}$ component of $L$ to $\alpha^{n_i}$,
where $\alpha=e^{2i\pi/q}$ and $n_i$ is an integer coprime to $q$. Consider $L$ as a $\nu$-colored link and set
$\omega=(\alpha^{n_1},\dots,\alpha^{n_\nu})$. Then,
\begin{eqnarray*}
\sigma(M,\chi)&=&\Big(\sigma_L(\omega)-\sum_{i<j}\Lambda_{ij}\Big)-\sign(\Lambda)+\frac{2}{q^2}\sum_{i,j}(q-n_i)n_j\Lambda_{ij},\\
\eta(M,\chi)&=&\eta_L(\omega).
\end{eqnarray*}
\end{thm}

Note that if all the $n_i$'s are equal, then this formula together with Proposition \ref{prop:lk} give back Casson and Gordon's
\cite[Lemma 3.1]{Cass-Gord2}. See also Gilmer \cite[Theorem 3.6]{Gi1}. On the other hand, if the matrix $\Lambda$ is zero, then
$\sigma_L(\omega)=\sigma(M,\chi)$. Hence our signature extends (a special case of) the invariant introduced by Levine in \cite{Le3}
for links with $\Lambda_{ij}=0$ for all $i\neq j$. This also relates $\sigma_L(\omega)$ with link signatures of Gilmer \cite{Gi1}.

Let us point out an interesting feature of this result before giving its proof. Recall that all the signatures $\sigma_L(\omega)$
of a fixed colored link $L$ are given by the signature of a single matrix $H$ evaluated at $\omega$. Using Theorem \ref{as}, one can
compute the Casson-Gordon invariants of all $3$-manifolds obtained by surgery on $L$, for many characters of finite order,
by using this single matrix.
   
\begin{proof}
Let $F=F_1\cup\dots\cup F_\nu$ be the surface obtained by pushing a connected $C$-complex for $L$ in $B^4$. Denote by $\{p_\ell\}_\ell$
the finite set of double points of $F$ (coming form the clasps of the C-complex) and by $\{B_\ell\}_\ell$ a set of small disjoint closed
$4$-balls such that $p_\ell\in int\,B_\ell$ for all $\ell$. We shall denote by $F_{j(\ell)}$ and $F_{k(\ell)}$ the components of
$F$ which intersect at $p_\ell$. Note that $F\cap B_\ell$ consists of two transverse discs with boundary a Hopf link
$K_\ell\subset S_\ell^3=\partial B_\ell$. Let $\epsilon_\ell=\pm 1$ denote the linking number of the components of $K_\ell$ (that is,
the algebraic intersection number of $F_{j(\ell)}$ and $F_{k(\ell)}$ at $p_\ell$).
Set $\Omega=B^4\setminus\bigsqcup_\ell int\,B_\ell$. Clearly, the character
$\chi$ restricted to the exterior of $L$ in $S^3$ extends to a character on the exterior of $F$ in $B^4$, which itself restricts to
the exterior of $K_\ell$ in $S^3_\ell$. This restriction maps the meridians of $K_\ell$ to $\chi(m_{j(\ell)})=\alpha^{n_{j(\ell)}}$ and
$\chi(m_{k(\ell)})=\alpha^{n_{k(\ell)}}$.

Now, let $U$ be the $4$-manifold obtained by attaching $2$-handles $B^2\times B^2$ to $\Omega$
as follows. First, attach $2$-handles to $\Omega$ along a tubular neighborhood of $L\subset S^3$ according to the framings
$\{\Lambda_{ii}\}_i$. Then, for all $\ell$, perform a surgery along $K_\ell\subset S^3_\ell$ according to framings $f^\ell_{j(\ell)}$ and
$f^\ell_{k(\ell)}$ which satisfy the following property: the character on the exterior of $K_\ell$ in $S^3_\ell$ extends to the
$3$-manifold $M_\ell$ obtained by surgery on the framed link $K_\ell$. This is the case if and only if the congruences
$f_{j(\ell)}^\ell n_{j(\ell)}+\epsilon_\ell n_{k(\ell)}\equiv\epsilon_\ell n_{j(\ell)}+f^\ell_{k(\ell)}n_{k(\ell)}\equiv
0\;(\hbox{mod } q)$ hold. By hypothesis, $n_{j(\ell)}$ and $n_{k(\ell)}$ are invertible modulo $q$, so such framings exist.

Let $F'\subset int\,U$ be the smooth closed surface with $\nu$ connected components obtained from $F\cap\Omega$ by gluing the cores
of the $2$-handles. Let $U_{F'}$ be the exterior of $F'$ in $U$. From $\chi$ (that is, from $\omega$), one easily constructs a character
on $H_1(U_{F'})$ inducing a twisted intersection form $\varphi^\omega$ on $H_2^\omega(U_{F'};\Cc)$. Since
$\partial U=M\sqcup\bigsqcup_\ell(-M_\ell)$, Gilmer's \cite[Proposition 3.5]{Gi1} gives 
\[
\sigma(M,\chi)-\sum_\ell\sigma(M_\ell,\chi)=\sign(\varphi^\omega)-\sign(U)+\frac{2}{q^2}\sum_{i=1}^\nu(q-n_i)n_i\big(F'_i\cdot F'_i\big).
\]
(Note that a Mayer-Vietoris argument shows that $\varphi^\omega$ is conjugate to the form related to the covering of $U$ branched along
$F'$, considered by Gilmer in his formula.) We shall now compute separately each term of this equation.

Recall that $\omega$ induces a $C_q\times\dots\times C_q$-covering of $B^4$ branched along $F$. By Proposition \ref{ccompl},
$\sigma_L(\omega)$ is equal to the signature of the corresponding intersection form $\overline{\psi}_F^\omega$. Furthermore, the
signature
corresponding to the covering of $B_\ell$ branched along $F\cap B_\ell$ is equal to zero. Therefore, Novikov additivity implies that
$\sign(\overline{\psi}_F^\omega)$ is equal to the signature of the twisted intersection form of the covering of $\Omega$ branched along
$F\cap\Omega$. Finally, a standard Mayer-Vietoris argument shows that adding $2$-handles to a $4$-manifold has no effect on its twisted
signature, so
\[
\sign(\varphi^\omega)=\sigma_L(\omega).
\]

One easily checks that a matrix for the intersection form on $H_2(U)$ is given by $\Lambda\oplus\bigoplus_\ell(-\Lambda_\ell)$, where
$\Lambda_\ell$ is the linking matrix of $K_\ell$. Therefore,
\[
\sign(U)=\sign(\Lambda)-\sum_\ell\sign(\Lambda_\ell).
\]

Using \cite[Proposition 3.8 and p. 367]{Gi1}, we obtain $\sigma(M_\ell,\chi)=-\epsilon_\ell-\sign(\Lambda_\ell)+\frac{2}{q^2}s_\ell$,
where the integer $s_\ell$ is given by
\[
s_\ell=(q-n_{j(\ell)})n_{j(\ell)}f^\ell_{j(\ell)}+(q-n_{k(\ell)})n_{k(\ell)}f^\ell_{k(\ell)}+
\epsilon_\ell\left((q-n_{j(\ell)})n_{k(\ell)}+(q-n_{k(\ell)})n_{j(\ell)}\right).
\]

Finally, $F_i'\cdot F_i'=\Lambda_{ii}-\sum_\ell f^\ell_i$, where the sum is on all indices $\ell$ such that $j(\ell)=i$ or $k(\ell)=i$.
The first equality of the theorem now follows from the fact that $\sum_\ell\epsilon_\ell=\sum_{i<j}\Lambda_{ij}$.

Since $M$ is obtained from the link complement $X$ by adjoining tori, an easy Mayer-Vietoris argument gives the equality between the
nullities.
\end{proof}

\section{Concordance and existence of surfaces in $B^4$}\label{section:slice}

The properties of $\sigma_L(\omega)$ and $\eta_L(\omega)$ studied in this section do not hold for all
$\omega$ in $T_*^\mu$. We shall denote by $T^\mu_\pp$ the dense subset of $T_*^\mu$ constitued by the 
elements $\omega=(\omega_1,\dots,\omega_\mu)$ which satisfy the following condition: there exists a prime $p$ such that
for all $i$, the order of $\omega_i$ is a power of $p$.

We first prove the invariance of the restriction of $\sigma_L$ and $\eta_L$ to $T^\mu_\pp$ under (colored)
concordance. Then, we show that the signature and nullity provide a lower bound for the genus of a surface
in $B^4$ with boundary $L$ (Theorem \ref{thm:MT}). Finally subsection \ref{sub:KT} deals with an analogous
result concerning surfaces in $S^4$ whose intersection with a standardly embedded $3$-sphere in equal to the colored
link $L$. These results generalize celebrated theorems of Murasugi-Tristram and Kauffman-Taylor.
 
\begin{definition}
Two colored links $L$ and $L^\prime$ with $\nu$ components
are said to be \emph{concordant} if there exists 
a collection of smooth disjoint cylinders $T_1,\dots,T_\nu$ properly embedded in $S^3 \times [0,1]$, such that for all $i$,
$T_i$ is a concordance between components of $L$ and $L'$ of the same color.
\label{conc}
\end{definition} 

\begin{thm}
For all $\omega \in T^\mu_\pp$, $\sigma_L(\omega)$ and $\eta_L(\omega)$ are concordance invariants.
\end{thm}

This result follows from the fact that the exterior of the concordance is a homology
cobordism. The detailed proof can be found in \cite[Theorem 4.15]{Flo} for the case of colored
links with $lk(L_i,L_j)=0$ for all $i\neq j$. It obviously extends to the general case.
Note that this theorem can also be viewed as a consequence of \cite[Theorem 9]{Gil-Liv}.
  
\subsection{Surfaces with double points and Murasugi-Tristram inequality}\label{sub:MT}
 
\begin{thm}\label{thm:MT}
Suppose that $F=F_1\cup\dots\cup F_\mu$ in $B^4$ has boundary $L$ (in the sense of Section \ref{section:4-dim}).
Set $\beta_1 = \sum_i\rk H_1(F_i)$, and let $c$ be the number of double points of $F$. Then, for all $\omega \in T^\mu_\pp$, 
\[
|\sigma_L(\omega)| + |\eta_L(\omega) - \mu + 1 | \leq \beta_1 + c.
\]
\end{thm}

The case $c=0$ can be found in \cite[Theorem 5.19]{Flo}. The proof of this generalization is very similar.
We refer to the upcoming paper \cite{Ore3} for an interesting application of this result to the study of real algebraic plane curves.

\begin{proof}

Let $W_F$ denote the exterior of $F$ in $B^4$. The non-vanishing homology groups of $W_F$ are given by
\[
H_0(W_F)=\Zz,\  H_1(W_F)=\Zz^\mu\ \hbox{ and }  \ H_2(W_F)=\Zz^{\beta_1+c}.
\]
Set $\beta_k^\omega= \dim H_k^\omega(W_F,\Cc)$. Clearly, $\beta_0^\omega=\beta_4^\omega=0$.
Since $\omega$ belongs to $T^\mu_\pp$, the order of the associated cyclic covering of $W_F$ is a power of prime.
Therefore, we can make use of Gilmer's results \cite{Gi1}. In particular, by \cite[Proposition 1.4]{Gi1}, $\beta^\omega_3=0$.
It follows that
\[
\beta_2^\omega-\beta_1^\omega = \chi(W_F) = 1-\mu+\beta_1+c.
\]
Therefore,
\[
|\sign(\varphi^\omega)|+\nuli(\varphi^\omega)\leq\beta_2^\omega=\beta_1^\omega+1-\mu+\beta_1+c.
\]
By Theorem \ref{main1}, $\sign(\varphi^\omega)=\sigma_L(\omega)$ and $\nuli(\varphi^\omega)+\beta_1^\omega=\eta_L(\omega)$,
leading to
\[
|\sigma_L(\omega)|+\eta_L(\omega)\leq 2\beta_1^\omega+1-\mu+\beta_1+c.
\]
By \cite[Proposition 1.5]{Gi1}, $\beta_1^\omega \leq \mu -1$, so $|\sigma_L(\omega)|+\eta_L(\omega)-\mu+1\leq \beta_1+c$,
giving the first part of the inequality. On the other hand, $\beta_1^\omega=\eta_L(\omega)-\nuli\varphi^\omega\le \eta_L(\omega)$, so
$|\sigma_L(\omega)|-\eta_L(\omega)+\mu-1\leq \beta_1+c$, completing the proof of the theorem.
\end{proof}

\subsection{The slice genus}\label{sub:KT}

\begin{definition}
Let $S^3$ denote the standard embedding of the $3$-sphere in $S^4$. The {\it slice genus} $g_s(L)$ of a $\mu$-colored link $L$ is
the minimal genus of a closed oriented smooth surface $P=P_1\sqcup\dots\sqcup P_\mu\subset S^4$ such that $P_i\cap S^3=L_i$ for all $i$.
A $\mu$-colored link is said to be {\it slice} if its slice genus is zero.
\end{definition}

Note that for such a surface to exist, we must have $lk(L_i,L_j)=0$ for all $i\neq j$.
This definition should be understood as a unification of several well-known notions of `sliceness.'
Indeed, consider the case $\mu=1$. A $1$-colored link is slice if it is the cross-section
of a smooth $2$-sphere in $S^4$, that is, using Fox's terminology \cite{Fo2}, if it is {\em slice in the ordinary sense\/}.
On the other hand, consider a $\nu$-colored link with $\nu$ components. Such a colored link is slice if it is the cross-section
of $\nu$ smooth disjoint $2$-spheres in $S^4$. According to Fox, such a link is called {\em slice in the strong sense\/}.

The signature and nullity provide a lower bound for the slice genus of a 
colored link. Indeed, we have the following generalization of \cite[Theorem 3.13]{Kauff-Tay}.

\begin{thm}
For all $\omega$ in $T^{\mu}_\pp$,
\[
|\sigma_L(\omega)|\leq g_s(L)+\min(0,\eta_L(\omega)+1-\mu).
\]
\end{thm}

\begin{proof}
Consider a closed oriented smooth surface $P=P_1\sqcup\dots\sqcup P_\mu$ in $S^4$ such that $P_i\cap S^3=L_i$ for all $i$. 
It may be assumed that each $P_i$ is connected. Let $W$ be the exterior of $P$ in $S^4$. By duality, the homology of $W$ is given by
\[
H_0(W)= \Zz,\ H_1(W)=\Zz^\mu,\ H_2(W)=\Zz^{2g} \hbox{ and } H_3(W)=\Zz^\mu,
\]
where $g$ denotes the genus of $P$. 
As in the previous section, any $\omega$ in $T^{\mu}_\pp$ induces a character $\chi_\omega\colon H_1(W)\to\Cc^*$ of prime power order
sending the meridian of $P_i$ to $\omega_i$. For simplicity, we simply write $H^\omega_*(W)$ for $H^\omega_*(W;\Cc)$.

Let $X$ be the exterior of $L$ in $S^3$. The sphere $S^3$ standardly embedded in $S^4$ splits $P$ into two
surfaces $F_1$ and $F_2$ with $F_1\cap F_2=L$. The manifold $W$ can be described as a union $W_1\cup W_2$ with $W_1\cap W_2=X$,
where $W_i$ is the complement of an open tubular neighborood of $F_i$ in $B^4$ for $i=1,2$.
The character $\chi_\omega$ restricts to characters on $H_1(W_i)$.
Let $\varphi_i^\omega$ be the intersection form on the corresponding twisted homology $H^\omega_2(W_i)$.
By Theorem \ref{main1}, $\sigma_L(\omega)=\sign(\varphi_i^\omega)$ for $i=1,2$. Clearly,
$\varphi_i^\omega$ is dual to the inclusion homomorphism $j_i\colon H_2^\omega(W_i)\to H_2^\omega(W_i,M_i)$,
where $M_i$ stands for $\partial W_i$. It follows that $|\sigma_L(\omega)|\leq\dim K_i$, where $K_i= H_2^\omega(W_i)/\ker(j_i)$. 

By the Mayer-Vietoris exact sequence with twisted coefficients, we have an isomorphism
\[
H_2^\omega(W_1,M_1) \oplus H_2^\omega(W_2,M_2)\simeq H^\omega_2(W,M_1 \cup M_2)
\]
which fits into the following commutative diagram, where the lines are exact:
\[
\xymatrix{
0 \ar[r] & K_1 \oplus K_2 \ar[r] & H_2^\omega(W_1,M_1)\oplus H_2^\omega(W_2,M_2)
\ar[d]^{\simeq}\ar[r]^{\qquad\partial_1\oplus\partial_2}&H^\omega_1(M_1)\oplus 
H^\omega_1(M_2) \ar[d]^\psi\\
& H_2^\omega(W) \ar[r] & H_2^\omega(W,M_1 \cup M_2)\ar[r]^\partial & H_1^\omega(M_1 \cup M_2).
}
\]
Therefore, $K_1\oplus K_2=Ker(\partial_1\oplus\partial_2)\subset Ker(\psi\circ(\partial_1\oplus\partial_2))\simeq Ker(\partial)$.
This implies that $\dim K_1+\dim K_2\leq\beta^\omega_2$, so $2|\sigma_L(\omega)|\leq\beta^\omega_2$. The equation
\[
1-2\mu+2g=\chi(W)=-\beta^\omega_1+\beta^\omega_2-\beta^\omega_3
\]
leads to
\[
2 | \sigma_L(\omega) | \leq 1 -2 \mu + 2g + \beta^\omega_1 + \beta^\omega_3.
\]
By \cite[Propositions 1.4 and 1.5]{Gi1} , $\beta^\omega_3\leq\dim H_3(W)=\mu$ and $\beta^\omega_1(W)\leq\dim H_1(W)-1=\mu-1$.
It follows that $|\sigma_L(\omega)|\leq g$ and the first part of the inequality is proved.

To check the second part, let us denote by $\mu_i$ be the number of 
components of $F_i$ and by $\beta_1(F_i)$ the rank of $H_1(F_i)$. Since $P= F_1 \cup F_2$ and $F_1\cap F_2=L$,
the additivity of the Euler characteristic implies
\[
\mu_1+\mu_2-\beta_1(F_1)-\beta_1(F_2)=2\mu-2g.
\]
We then apply Theorem \ref{thm:MT} to the surfaces $F_1$ and $F_2$, giving
\[
|\sigma_L(\omega)|-\eta_L(\omega)+\mu_i-1\leq\beta_1(F_i),
\]
for $i=1,2$. These three equations easily give
\[
|\sigma_L(\omega)|\leq g+\eta_L(\omega)+1-\mu,
\]
which implies the second part of the inequality.
\end{proof}

\begin{cor}
If a $\mu$-colored link $L$ is slice, then $\sigma_L(\omega)=0$ and $\eta_L(\omega)\ge\mu-1$ for all $\omega$ in $T^{\mu}_\pp$.\qed
\end{cor}

The case $\mu=1$ gives the following result.

\begin{cor}\label{cor:ordinary}
If a link is slice in the ordinary sense, then $\sigma_L(\omega)=0$ for all
$\omega$ root of the unity of prime power order.\qed
\end{cor}

On the other hand, if $\mu$ is the maximum number of colors, we get:

\begin{cor}\label{cor:strong}
If a link with $\nu$ components is slice in the strong sense, then
$\sigma_L(\omega_1,\dots,\omega_\nu)=0$ and $\eta_L(\omega_1,\dots,\omega_\nu)\ge\nu-1$
for all $(\omega_1,\dots,\omega_\nu)\in T^\nu_\pp$.\qed
\end{cor}

Finally, let $\Sigma_L$ denote the zero in $T_*^\mu$ of the first non-vanishing Alexander ideal of $L$.
By Theorem \ref{thm:semi}, the signature
and nullity of $L$ are continuous functions on $T^\mu_*\setminus\Sigma_L$. By density of $T^\mu_\pp$,
all the results stated in this section hold for $\omega$ in $T^\mu_*\setminus\Sigma_L$.
\medskip

Let us conclude this paper with one last didactic example.
\begin{ex}
In \cite{Fo2}, Fox presents the link illustrated below.

\noindent\parbox{0.65\textwidth}{\noindent It is a very simple link which is slice in the
ordinary sense, but not in the strong sense. We shall compute the signatures of this link in order to test the results of this section.
Let us order its components as illustrated to obtain a $3$-colored link $L$. There
is an obvious C-complex for $L$ which has the homotopy type of a circle. The corresponding Seifert matrices are given by $A^\eps=(-1)$
if $\pm\eps=(1,-1,1)$, and $A^\eps=(0)$ else.}\hfill\parbox{0.3\textwidth}{\epsfig{figure=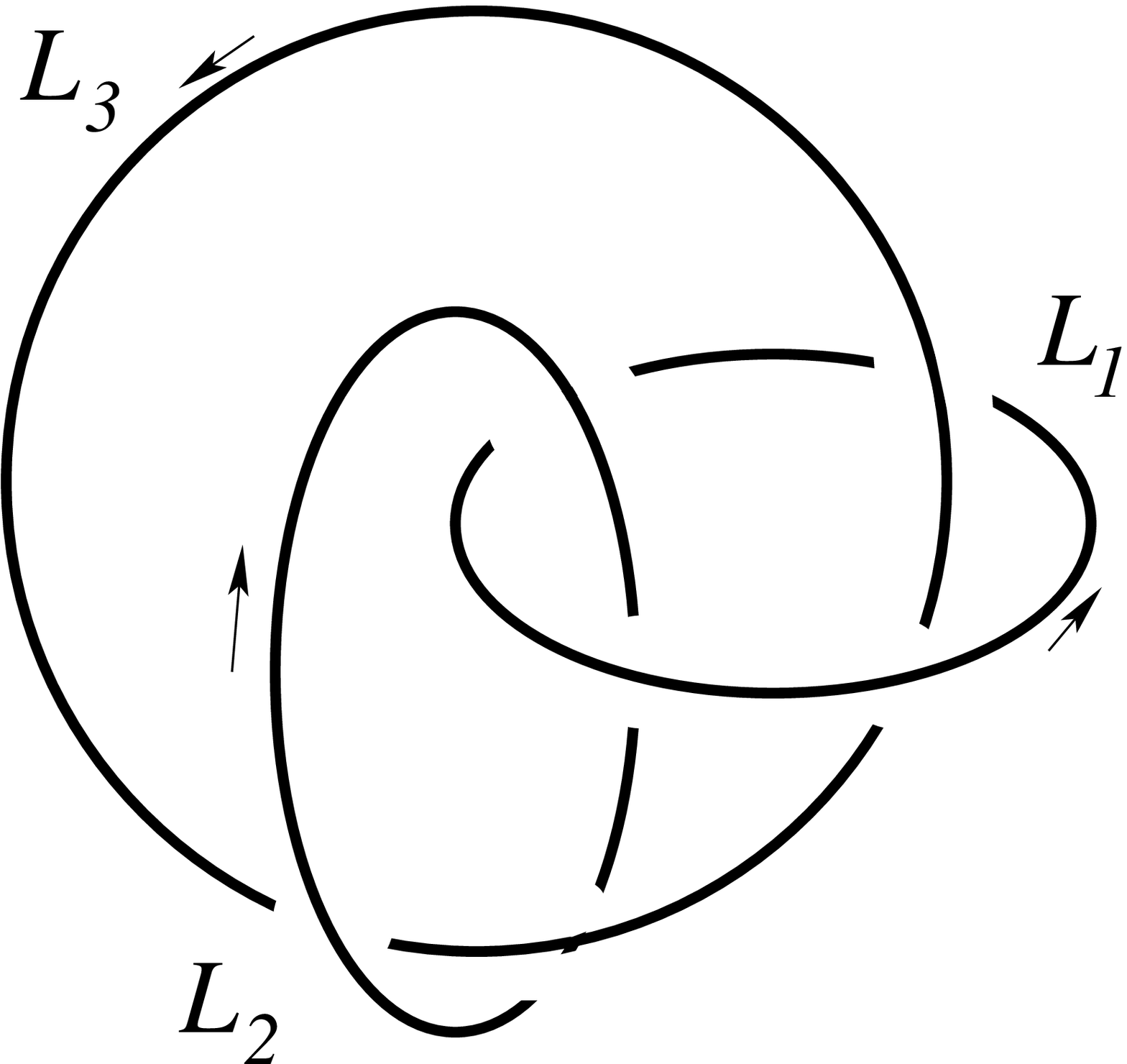,height=3cm}}
Therefore, we obtain
\[
\sigma_L(\omega_1,\omega_2,\omega_3)=sgn\,\Re((\omega_1-1)(\overline\omega_2-1)(\omega_3-1)),
\]
where $sgn$ denotes the sign function.
By Corollary \ref{cor:strong}, the $3$-colored link $L$ is not slice in the strong sense (which is obvious since the linking numbers
don't vanish). On the other hand, Proposition \ref{prop:lk} implies that the Levine-Tristram signature of the underlying
link $L'$ is equal to
\[
\sigma_{L'}(\omega)=sgn\,\Re(\omega-1)+1=0.
\]
This is the expected result since $L'$ is slice in the ordinary sense.
\end{ex}

\bibliographystyle{amsplain}

\end{document}